\documentclass[10pt]{amsart}
\usepackage{amsmath}
\usepackage{amssymb}
\usepackage{doublespace}
\usepackage[final]{graphicx}
\usepackage{float}
\usepackage{subfigure}
\usepackage{warmread}

\usepackage[authoryear]{natbib}

\usepackage[all,frame,ps,tips,line,poly,arrow,curve,import]{xy}
\SelectTips{cm}{10} \UseTips

\newcommand{\scaledfig}[2]{\scalebox{#1}{\includegraphics{#2}}}

\newcommand{\Div}{\operatorname{div}}
\newcommand{\curl}{\operatorname{curl}}

\def\ZZ{\mathbb{Z}}
\def\RR{\mathbb{R}}

\def\d{\mathbf{d}}
\def\X{\mathbf{X}}

\newcommand{\bfi}{\bfseries\itshape}

\DeclareMathOperator{\SO}{SO}

\restylefloat{figure}

\newtheorem{theorem}{Theorem}[section]
\newtheorem{corollary}[theorem]{Corollary}

\newtheorem{lemma}[theorem]{Lemma}

\newtheorem{remark}{Remark}[section]
\newtheorem{example}{Example}[section]

\newtheorem{definition}{Definition}[section]
\numberwithin{equation}{section}

\renewcommand{\paragraph}[1]{\vspace*{0.1in}\noindent\textbf{#1}}

\title{Discrete Exterior Calculus}
\author{Mathieu Desbrun}
\address{158-79, Computer Science, Caltech, Pasadena, CA 91125.}
\email{mathieu@caltech.edu}
\author{Anil N. Hirani}
\address{Department of Computer Science, University of Illinois, Urbana, IL 61801.}
\email{hirani@cs.uiuc.edu}
\author{Melvin Leok}
\address{Department of Mathematics, University of Michigan, Ann Arbor, MI 48109.}
\email{mleok@umich.edu}
\author{Jerrold E. Marsden}
\address{107-81, Control and Dynamical Systems, Caltech, Pasadena, CA 91125.}
\email{marsden@cds.caltech.edu}

\allowdisplaybreaks

\begin{document}

\begin{abstract}
We present a theory and applications of discrete exterior calculus
on simplicial complexes of arbitrary finite dimension. This can be
thought of as calculus on a discrete space. Our theory includes
not only discrete differential forms but also discrete vector
fields and the operators acting on these objects.  This allows us
to address the various interactions between forms and vector
fields (such as Lie derivatives) which are important in
applications. Previous attempts at discrete exterior calculus have
addressed only differential forms.  We also introduce the notion
of a circumcentric dual of a simplicial complex. The importance of
dual complexes in this field has been well understood, but
previous researchers have used barycentric subdivision or
barycentric duals. We show that the use of circumcentric duals is
crucial in arriving at a theory of discrete exterior calculus that
admits both vector fields and forms.
\end{abstract}

\maketitle

\setcounter{tocdepth}{2} \tableofcontents

\section{Introduction}
\label{dec:sec:intro}

This work presents a theory of {\bfi discrete exterior calculus} (DEC) motivated by potential applications in computational methods for field theories such as elasticity, fluids, and electromagnetism. In addition, it provides much needed mathematical machinery to enable a systematic development of numerical schemes that mirror the approach of geometric mechanics.

This theory has a long history that we shall outline below in \S\ref{dec:sec:Previous}, but we aim
at a comprehensive, systematic, as well as useful, treatment. Many previous works, as we shall review, are incomplete both in terms of the objects that they treat as well as the types of meshes that they allow.

Our vision of this theory is that it should proceed {\it ab initio} as a discrete theory that parallels the continuous one. General views of the subject area of DEC are common in the literature (see, for instance, \cite{Ma2000}), but they usually stress the process of discretizing a continuous theory and the overall approach is tied to this goal. However, if one takes the point of view that the discrete theory can, and indeed should, stand in its own right, then the range of application areas naturally is enriched and increases.

Convergence and consistency considerations alone are inadequate to discriminate between the various choices of discretization available to the numerical analyst, and only by requiring, when appropriate, that the discretization exhibits discrete analogues of continuous properties of interest can we begin to address the question of what makes a discrete theory a canonical discretization of a continuous one.

\paragraph{Applications to Variational Problems.}\index{variational problems} One of the major application areas we envision is to variational problems, be they in mechanics or optimal control. One of the key ingredients in this
direction that we imagine will play a key role in the future is that of AVI's (asynchronous variational integrators) designed for the numerical integration of mechanical systems, as in \cite{LeMaOrWe2003}. These are integration algorithms  that respect some of the key features of the continuous theory, such as their multi-symplectic nature and exact conservation laws. They do so by discretizing the underlying variational principles of mechanics rather than discretizing the equations. It is well-known (see the reference just mentioned for some of the literature) that variational problems come equipped with a rich exterior calculus structure and so on the discrete level, such structures will be enhanced by the availability of a discrete exterior calculus. One of the objectives of this chapter is to fill this gap.

\paragraph{Structured Constraints.}\index{constraints!structured} There are many constraints in numerical algorithms that naturally involve differential forms, such as the divergence constraint for incompressibility of fluids, as well as the fact that differential forms are naturally the fields in electromagnetism, and some of Maxwell's equations are expressed in terms of the divergence and curl operations on these fields. Preserving, as in the mimetic differencing literature, such features directly on the discrete level is another one of the goals, overlapping with our goals for variational problems.

\paragraph{Lattice Theories.}\index{lattice theory} Periodic crystalline lattices are of important practical interest in material science, and the anisotropic nature of the material properties arises from the geometry and connectivity of the intermolecular bonds in the lattice. It is natural to model these lattices as inherently discrete objects, and an understanding of discrete curvature that arises from DEC is particularly relevant, since part of the potential energy arises from stretched bonds that can be associated with discrete curvature in the underlying relaxed configuration of the lattice. In particular, this could yield a more detailed geometric understanding of what happens at grain boundaries. Lattice defects can also be associated with discrete curvature when appropriately interpreted. The introduction of a discrete notion of curvature will lay the foundations for a better understanding of the role of geometry in the material properties of solids.

\paragraph{Some of the Key Theoretical Accomplishments.} Our development of discrete exterior calculus includes discrete differential forms, the Hodge star operator, the wedge product, the
exterior derivative, as well as contraction and the Lie derivative. For example, this approach leads to the proper definition of discrete divergence and curl operators and has already resulted in applications like a discrete Hodge type decomposition of 3D vector fields on irregular grids---see \cite{ToLoHiDe2003}.

\paragraph{Context.} We present the theory and some applications of DEC in the context of simplicial complexes of arbitrary finite dimension.

\paragraph{Methodology.} We believe that the correct way to proceed with this program is to develop, as we have already stressed, {\it ab initio}, a calculus on discrete manifolds which parallels the calculus on smooth manifolds of arbitrary finite dimension. Chapters 6 and 7 of \cite{AbMaRa1988} are a good source for the concepts and definitions in the smooth case. However we have tried to make this chapter as self-contained as possible. Indeed, one advantage of developing a calculus on discrete manifolds, as we do here, is pedagogical. By using concrete examples of discrete two- and three-dimensional spaces one can explain most of calculus on manifolds at least formally as we will do using the examples in this chapter. The machinery of Riemannian manifolds and general manifold theory from the smooth case is, strictly speaking, not required in the discrete world. The technical terms that are used in this introduction will be defined in subsequent sections, but they should be already familiar to someone who knows the usual exterior calculus on smooth manifolds.

\paragraph{The Objects in DEC.} To develop a discrete theory, one must define discrete differential forms along with vector fields and operators involving these. Once discrete forms and vector fields are
defined, a calculus can be developed by defining the discrete
exterior derivative
($\mathbf{d}$\index{$d$@$\mathbf{d}$|see{exterior derivative}}),
codifferential
($\boldsymbol{\delta}$\index{$\delta$@$\boldsymbol{\delta}$|see{codifferential}})
and Hodge star ($\ast$\index{$\ast$|see{Hodge star}}) for
operating on forms, discrete wedge product
($\wedge$\index{$\wedge$|see{wedge product}}) for combining forms,
discrete flat ($\flat$) and sharp ($\sharp$) operators for going
between vector fields and $1$-forms and discrete contraction
operator
($\mathbf{i}_X$\index{$i_x$@$\mathbf{i}_X$|see{contraction}}) for
combining forms and vector fields. Once these are done, one can
then define other useful operators. For example, a discrete Lie
derivative ($\pounds_X$) can be {\em defined} by requiring that
the Cartan magic (or homotopy) formula hold. A discrete divergence
in any dimension can be defined. A discrete Laplace--deRham
operator ($\Delta$) can be defined using the usual definition of
$\mathbf{d} \boldsymbol{\delta} + \boldsymbol{\delta} \mathbf{d}$.
When applied to functions, this is the same as the discrete
Laplace--Beltrami operator ($\nabla^2$), which is the defined as
$\operatorname{div} \circ \operatorname{curl}$. We define all
these operators in this chapter.

The discrete manifolds we work with are simplicial complexes. We will recall the standard formal definitions in \S\ref{dec:sec:Primal} but familiar examples of simplicial complexes are meshes of triangles embedded in $\mathbb{R}^3$ and meshes made up of tetrahedra occupying a portion of $\mathbb{R}^3$. We will assume that the angles and lengths on such discrete manifolds are computed in the embedding space $\mathbb{R}^N$ using the standard metric of that space. In other words, in this chapter we do not address the issue of how to discretize a given smooth Riemannian manifold, and how to embed it in $\mathbb{R}^N$, since there may be many ways to do this. For example, $\SO(3)$ can be embedded in $\mathbb{R}^9$ with a constraint, or as the unit quaternions in $\mathbb{R}^4$. Another potentially important consideration in discretizing the manifold is that the topology of the simplicial complex should be the same as the manifold to be discretized. This can be verified using the methods of computational homology (see, for example, \cite{KaMiMr2004}), or discrete Morse theory (see, for example, \cite{Fo2002, Wood2003}). For the purposes of discrete exterior calculus, only local metric information is required, and we will comment towards the end of \S\ref{dec:sec:Primal} how to address the issue of embedding in a local fashion, as well as the criterion for a good global embedding.

Our development in this chapter is for the most part formal in that we choose appropriate geometric definitions of the various objects and quantities involved.  For the most part, we do not prove that these definitions converge to the smooth counterparts. The definitions are chosen so as to make some important theorems like the generalized Stokes' theorem true by definition.  Moreover, in the cases where previous results are available, we have checked that the operators we obtain match the ones obtained by other means, such as variational derivations.

\section{History and Previous Work}
\label{dec:sec:Previous}

The use of simplicial chains and cochains as the basic building
blocks for a discrete exterior calculus has appeared in several
papers. See, for instance, \cite{SeSeSeAd2000}, \cite{Ad1996},
\cite{Bo2002c}, and references therein. These authors view forms as
linearly interpolated versions of smooth differential forms, a
viewpoint originating from \cite{Wh1957}, who introduced the
Whitney and deRham maps that establish an isomorphism between
simplicial cochains and Lipschitz differential forms.

We will, however, view discrete forms as real-valued linear functions on the space of chains. These are inherently discrete objects that can be paired with chains of oriented simplices, or their geometric duals, by the bilinear pairing of evaluation. In the next chapter, where we consider applications involving the curvature of a discrete space, we will relax the condition that discrete forms are real-valued, and consider group-valued forms.

Intuitively, this natural pairing of evaluation can be thought of as integration of the discrete form over the chain. This difference from the work of \cite{SeSeSeAd2000} and \cite{Ad1996} is apparent in the definitions of operations like the wedge product as well.

There is also much interest in a discrete exterior calculus in the
computational electromagnetism community, as represented by
\cite{Bo2001, Bo2002a, Bo2002b, Bo2002c}, \cite{GrKo2001},
\cite{Hi1999, Hi2001a, Hi2001b, Hi2002}, \cite{Ma1997, Ma2000}, \cite{NiWa1998},
\cite{Te2001}, and \cite{To2002}.

Many of the authors cited above, for example, \cite{Bo2002c}, \cite{SeSeSeAd2000}, and \cite{Hi2002}, also introduce the notions of dual complexes in order to construct the Hodge star operator. With the exception of Hiptmair, they use barycentric duals. This works if one develops a theory of discrete forms and does not introduce discrete vector fields.  We show later that to introduce discrete vector fields into the theory the notion of circumcentric duals seems to be important.

Other authors, such as \cite{Mo2000, MoSc2001, ScMo1999}, have incorporated vector fields into the cochain based approach to exterior calculus by identifying vector fields with cochains, and having them supported on the same mesh. This is ultimately an unsatisfactory approach, since dual meshes are essential as a means of encoding physically relevant phenomena such as fluxes across boundaries.

The use of primal and dual meshes arises most often as staggered meshes in finite volume and finite difference methods. In fluid computations, for example, the density is often a cell-centered quantity, which can either be represented as a primal object by being associated with the $3$-cell, or as a dual object associated with the $0$-cell at the center of the $3$-cell. Similarly, the flux across boundaries can be associated with the $2$-cells that make up the boundary, or the $1$-cell which is normal to the boundary.

Another approach to a discrete exterior calculus is presented in \cite{De1995}.  He defines a one-dimensional discretization of the real line in much the same way we would. However, to generalize to higher dimensions he introduces a tensor product of this space. This results in logically rectangular meshes. Our calculus, however, is defined over simplicial meshes. A further difference is that like other authors in this field, \cite{De1995} does not introduce vector fields into his theory.

A related effort for three-dimensional domains with logically rectangular meshes is that of \cite{MaHy2001}, who established a variational complex for difference equations by constructing a discrete homotopy operator. We construct an analogous homotopy operator for simplicial meshes in proving the discrete Poincar\'e lemma.

\section{Primal Simplicial Complex and Dual Cell Complex}
\label{dec:sec:Primal}
In constructing the discretization of a continuous problem in the context of our formulation of discrete exterior calculus, we first discretize the manifold of interest as a simplicial complex. While this is typically in the form of a simplicial complex that is embedded into Euclidean space, it is only necessary to have an abstract simplicial complex, along with a local metric defined on adjacent vertices. This abstract setting will be addressed further toward the end of this section.

We will now recall some basic definitions of simplices and simplicial complexes, which are standard from simplicial algebraic
topology. A more extensive treatment can be found in \cite{Mu1984}.

\begin{definition}
A $k$-\textbf{simplex}\index{simplex} is the convex span of $k+1$
geometrically independent points,
\begin{align*}
\sigma^k=[ v_{0},v_{1},\ldots,v_k] =\left\{
\sum_{i=0}^k\alpha^{i}v_{i}\,\Bigg|\,\alpha^{i}\geq0,\sum_{i=0}^{n}\alpha
^{i}=1\right\}.
\end{align*}
The points $v_0,\ldots, v_k$ are called the \textbf{vertices}\index{vertex} of
the simplex, and the number $k$ is called the \textbf{dimension}
of the simplex. Any simplex spanned by a (proper) subset of
$\{v_0,\ldots,v_k\}$ is called a \textbf{(proper) face}\index{face} of
$\sigma^k$. If $\sigma^l$ is a proper face of $\sigma^k$, we denote this by $\sigma^l \prec \sigma^k$.
\end{definition}

\begin{example}\index{simplex!examples}
Consider 3 non-collinear points $v_0, v_1$ and $v_2$ in $\RR^3$.
Then, these three points individually are examples of $0$-simplices, to which an orientation is assigned through the choice of a sign. Examples of $1$-simplices
are the oriented line segments $[v_0, v_1]$, $[v_1, v_2]$ and
$[v_0, v_2]$. By writing the vertices in that order we have given
orientations to these $1$-simplices, i.e., $[v_0, v_1]$ is oriented
from $v_0$ to $v_1$. The triangle $[v_0, v_1, v_2]$ is a $2$-simplex
oriented in counterclockwise direction. Note that the orientation
of $[v_0, v_2]$ does not agree with that of the triangle.
\end{example}

\begin{definition}
A \textbf{simplicial complex}\index{simplicial!complex}\index{complex!simplicial|see{simplicial, complex}} $K$ in
$\mathbb{R}^N$ is a collection of simplices in $\mathbb{R}^N$, such
that,
\begin{enumerate}
\item Every face of a simplex of $K$ is in $K$.
\item The intersection of any two simplices of $K$ is a face of
each of them.
\end{enumerate}
\end{definition}

\begin{definition}
A \textbf{simplicial triangulation}\index{simplicial!triangulation} of a polytope $|K|$\index{polytope} is a simplicial complex $K$ such that the union of the simplices of $K$ recovers the polytope $|K|$.
\end{definition}

\begin{definition}
If $L$ is a subcollection of $K$ that contains all faces of its
elements, then $L$ is a simplicial complex in its own right, and
it is called a \textbf{subcomplex}\index{subcomplex}\index{complex!sub|see{subcomplex}} of $K$. One subcomplex of $K$
is the collection of all simplices of $K$ of dimension at most
$k$, which is called the \textbf{$k$-skeleton}\index{skeleton} of $K$, and is
denoted $K^{(k)}$.
\end{definition}

\paragraph{Circumcentric Subdivision.}
We will also use the notion of a circumcentric dual or Voronoi
mesh of the given primal mesh. We will point to the importance of
this choice later on in \S\ref{dec:sec:Maps} and
\ref{dec:sec:Divergence}. We call the Voronoi dual a circumcentric
dual since the dual of a simplex is its circumcenter (equidistant
from all vertices of the simplex).

\begin{definition}
The \textbf{circumcenter}\index{circumcenter} of a $k$-simplex
$\sigma^k$ is given by the center of the $k$-circumsphere, where
the $k$-circumsphere is the unique $k$-sphere that has all $k+1$
vertices of $\sigma^k$ on its surface.  Equivalently, the
circumcenter is the unique point in the $k$-dimensional affine
space that contains the $k$-simplex that is equidistant from all
the $k+1$ nodes of the simplex. We will denote the circumcenter of
a simplex $\sigma^k$ by $c(\sigma^k)$.
\end{definition}

The circumcenter of a simplex $\sigma^k$ can be obtained by taking
the intersection of the normals to the $(k-1)$-dimensional faces
of the simplex, where the normals are emanating from the
circumcenter of the face. This allows us to recursively compute
the circumcenter.

If we are given the nodes which describe the primal mesh, we can
construct a simplicial triangulation by using the Delaunay
triangulation, since this ensures that the circumcenter of a
simplex is always a point within the simplex. Otherwise we assume
that a nice mesh has been given to us, i.e., it is such that the
circumcenters lie within the simplices. While this is not be
essential for our theory it makes some proofs simpler.  For some
computations the Delaunay triangulation is desirable in that it
reduces the maximum aspect ratio of the mesh, which is a factor in
determining the rate at which the corresponding numerical scheme
converges. But in practice there are many problems for which
Delaunay triangulations are a bad idea. See, for example,
\cite{Sc2002}.  We will address such computational issues in a
separate work.

\begin{definition}
The \textbf{circumcentric subdivision}\index{circumcentric!
subdivision} of a simplicial complex is given by the collection of
all simplices of the form 
\[[ c(\sigma_0),\ldots, c(\sigma_k)],\]
where $\sigma_0\prec\sigma_1\prec\ldots\prec\sigma_k$, or
equivalently, that $\sigma_i$ is a proper face of $\sigma_j$ for
all $i<j$.
\end{definition}

\paragraph{Circumcentric Dual.}\index{circumcentric!dual}\index{dual cell!circumcentric|see{circumcentric, dual}}
We construct a circumcentric dual to a $k$-simplex using the
circumcentric duality operator, which is introduced below.

\begin{definition}
The \textbf{circumcentric duality operator}\index{circumcentric!duality operator}\index{$\star$|see{circumcentric, duality operator}} is given by
\begin{align*}
\star\left(\sigma^k\right)
=\sum_{\sigma^k\prec\sigma^{k+1}\prec\ldots\prec\sigma^n}\epsilon_{\sigma^k,\ldots,\sigma^n}\left[
c(\sigma^k), c(\sigma^{k+1}),\ldots,c(\sigma^n)  \right],
\end{align*}
where the $\epsilon_{\sigma^k,\ldots,\sigma^n}$ coefficient
ensures that the orientation of $\left[ c(\sigma^k),
c(\sigma^{k+1}),\ldots,c(\sigma^n)  \right]$ is consistent with
the orientation of the primal simplex, and the ambient volume-form.

Orienting $\sigma^k$ is equivalent to choosing a ordered basis,
which we shall denote by $dx^1\wedge\ldots\wedge dx^k$. Similarly,
$\left[ c(\sigma^k), c(\sigma^{k+1}),\ldots,c(\sigma^n)  \right]$
has an orientation denoted by $dx^{k+1}\wedge\ldots\wedge dx^n$.
If the orientation corresponding to $dx^1\wedge\ldots\wedge dx^n$
is consistent with the volume-form on the manifold, then
$\epsilon_{\sigma^k,\ldots,\sigma^n}=1$, otherwise it takes the
value $-1$.
\end{definition}

We immediately see from the construction of the circumcentric
duality operator that the dual elements can be realized as a
submesh of the first circumcentric subdivision, since it consists
of elements of the form $[c({\sigma}_0),\ldots, c({\sigma}_k)]$,
which are, by definition, part of the first circumcentric
subdivision.

\begin{example}
The circumcentric duality operator maps a $0$-simplex into the
convex hull generated by
the circumcenters of $n$-simplices that contain the $0$-simplex,%
\begin{align*}
\star(\sigma^0)  =\left\{\sum\alpha_{\sigma^n}c\left(
\sigma^n\right)\,\Big|\,\alpha_{\sigma^n}\geq
0,\sum\alpha_{\sigma^n}=1, \sigma^0\prec \sigma^n\right\} ,
\end{align*}
and the circumcentric duality operator maps a $n-$simplex into the
circumcenter of the
$n-$simplex,%
\begin{align*}
\star(\sigma^n)  =c(\sigma^n) .
\end{align*}
\end{example}

\begin{figure}[H]
\begin{center}
\subfigure[Primal]{\includegraphics[scale=0.3,clip=true]{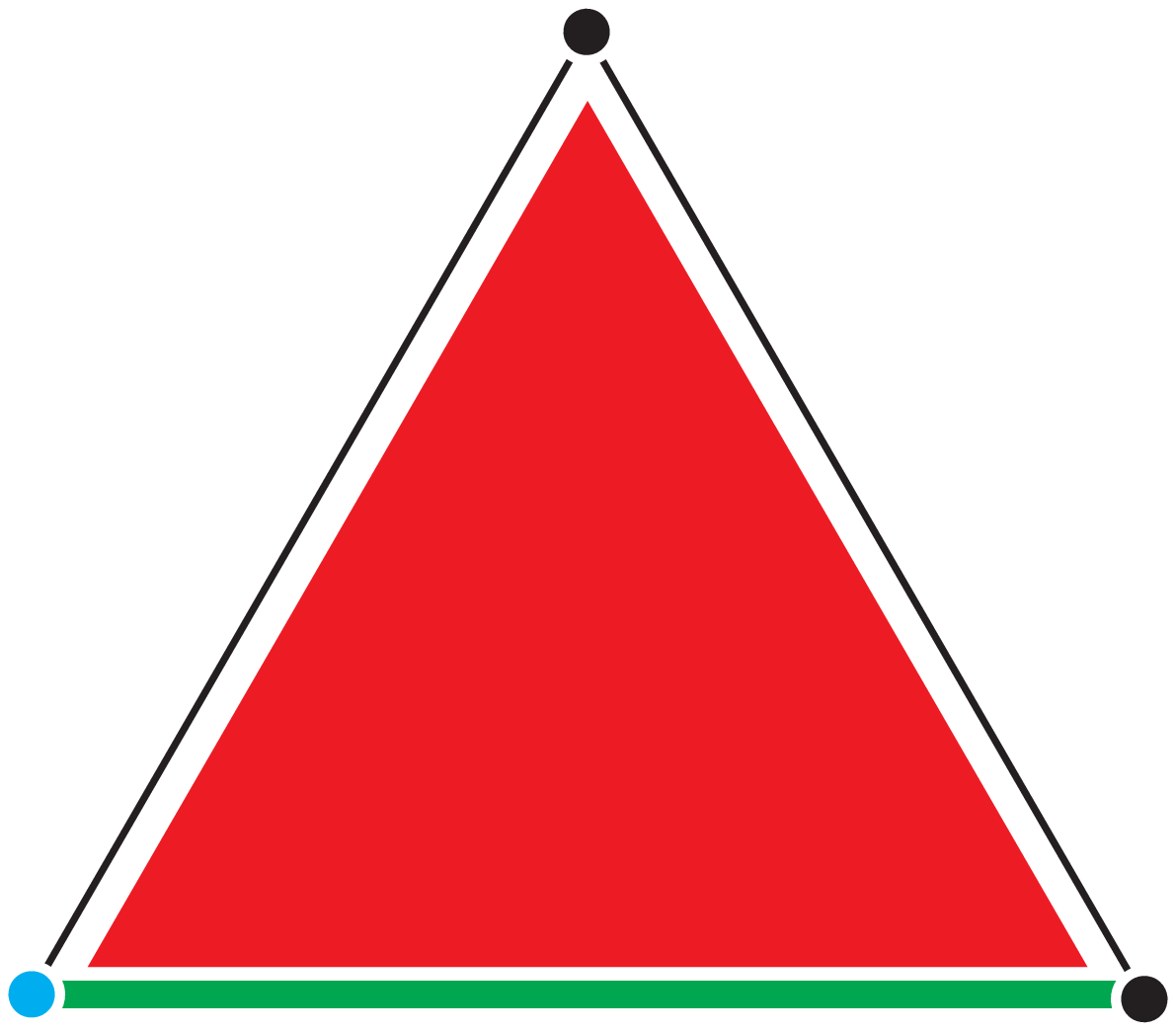}}\qquad
\subfigure[Dual]{\includegraphics[scale=0.3,clip=true]{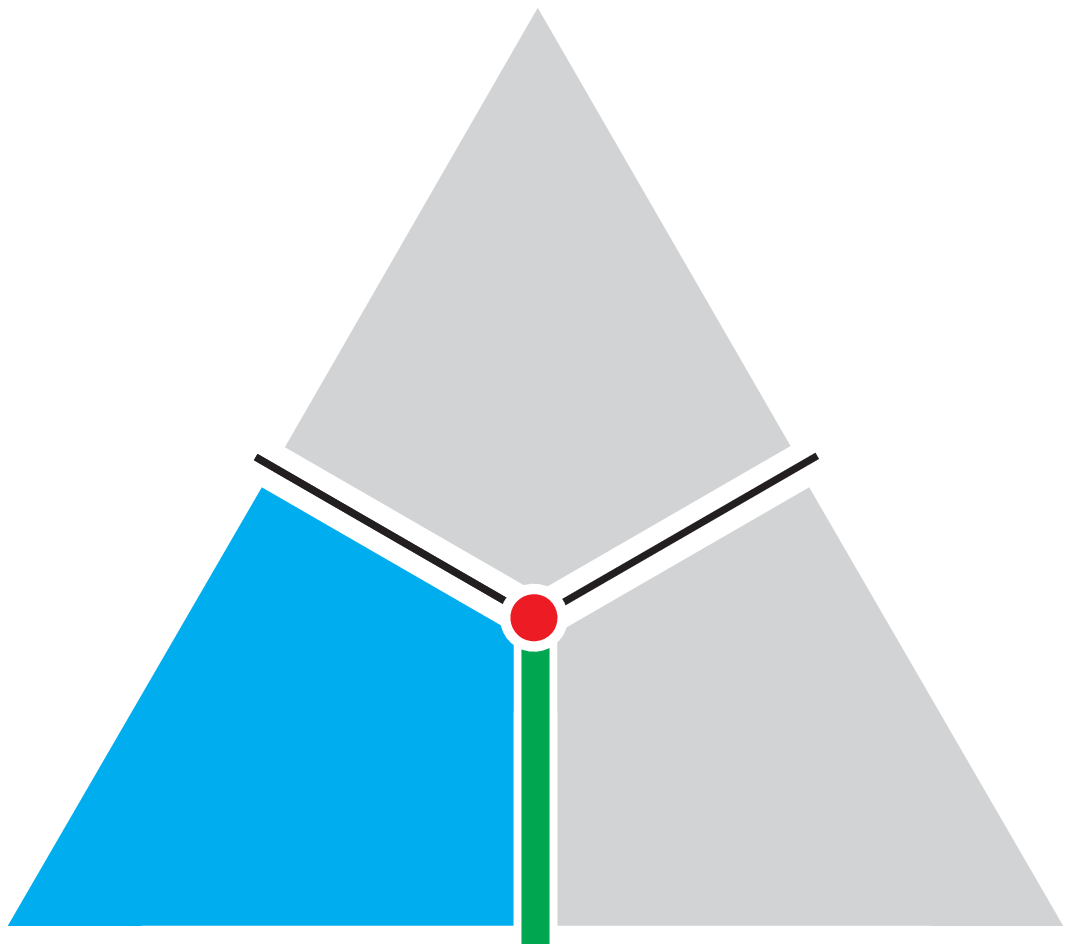}}\qquad
\subfigure[First subdivision]{\includegraphics[scale=0.3,clip=true]{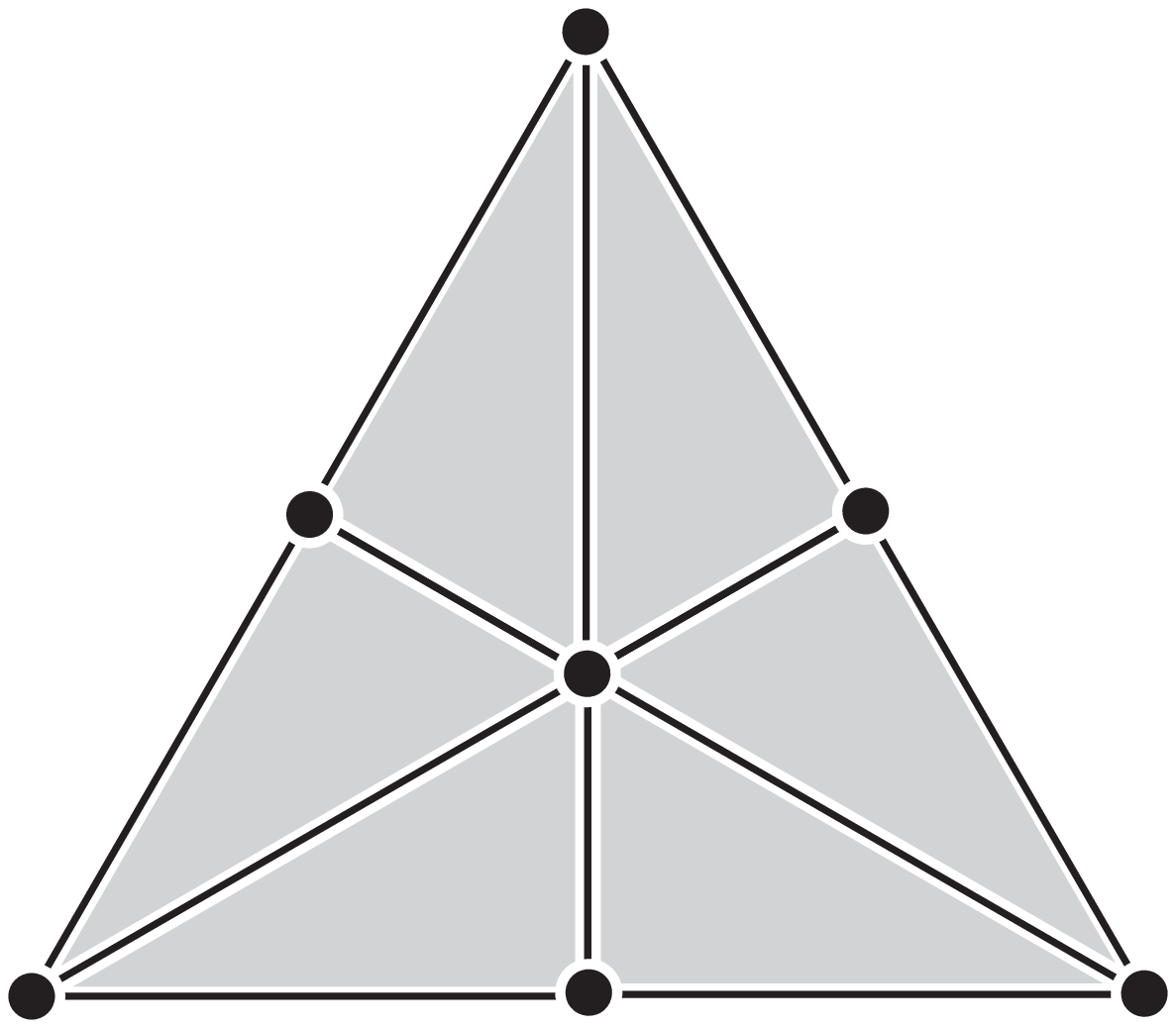}}
\end{center}
\caption{\label{dec:fig:primal_dual_refine}Primal, and dual meshes, as chains in the first circumcentric subdivision.}
\end{figure}

This is more clearly illustrated in Figure~\ref{dec:fig:primal_dual_refine}, where
the primal and dual elements are color coded to represent the dual
relationship between the elements in the primal and dual mesh.

The choice of a circumcentric dual is significant, since it allows
us to recover geometrically important objects such as normals to
$(n-1)$-dimensional faces, which are obtained by taking their
circumcentric dual, whereas, if we were to use a barycentric dual,
the dual to a $(n-1)$-dimensional face would not be normal to it.

\paragraph{Orientation of the Dual Cell.}\index{dual cell!orientation|see{orientation, dual cell}}\index{orientation!dual cell}
Notice that given an oriented simplex $\sigma^k$, which is represented by
$[v_0,\ldots, v_k]$, the orientation is equivalently represented
by $(v_1-v_0)\wedge (v_2-v_1) \wedge \ldots \wedge (v_k-v_{k-1})$,
which we denote by,
\[ [v_0,\ldots, v_k]\sim (v_1-v_0)\wedge (v_2-v_1) \wedge \ldots \wedge
(v_k-v_{k-1}), \] which is an equivalence at the level of orientation. It would be nice to express our criterion for determining the orientation of the dual cell in terms of the $(k+1)$-vertex representation.

To determine the orientation of the $(n-k)$-simplex given by
$[c(\sigma^k), c(\sigma^{k+1}),\ldots, c(\sigma^n)]$, or
equivalently, $dx^{k+1}\wedge\ldots\wedge dx^n$, we consider the
$n$-simplex given by $[c(\sigma^0),\ldots, c(\sigma^n)]$, where
$\sigma^0\prec\ldots\prec\sigma^k$. This is related to the
expression $dx^1\wedge\ldots\wedge dx^n$, up to a sign determined
by the relative orientation of $[c(\sigma^0),\ldots, c(\sigma^k)]$
and $\sigma^k$. Thus, we have that
\[
dx^1\wedge\ldots\wedge dx^n \sim
\operatorname{sgn}([c(\sigma^0),\ldots, c(\sigma^k)],\sigma^k)
[c(\sigma^0),\ldots, c(\sigma^n)]\, .
\]
Then, we need to check that $dx^1\wedge\ldots\wedge dx^n$ is
consistent with the volume-form on the manifold, which is
represented by the orientation of $\sigma^n$. Thus, we have that
the correct orientation for the $[c(\sigma^k),
c(\sigma^{k+1}),\ldots, c(\sigma^n)]$ term is given by,
\[
\operatorname{sgn}([c(\sigma^0),\ldots,
c(\sigma^k)],\sigma^k)\cdot\operatorname{sgn}([c(\sigma^0),\ldots,
c(\sigma^n)],\sigma^n).
\]
These two representations of the choice of orientation for the
dual cells are equivalent, but the combinatorial definition above
might be preferable for the purposes of implementation.

\begin{example}\index{orientation!example}
We would like to compute the orientation of the dual of a
$1$-simplex, in two dimensions, given the orientation of the two
neighboring $2$-simplices.

Given a simplicial complex, as shown in Figure~\ref{dec:subfig:dual:complex}, we consider a $2$-simplex of the form $[c(\sigma^0),c(\sigma^1),c(\sigma^2)]$, which is illustrated in Figure~\ref{dec:subfig:dual:2simplex}.
\begin{figure}[htbp]
\begin{center}
\subfigure[Simplicial complex]{\quad\qquad\includegraphics[scale=0.15]{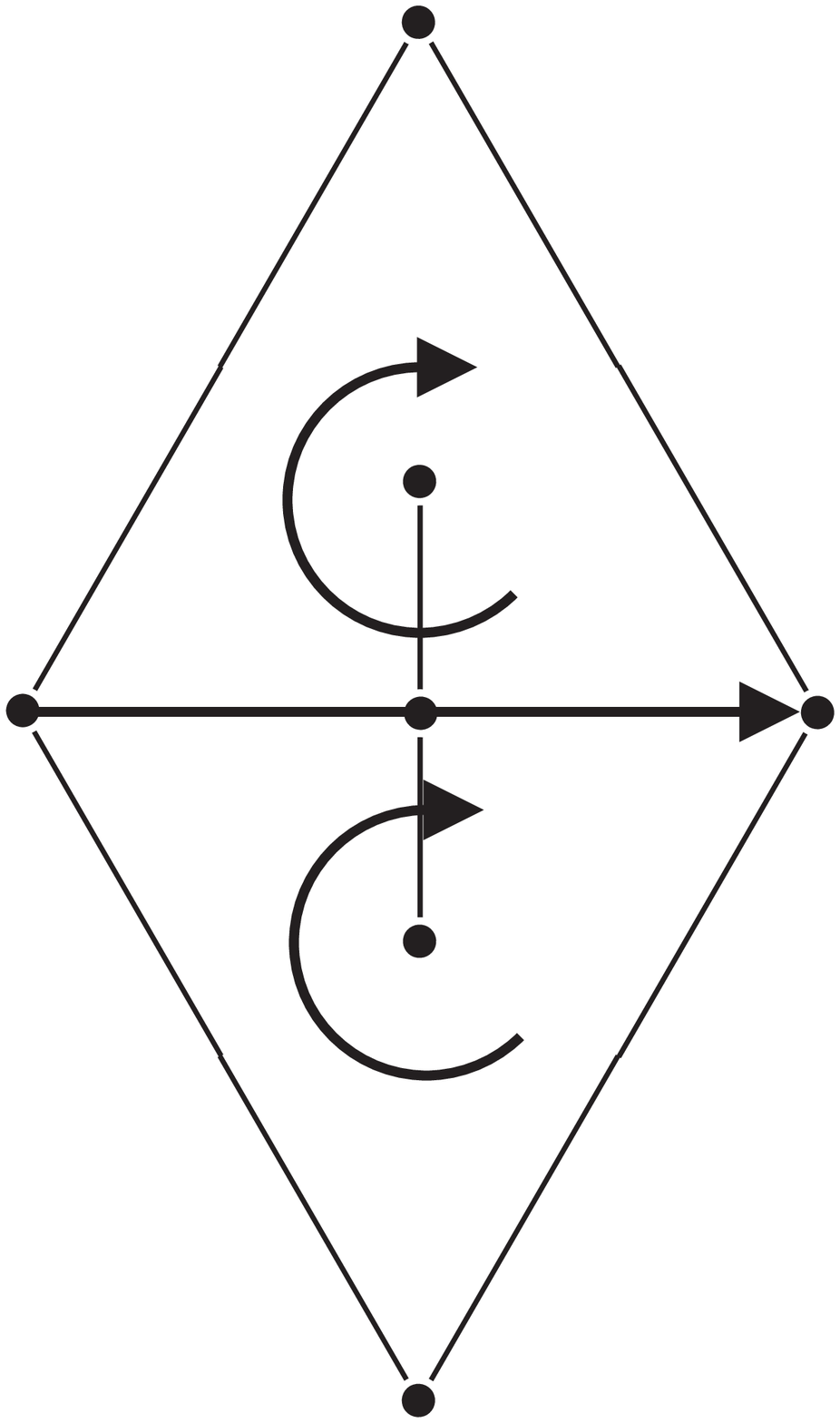}\label{dec:subfig:dual:complex}\quad\qquad}
\subfigure[$2$-simplex]{\quad\qquad\includegraphics[scale=0.15]{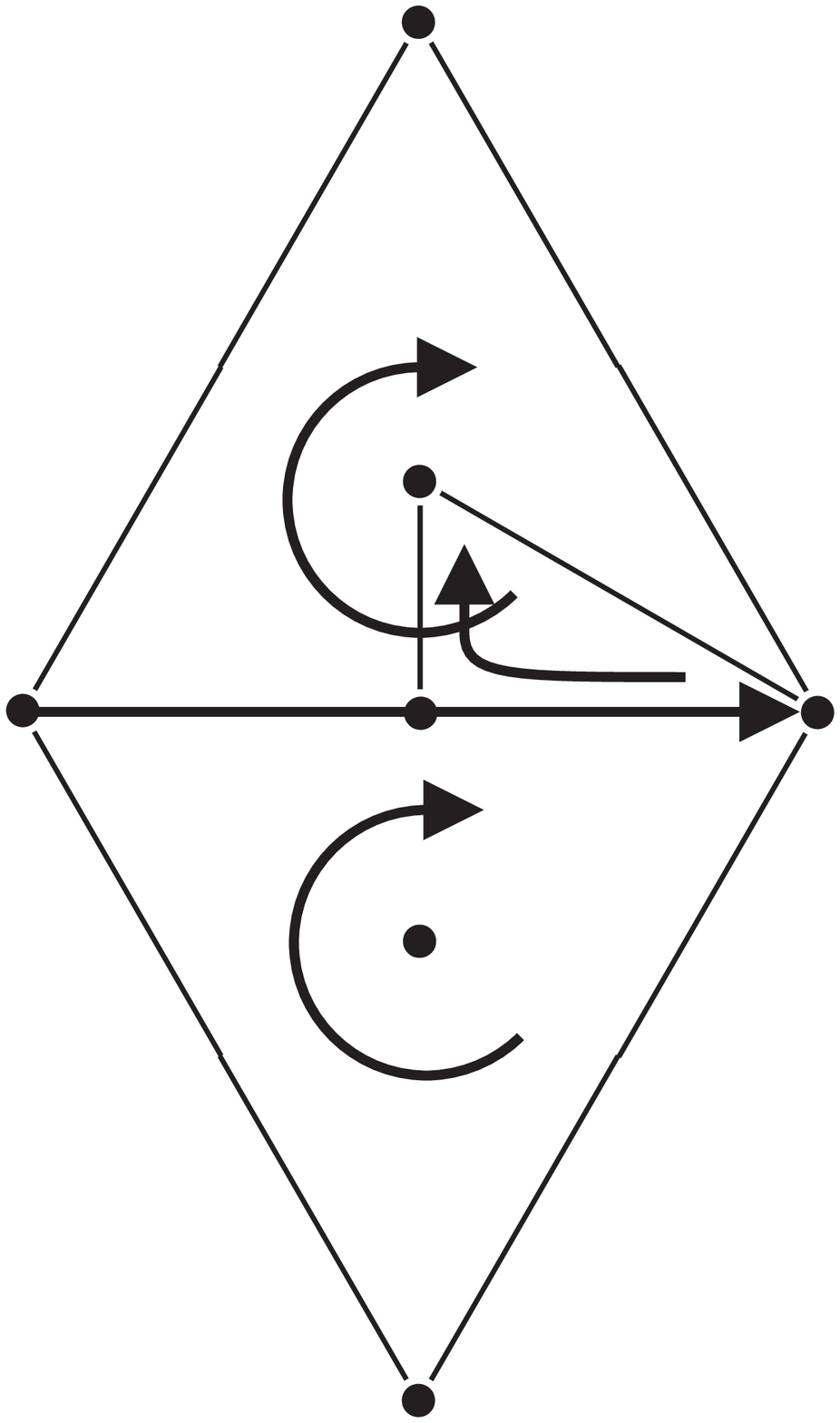}\label{dec:subfig:dual:2simplex}\quad\qquad}
\subfigure[$\star\sigma^1$]{\quad\qquad\includegraphics[scale=0.15]{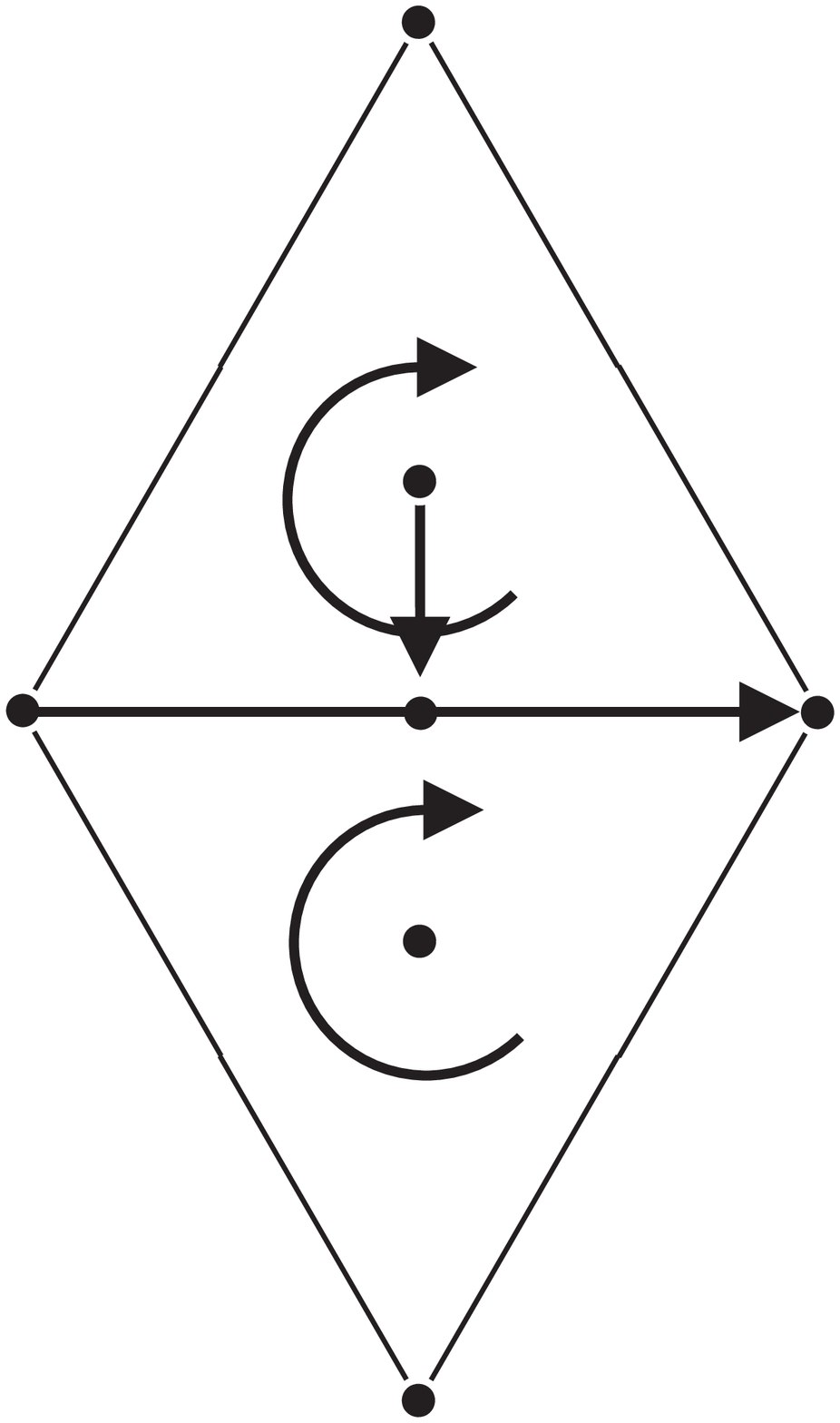}\label{dec:subfig:dual:dual}\quad\qquad}
\end{center}
\caption{\label{dec:fig:dual_orientation}Orienting the dual of a cell.}
\end{figure}

Notice that the orientation is consistent with the given orientation of the
$2$-simplex, but it is not consistent with the orientation of the
primal $1$-simplex, so the orientation should be reversed, to
give the dual cell illustrated in Figure~\ref{dec:subfig:dual:dual}.

We summarize the results for the induced orientation of dual cells for the other $2$-simplices of the form $[c(\sigma^0),c(\sigma^1),c(\sigma^2)]$,
in Table~\ref{dec:table:orientation}.
\begin{table}[h!]
\caption{\label{dec:table:orientation}\index{orientation!example}Determining the induced orientation of a dual cell.}
\begin{center}
\begin{tabular}{|c|c|c|c|c|}
  \hline
  \raisebox{8ex}[0pt]{$[c(\sigma^0),c(\sigma^1),c(\sigma^2)]$} &
  \includegraphics[scale=0.15]{orient1a} &
  \includegraphics[scale=0.15]{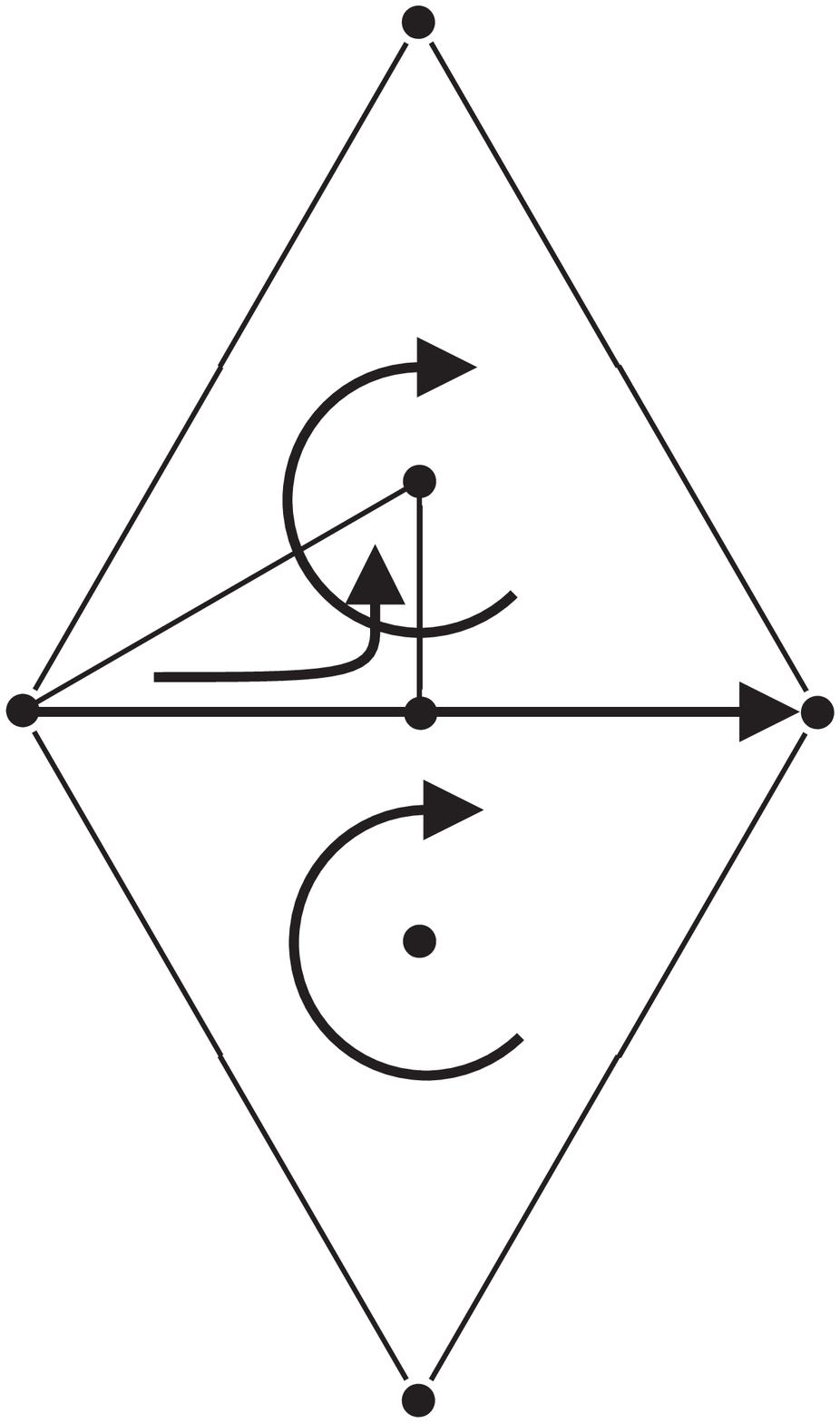} &
  \includegraphics[scale=0.15]{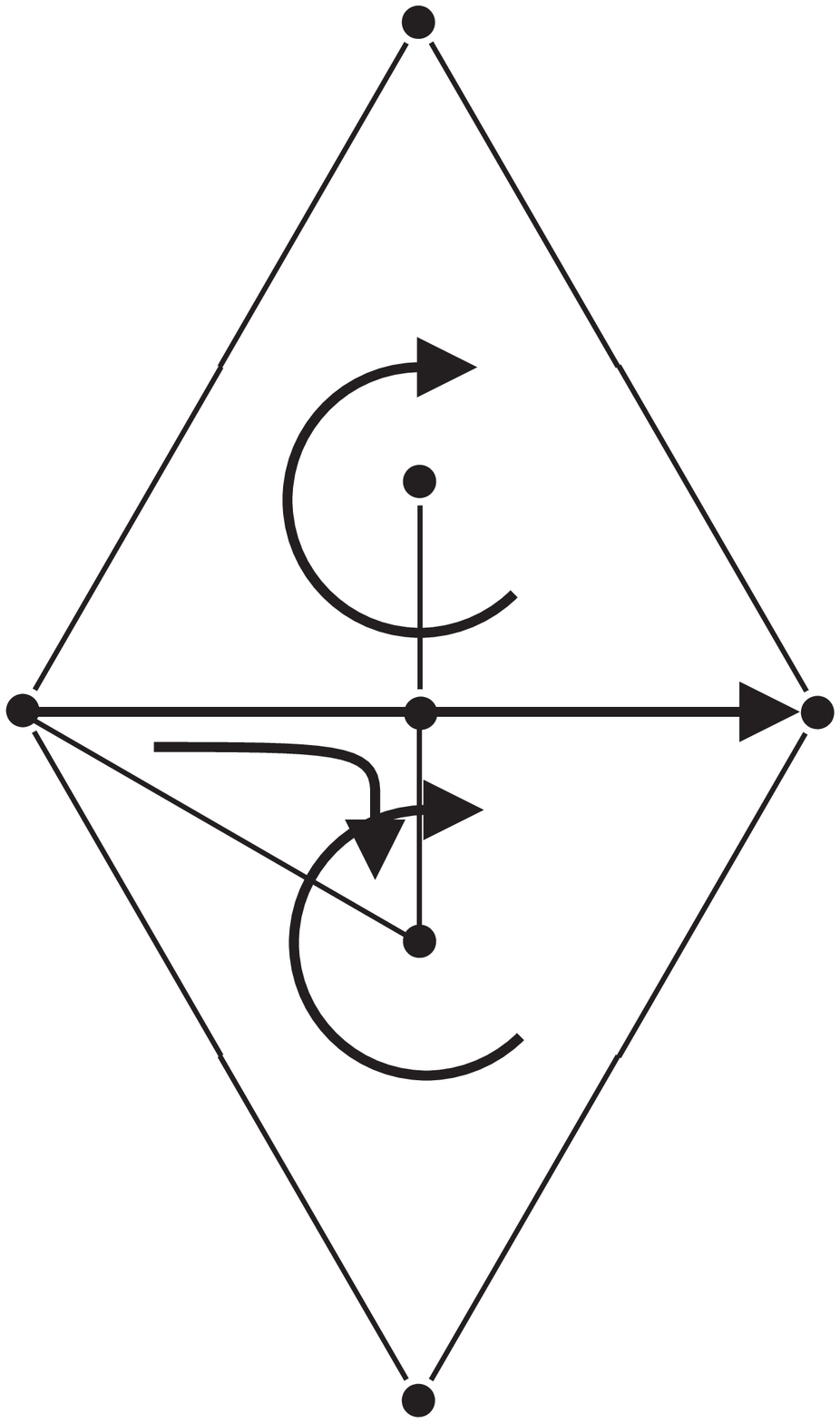} &
  \includegraphics[scale=0.15]{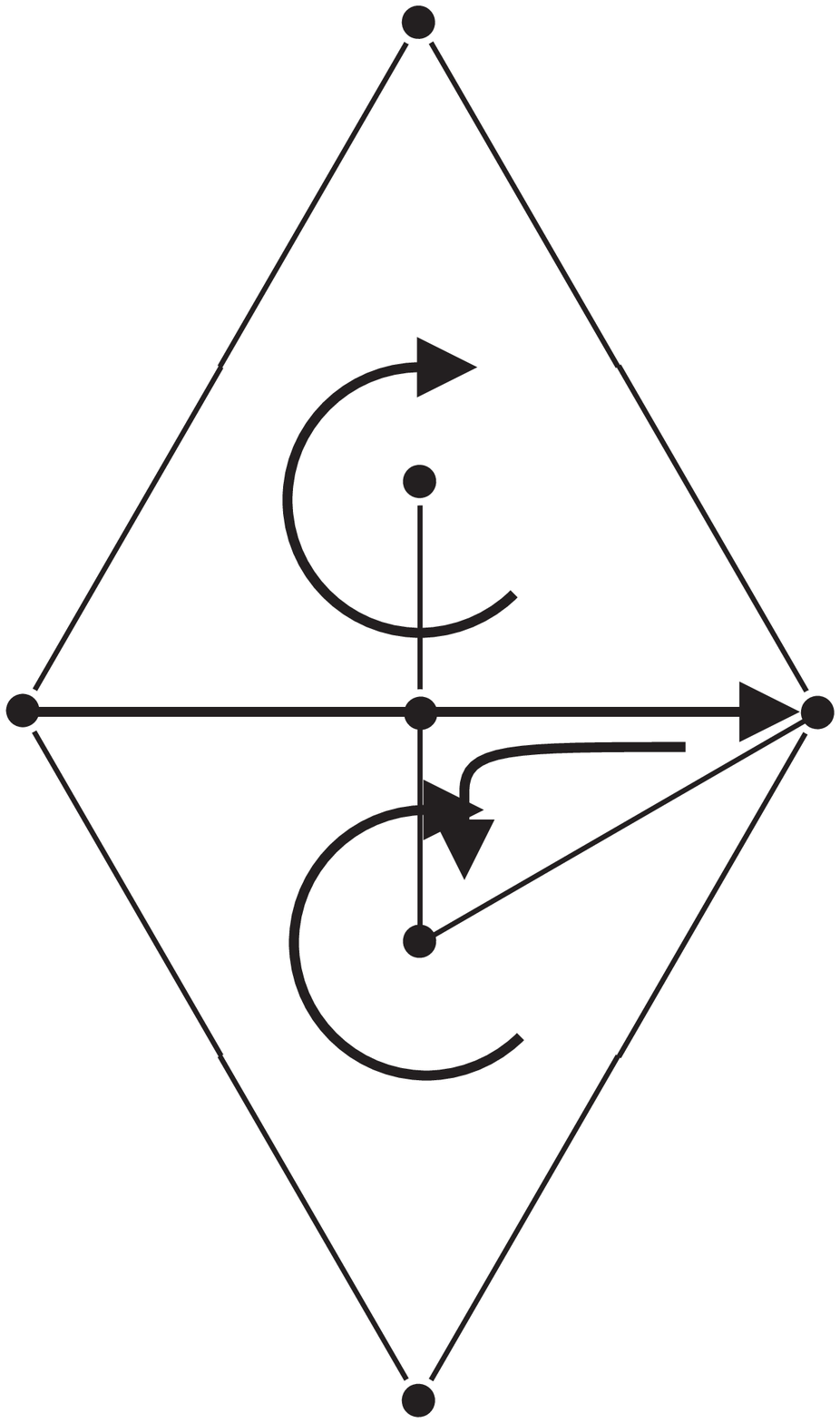}\\ \hline
  \rule[-3mm]{0mm}{8mm}
  $\operatorname{sgn}([c(\sigma^0),c(\sigma^1)],\sigma^1)$ &
  $-$ & $+$ & $+$ & $-$ \\ \hline
  \rule[-3mm]{0mm}{8mm}
  $\operatorname{sgn}([c(\sigma^0),c(\sigma^1),c(\sigma^2)],\sigma^2)$
  & $+$ & $-$ & $+$ & $-$ \\ \hline
  \raisebox{5ex}[0pt]{\parbox[b]{4.3cm}{$\operatorname{sgn}([c(\sigma^0),c(\sigma^1)],\sigma^1)$ \\
  $\cdot\operatorname{sgn}([c(\sigma^0),c(\sigma^1),c(\sigma^2)],\sigma^2)$ \\
          $\cdot[c(\sigma^1),c(\sigma^2)]$
         }} &
  \includegraphics[scale=0.15]{orient1b} &
  \includegraphics[scale=0.15]{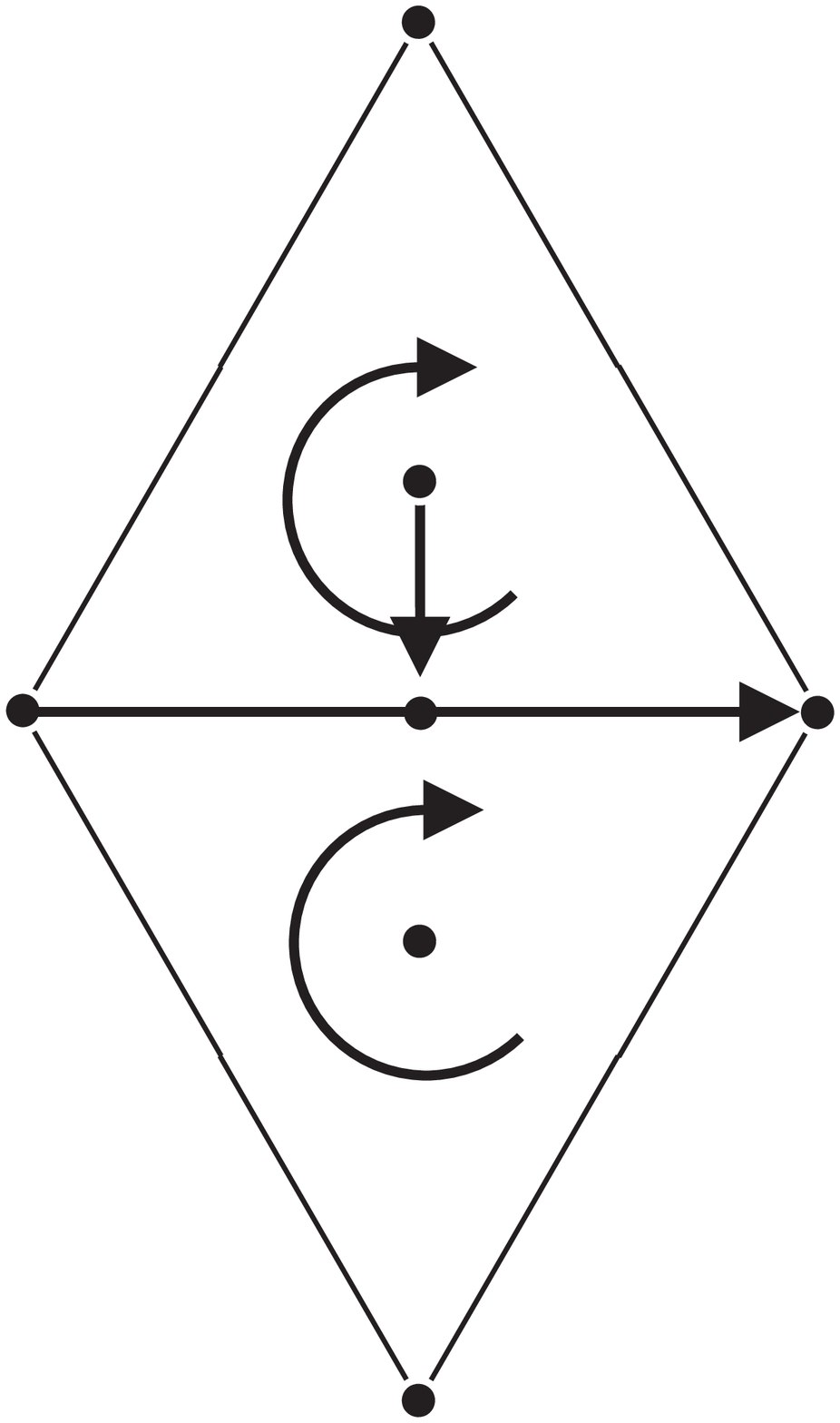} &
  \includegraphics[scale=0.15]{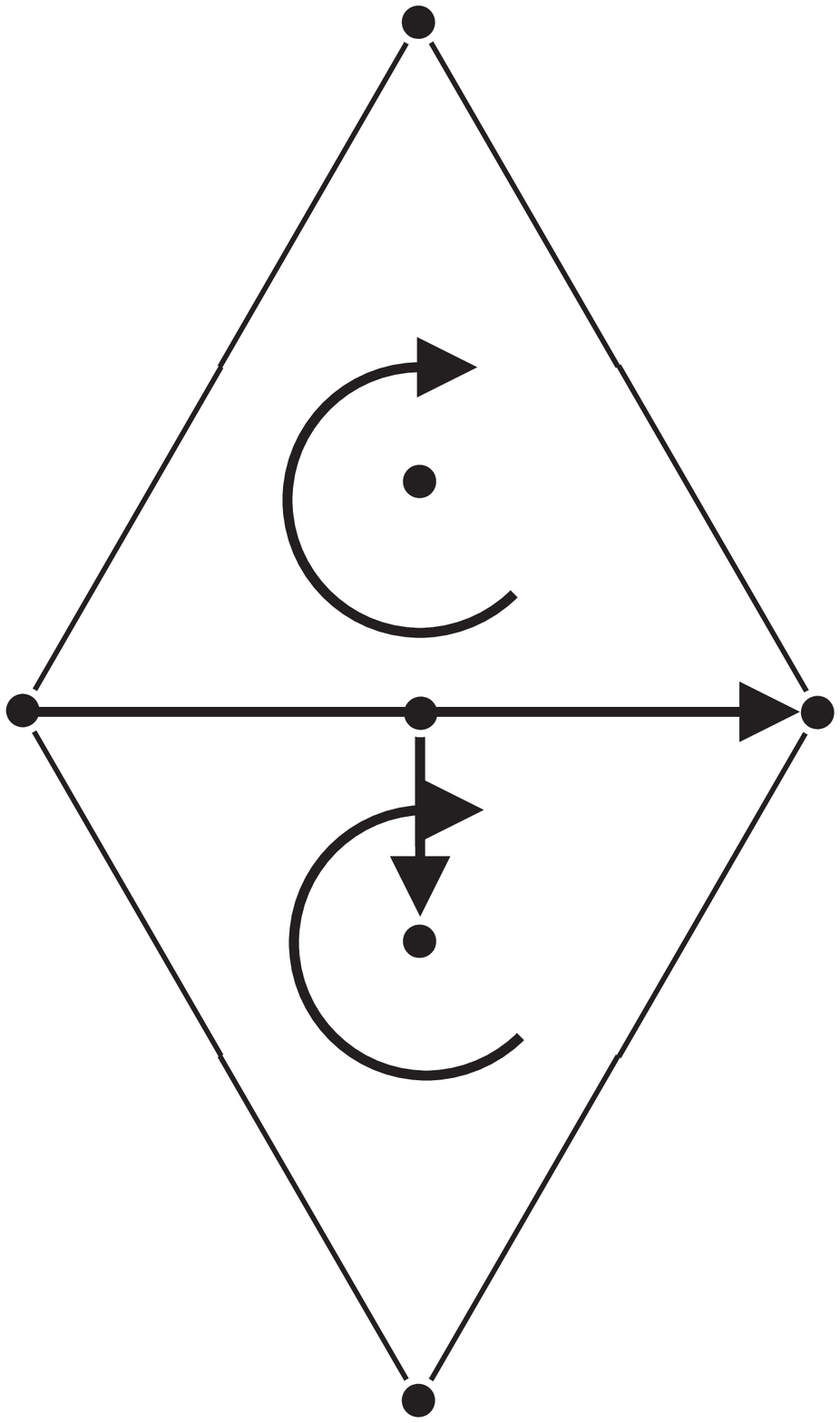} &
  \includegraphics[scale=0.15]{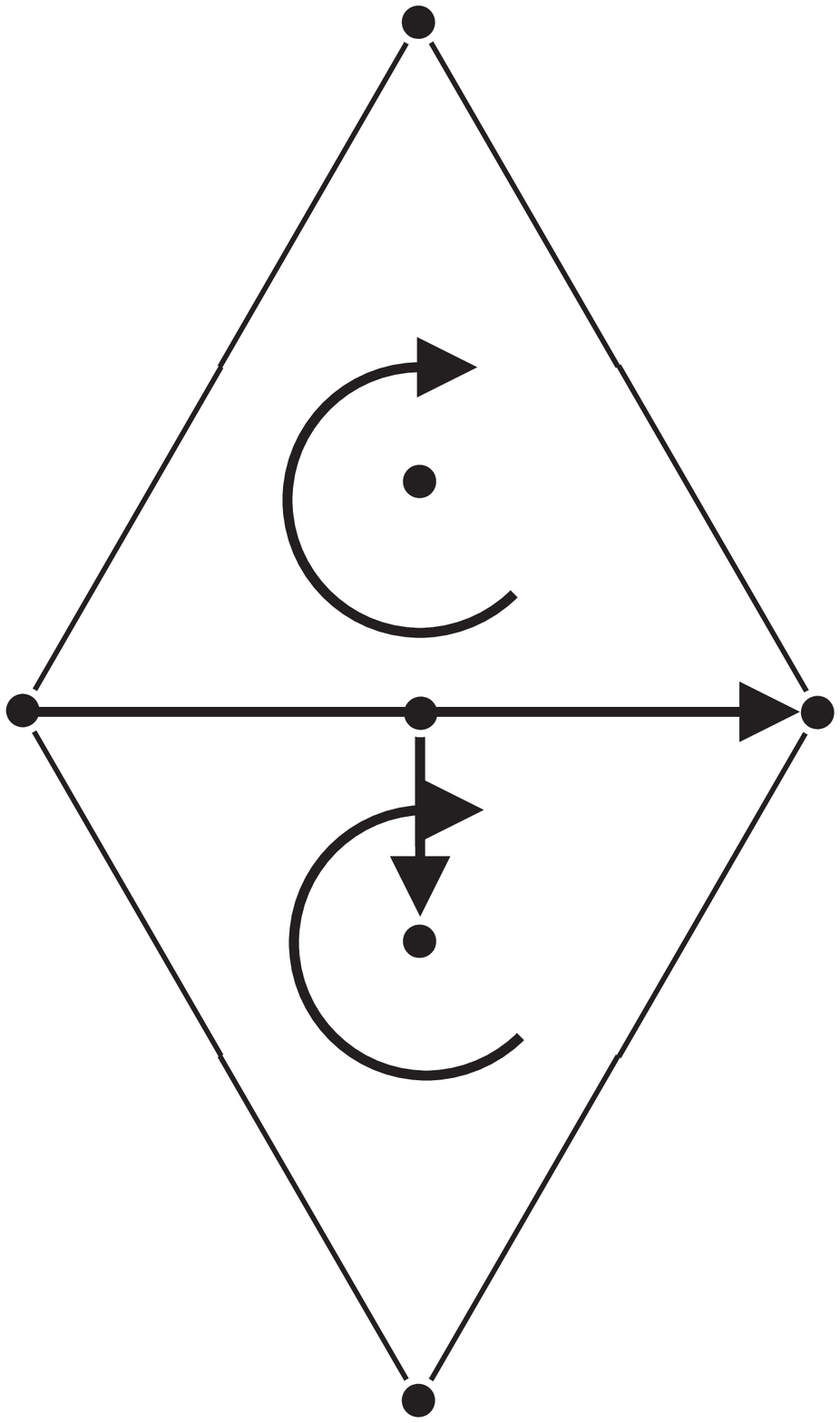}\\ \hline
  \end{tabular}
  \end{center}\end{table}
\end{example}

\paragraph{Orientation of the Dual of a Dual Cell.}\index{orientation!dual of dual}
While the circumcentric duality operator is a map from the primal
simplicial complex to the dual cell complex, we can formally
extend the circumcentric duality operator to a map from the dual
cell complex to the primal simplicial complex. However, we need to
be slightly careful about the orientation of primal simplex we
recover from applying the circumcentric duality operator twice.

We have that, $\star\star(\sigma^k)=\pm\sigma^k$, where the sign
is chosen to ensure the appropriate choice of orientation. If, as
before, $\sigma^k$ has an orientation represented by
$dx^1\wedge\ldots\wedge dx^k$, and $\star\sigma^k$ has an
orientation represented by $dx^{k+1}\wedge\ldots\wedge dx^n$, then
the orientation of $\star\star(\sigma^k)$ is chosen so that
$dx^{k+1}\wedge\ldots\wedge dx^n\wedge dx^1\wedge\ldots\wedge
dx^k$ is consistent with the ambient volume-form. Since, by
construction, $\star(\sigma^k)$, $dx^1\wedge\ldots\wedge dx^n$
has an orientation consistent with the ambient volume-form, we
need only compare $dx^{k+1}\wedge\ldots\wedge dx^n\wedge
dx^1\wedge\ldots\wedge dx^k$ with $dx^1\wedge\ldots\wedge dx^n$.
Notice that it takes $n-k$ transpositions to get the $dx^1$ term
in front of the $dx^{k+1}\wedge\ldots\wedge dx^n$ terms, and we
need to do this $k$ times for each term of $dx^1\wedge\ldots\wedge
dx^k$, so it follows that the sign is simply given by
$(-1)^{k(n-k)}$, or equivalently,
\begin{equation} \label{E:starstarsigma}
\star\star(\sigma^k)= (-1)^{k(n-k)}\sigma^k.
\end{equation}
A similar relationship holds if we use a dual cell instead of the
primal simplex $\sigma^k$.

\paragraph{Support Volume of a Primal Simplex and Its Dual Cell.}
We can think of a cochain as being constructed out of a basis
consisting of cosimplices or cocells with value $1$ on a single
simplex or cell, and $0$ otherwise. The way to visualize this
cosimplex is that it is associated with a differential form that
has support on what we will refer to as the {\bfi support
volume}\index{support volume} associated with a given simplex or
cell.

\begin{definition}
The \textbf{support volume}\index{support volume} of a simplex $\sigma^k$ is a
$n$-volume given by the convex hull of the geometric union of the
simplex and its circumcentric dual. This is given by
\[V_{\sigma^k}=\operatorname{convex hull}(\sigma^k,\star\sigma^k)\cap|K|.\]

The intersection with $|K|$ is necessary to ensure that the
support volume does not extend beyond the polytope $|K|$ which
would otherwise occur if $|K|$ is nonconvex.

We extend the notion of a support volume to a dual cell
$\star\sigma^k$ by similarly defining
\[V_{\star\sigma^k}=\operatorname{convex hull}(\star\sigma^k,\star\star\sigma^k)\cap|K|=V_{\sigma^k}.\]
\end{definition}

To clarify this definition, we will consider some examples of
simplices, their dual cells, and their corresponding support
volumes. For two-dimensional simplicial complexes, this is illustrated in Table~\ref{dec:table:2d_primal_dual_support}.

\begin{table}[h!]
\caption{\label{dec:table:2d_primal_dual_support}\index{simplex!example}\index{cell!example}\index{support volume!example}Primal simplices, dual cells, and support volumes in two dimensions.}
\begin{center}
\begin{tabular}{|c|c|c|}
  \hline
  Primal Simplex& Dual Cell & Support Volume\\ \hline
  \includegraphics[scale=0.2,clip=true]{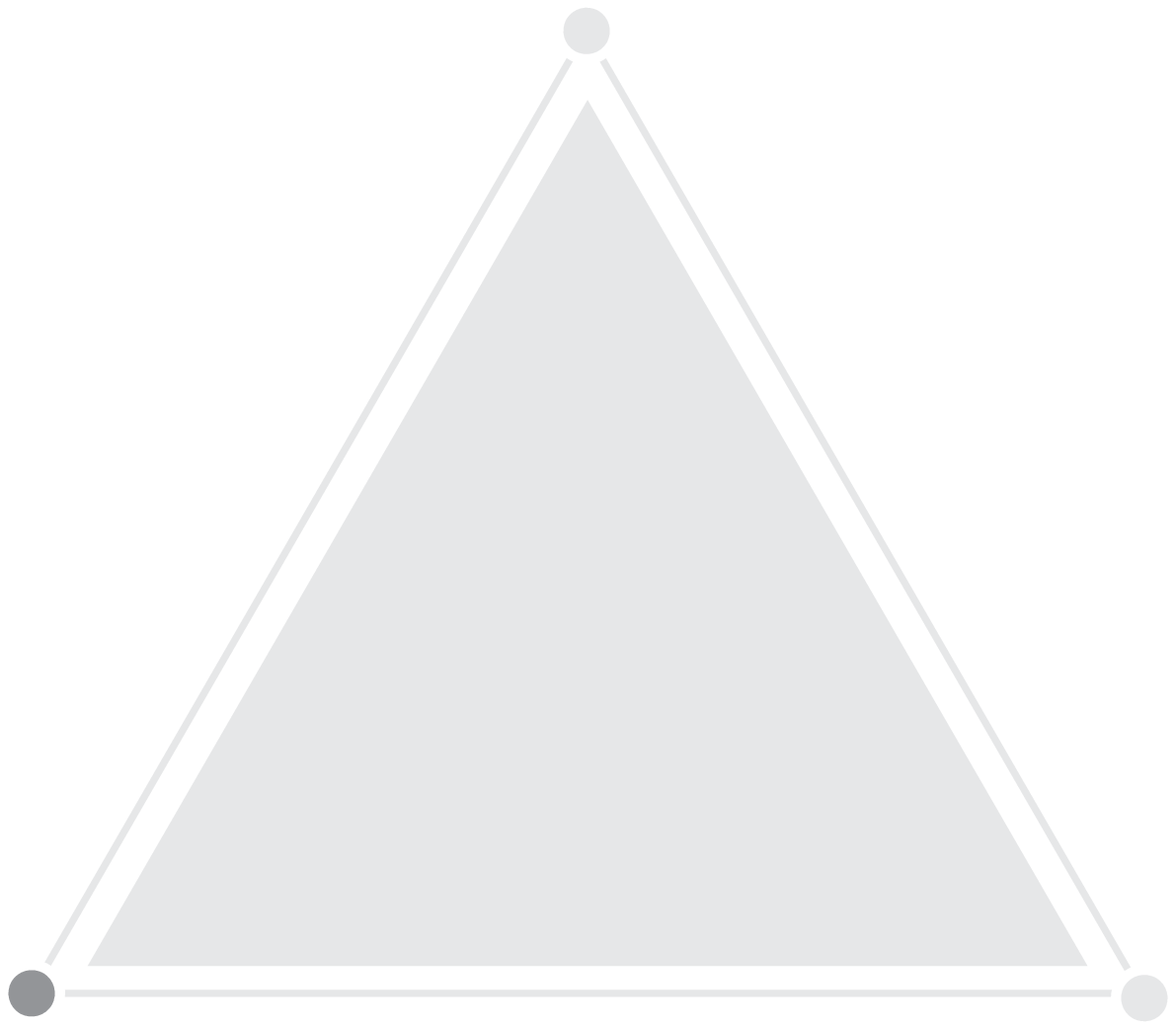} &
  \includegraphics[scale=0.2,clip=true]{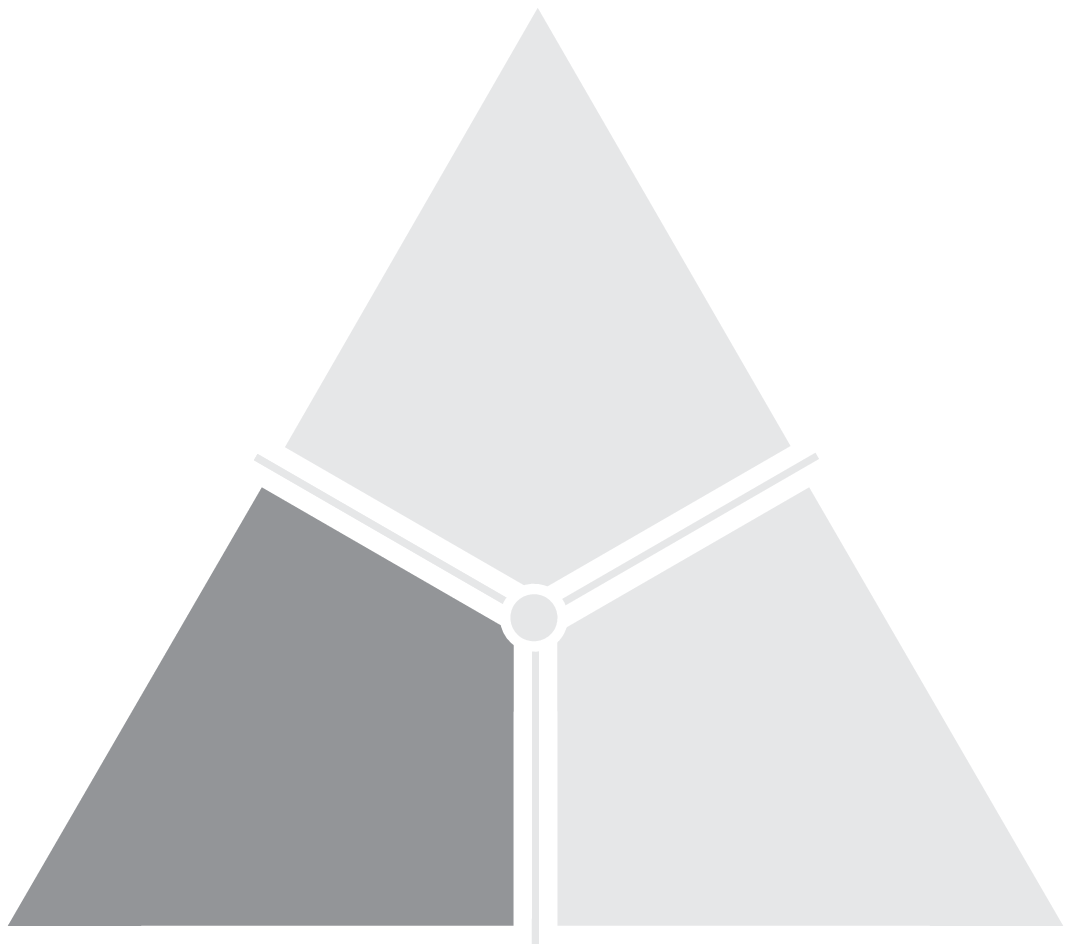} &
  \includegraphics[scale=0.2,clip=true]{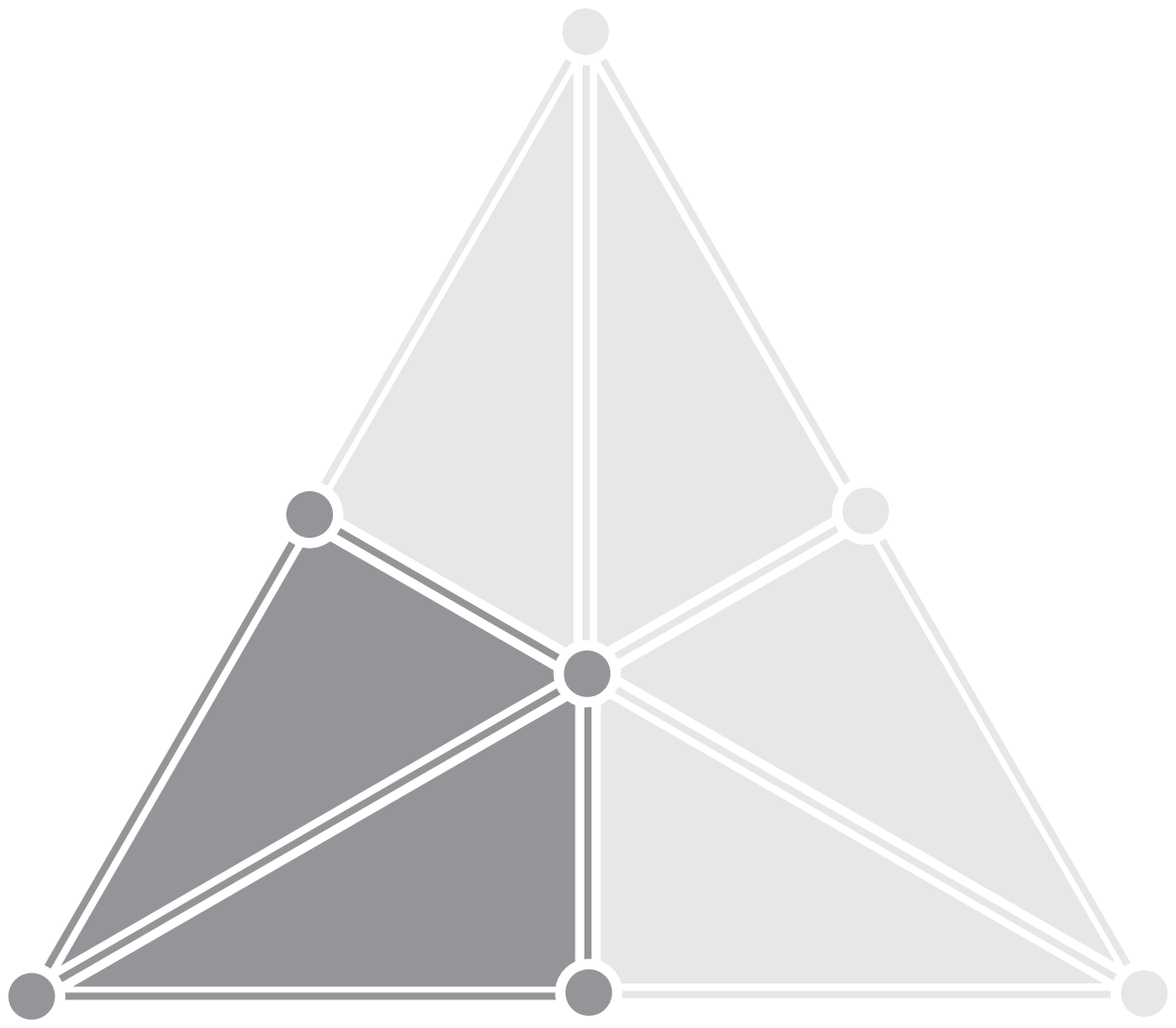}\\
  $\sigma^0$, $0$-simplex
  & $\star\sigma^0$, $2$-cell & $V_{\sigma^0}=V_{\star\sigma^0}$\\ \hline
  \includegraphics[scale=0.2,clip=true]{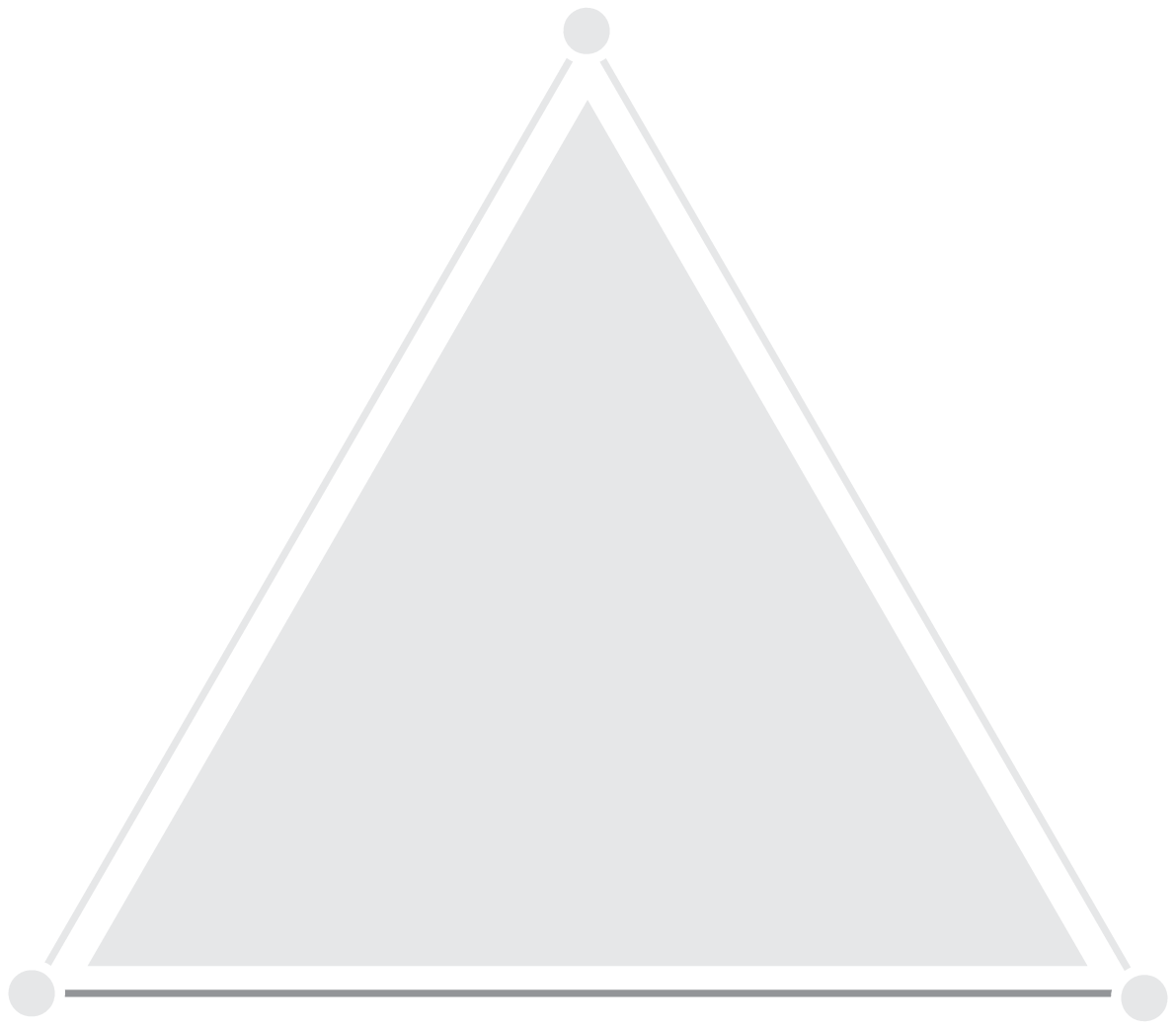} &
  \includegraphics[scale=0.2,clip=true]{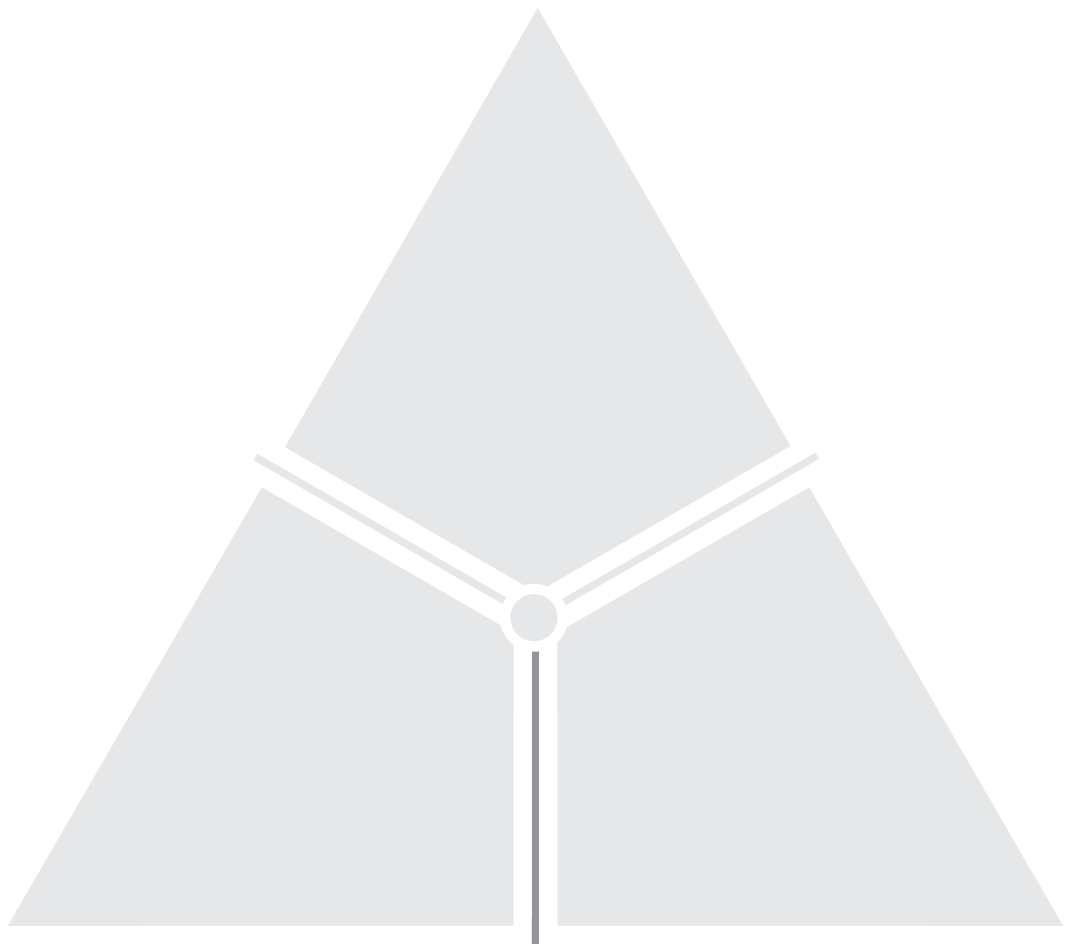} &
  \includegraphics[scale=0.2,clip=true]{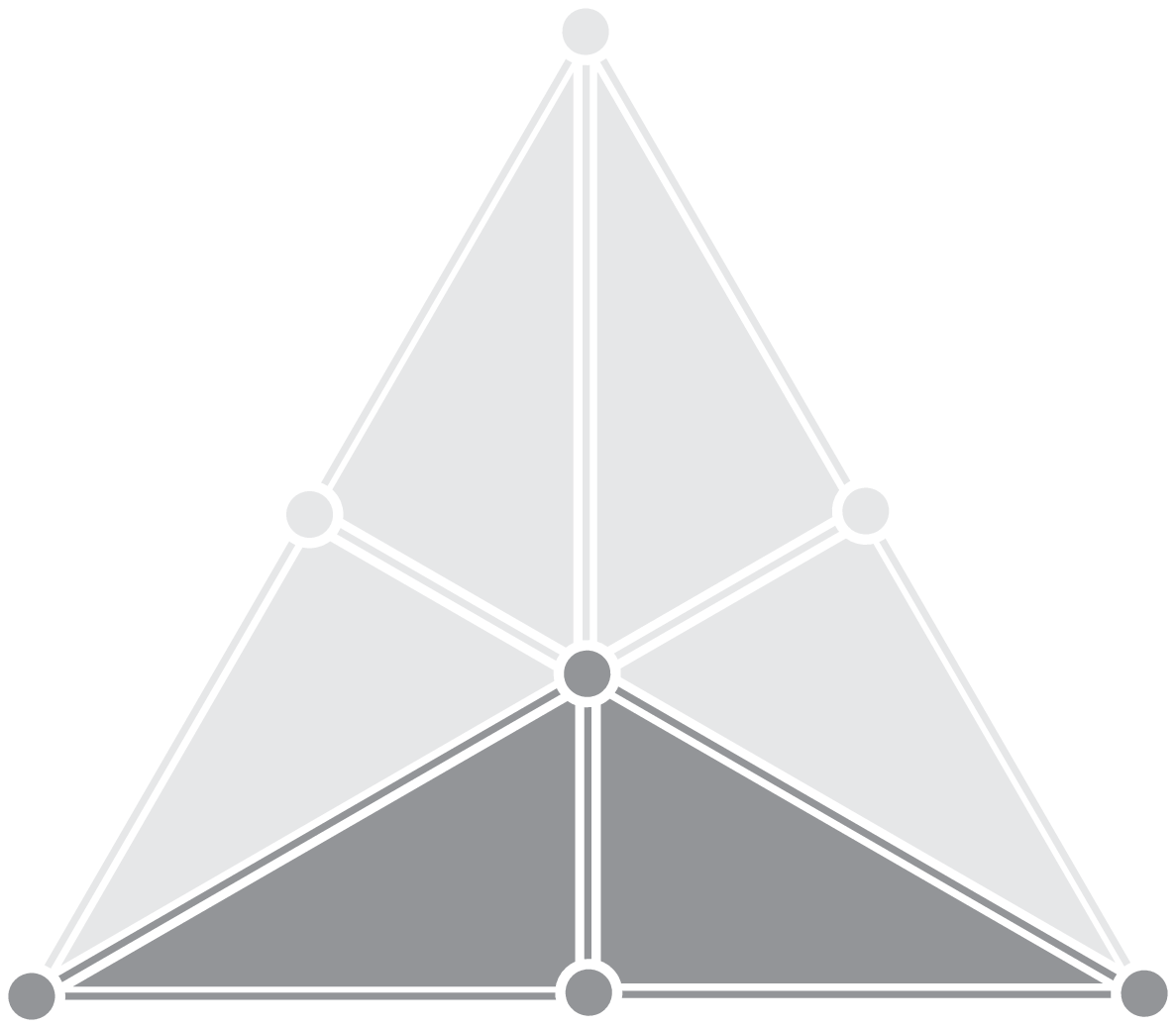}\\
  $\sigma^1$, $1$-simplex
  & $\star\sigma^1$, $1$-cell & $V_{\sigma^1}=V_{\star\sigma^1}$\\ \hline
  \includegraphics[scale=0.2,clip=true]{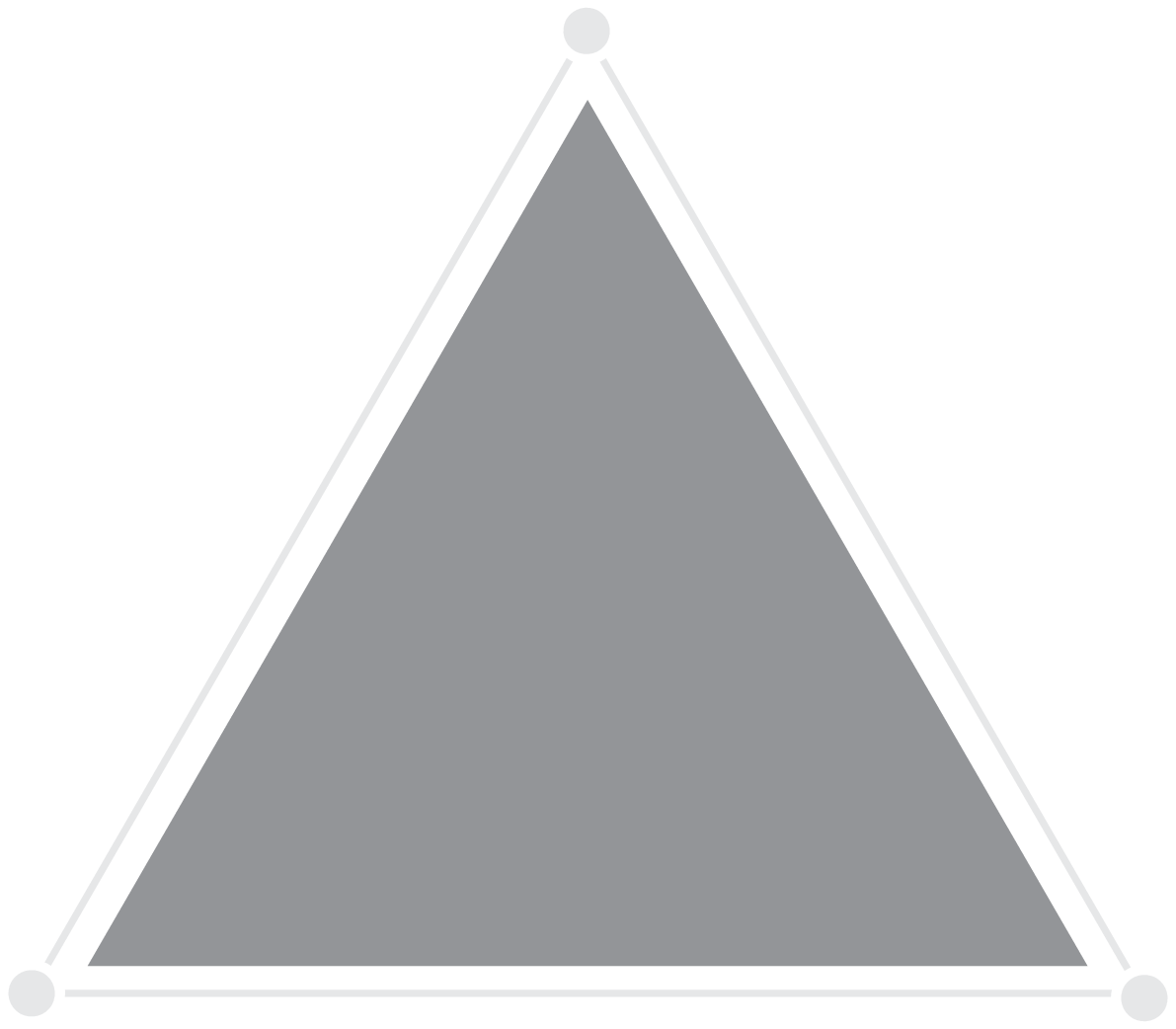} &
  \includegraphics[scale=0.2,clip=true]{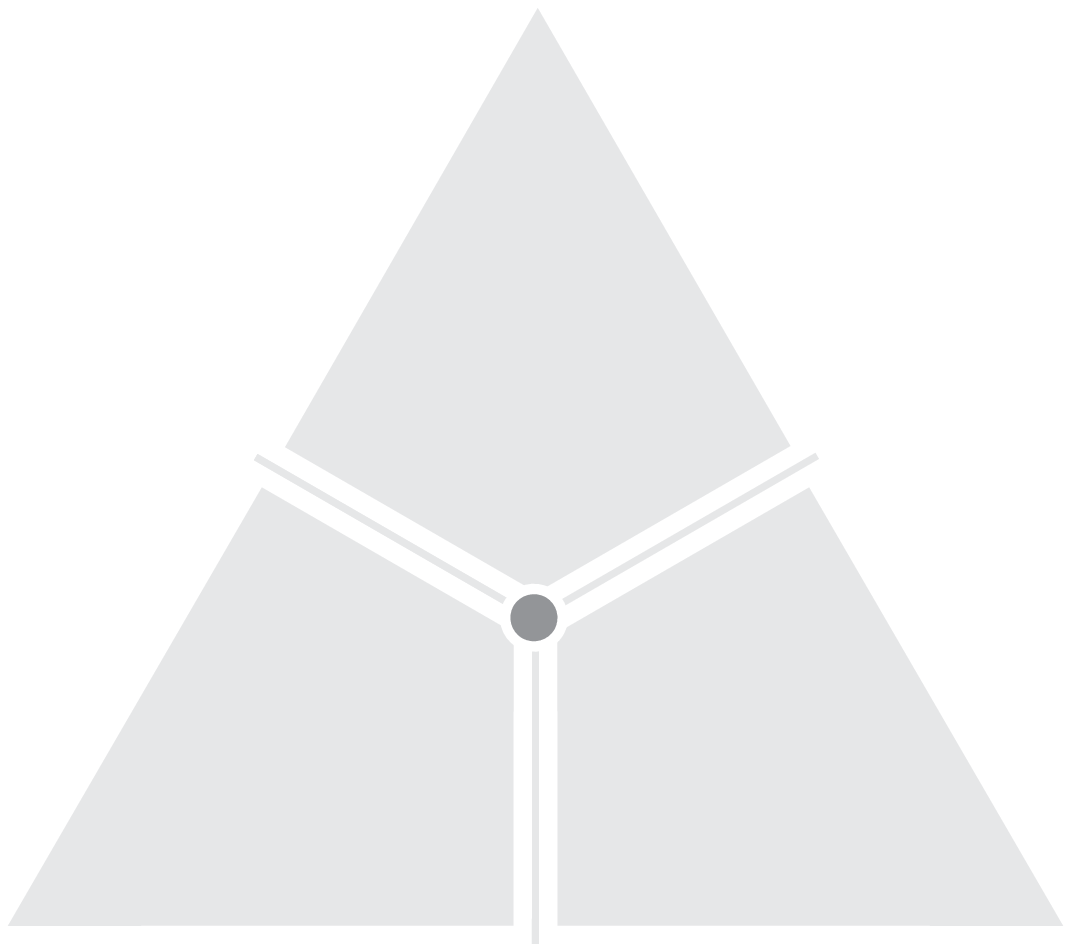} &
  \includegraphics[scale=0.2,clip=true]{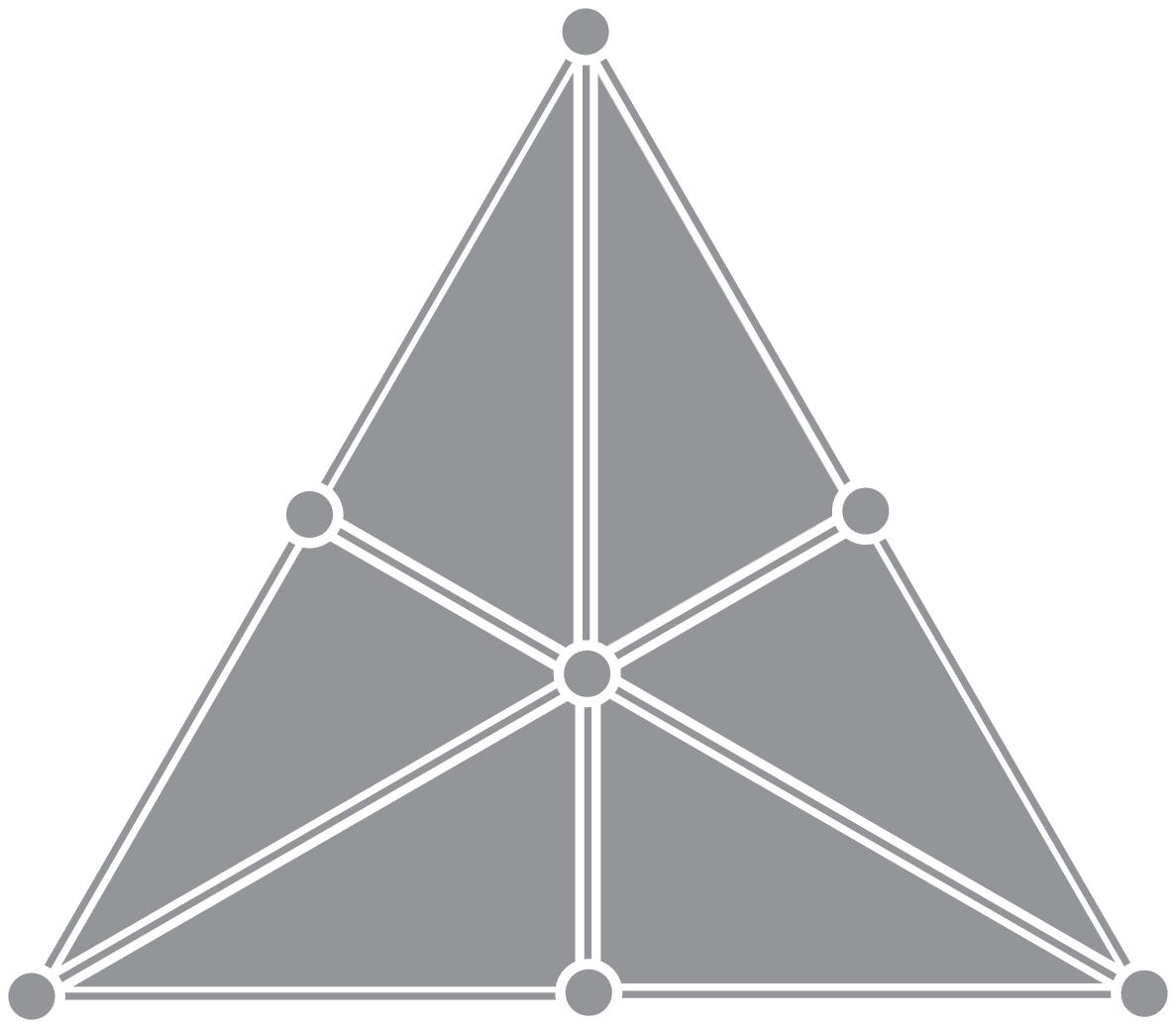} \\
  $\sigma^2$, $2$-simplex
  & $\star\sigma^2$, $0$-cell & $V_{\sigma^2}=V_{\star\sigma^2}$\\ \hline
\end{tabular}
\end{center}
\end{table}
The support volume has the nice property that at each dimension,
it partitions the polytope $|K|$ into distinct non-intersecting
regions associated with each individual $k$-simplex. For any two
distinct $k$-simplices, the intersection of their corresponding
support volumes have measure zero, and the union of the support
volumes of all $k$-simplices recovers the original polytope $|K|$.

Notice, from our construction, that the support volume of a
simplex and its dual cell are the same, which suggests that there
is an identification between cochains on $k$-simplices and
cochains on $(n-k)$-cells. This is indeed the case, and is a
concept associated with the Hodge star for differential forms.

Examples of simplices, their dual cells, and the
corresponding support volumes in three dimensions are given in Table~\ref{dec:table:3d_primal_dual_support}.

\begin{table}[h!]
\caption{\label{dec:table:3d_primal_dual_support}\index{simplex!example}\index{cell!example}\index{support volume!example}Primal simplices, dual cells, and support volumes in three dimensions.}
\begin{center}
\begin{tabular}{|c|c|c|}
  \hline
  Primal Simplex& Dual Cell & Support Volume\\ \hline
  \includegraphics[viewport=0 0 201 253,clip,scale=0.35]{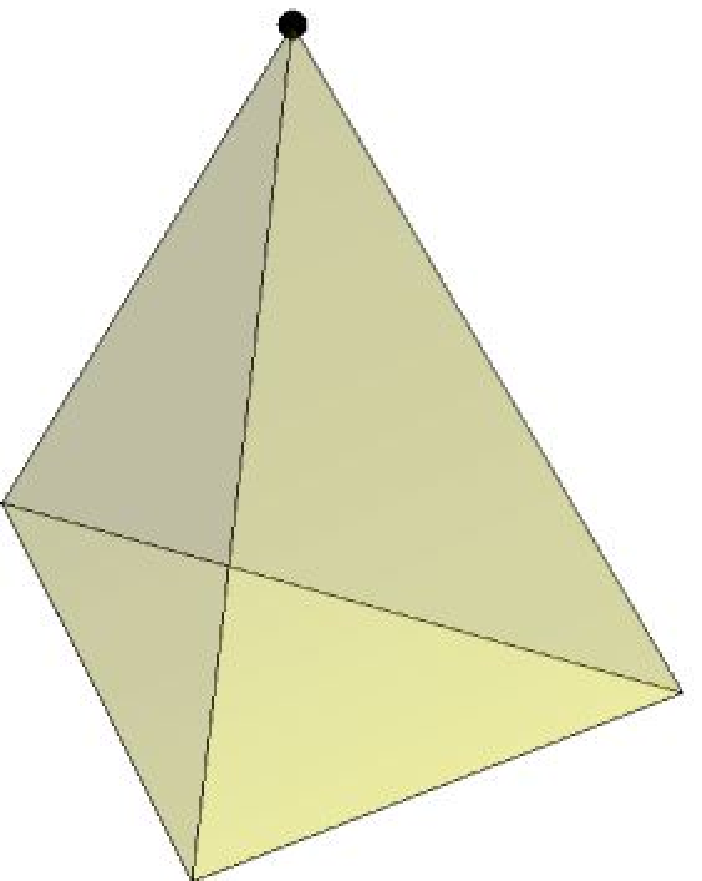} &
  \includegraphics[viewport=0 0 201 253,clip,scale=0.35]{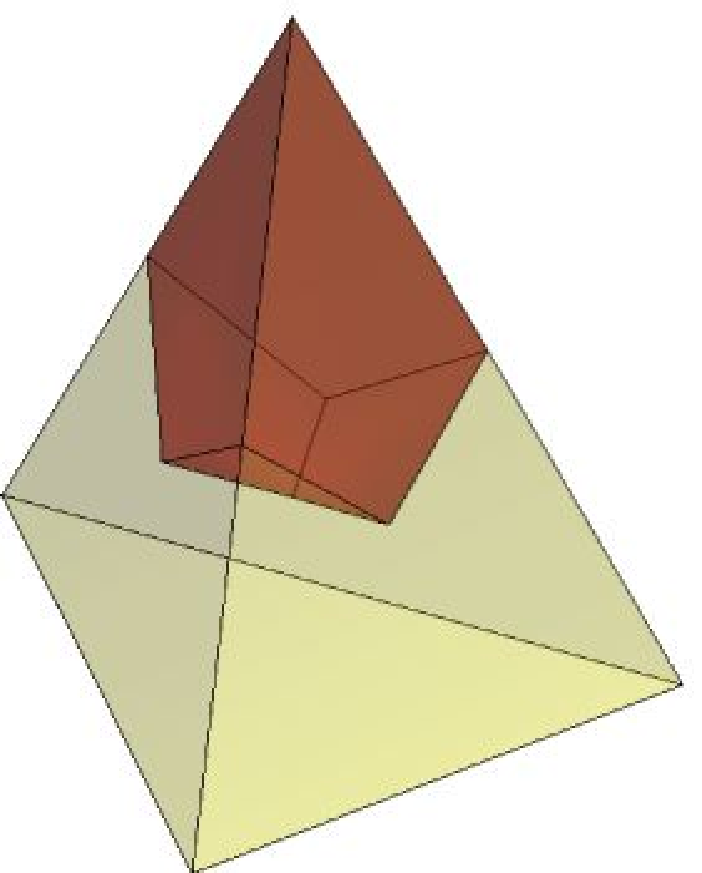} &
  \includegraphics[viewport=0 0 201 253,clip,scale=0.35]{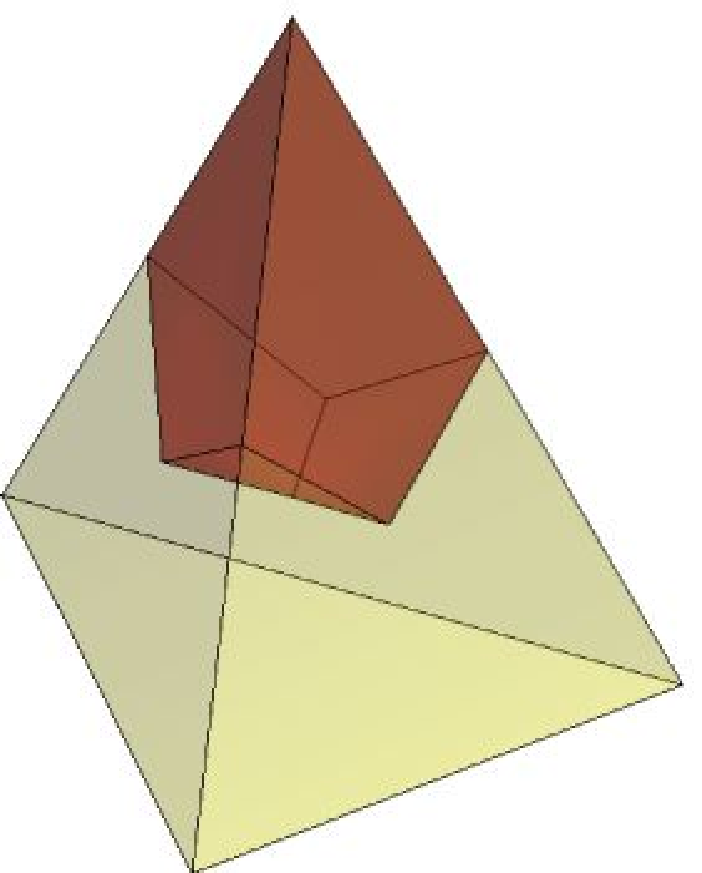}\\ $\sigma^0$,
  $0$-simplex & $\star\sigma^0$, $3$-cell & $V_{\sigma^0}=V_{\star\sigma^0}$\\
  \hline
%
  \includegraphics[viewport=0 0 205 253,clip,scale=0.35]{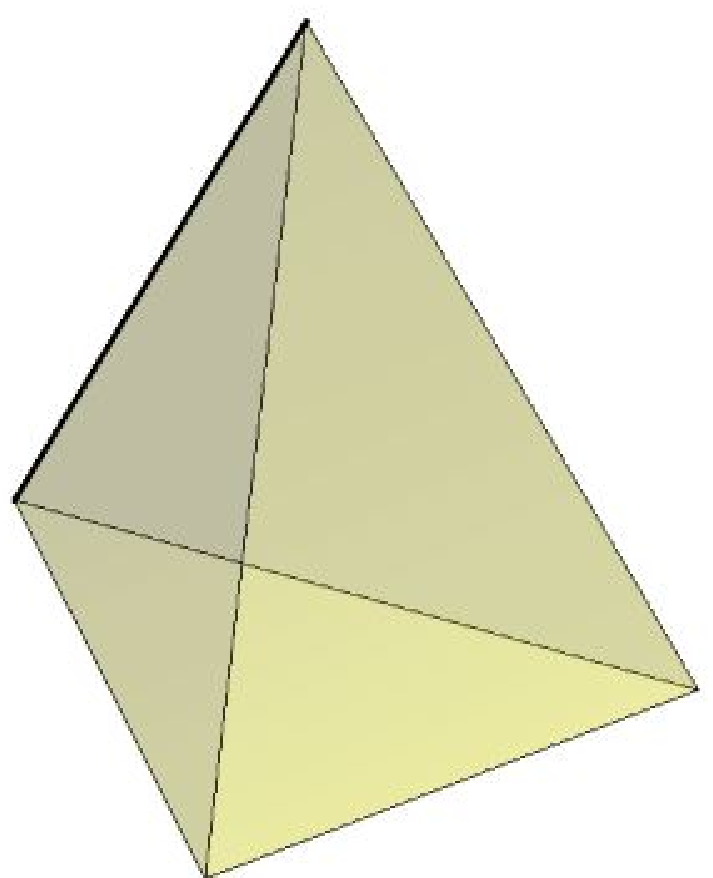} &
  \includegraphics[viewport=0 0 201 253,clip,scale=0.35]{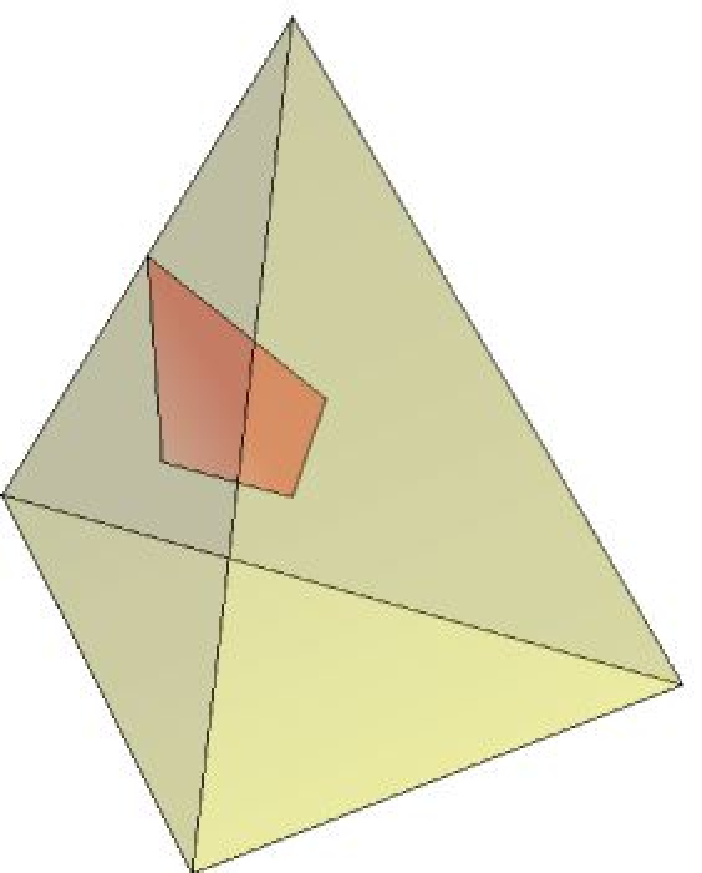} &
  \includegraphics[viewport=0 0 201 253,clip,scale=0.35]{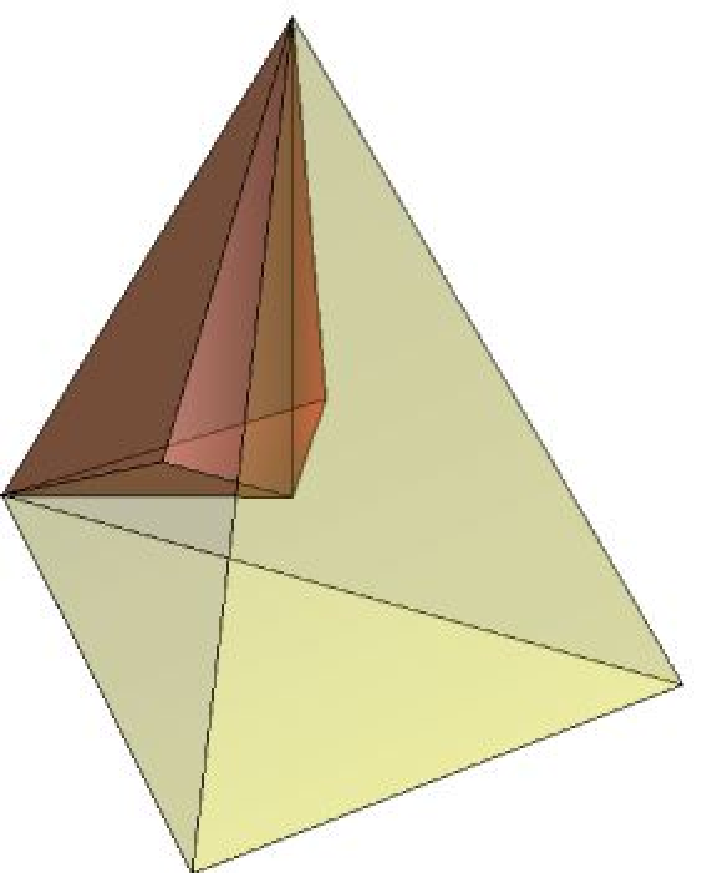}\\
  $\sigma^1$, $1$-simplex
  & $\star\sigma^1$, $2$-cell & $V_{\sigma^1}=V_{\star\sigma^1}$\\ \hline
  \includegraphics[viewport=0 0 206 255,clip,scale=0.35]{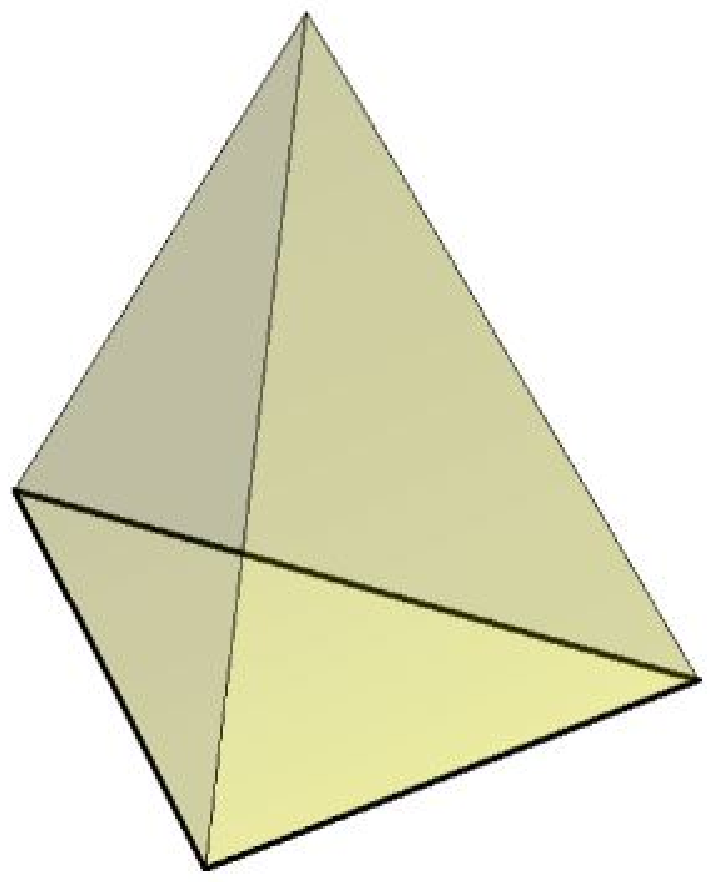} &
  \includegraphics[viewport=0 0 201 253,clip,scale=0.35]{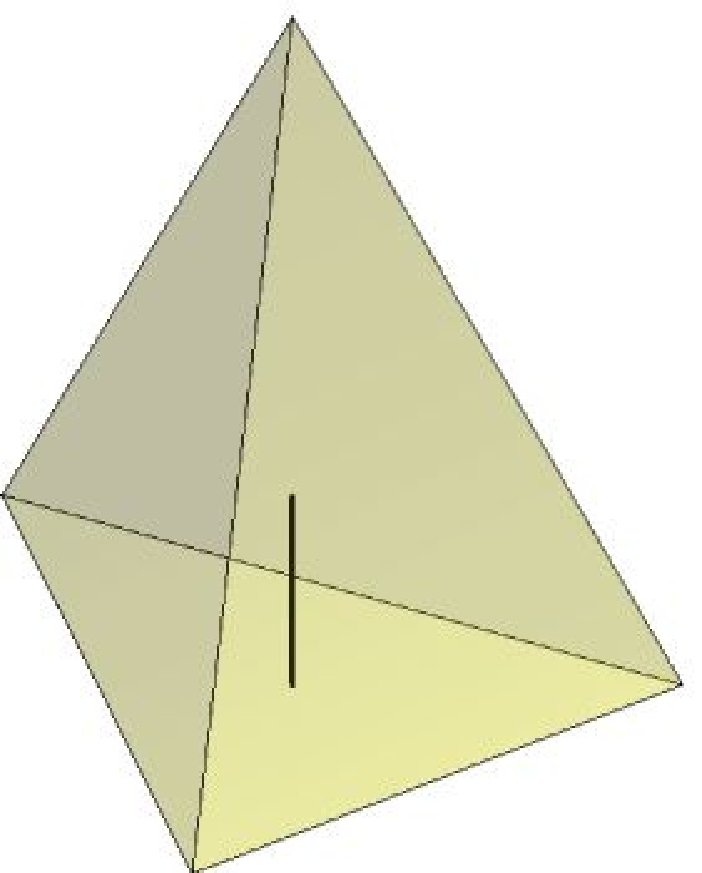} &
  \includegraphics[viewport=0 0 201 253,clip,scale=0.35]{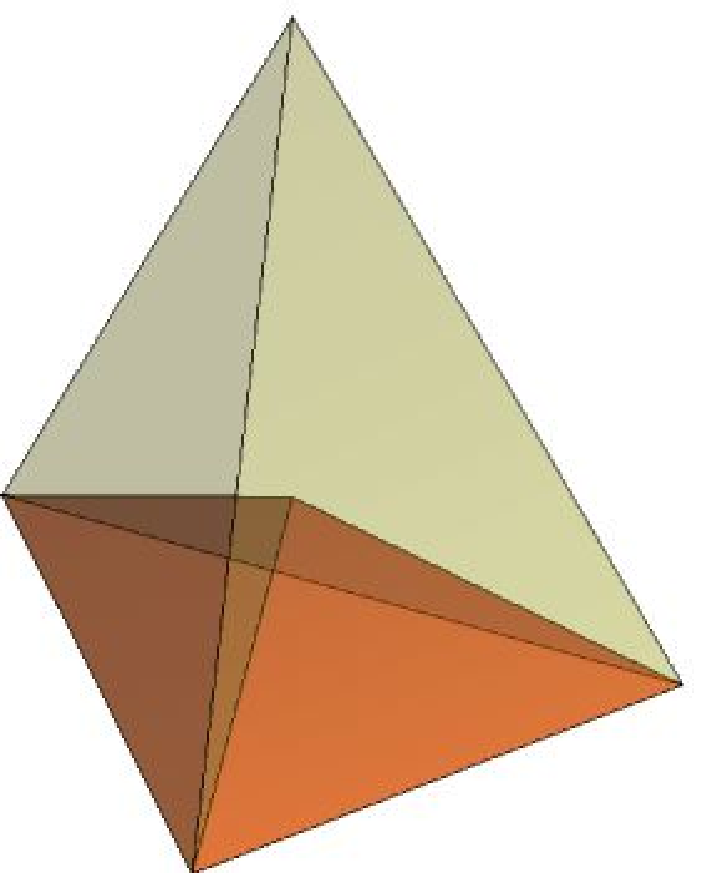} \\ $\sigma^2$,
  $2$-simplex & $\star\sigma^2$, $1$-cell & $V_{\sigma^2}=V_{\star\sigma^2}$\\
  \hline
  \includegraphics[viewport=0 0 206 257,clip,scale=0.35]{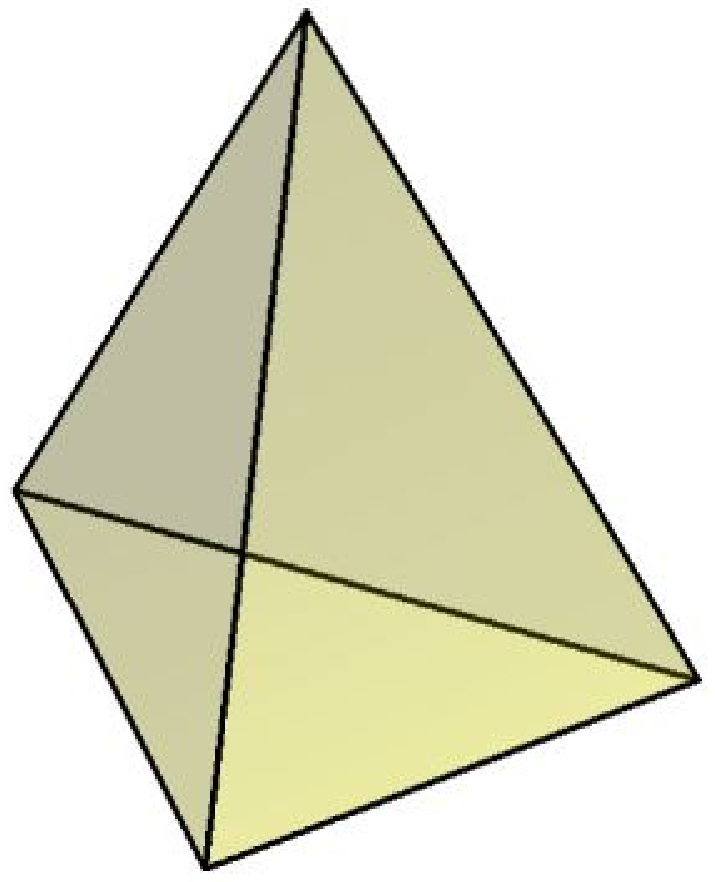} &
  \includegraphics[viewport=0 0 201 253,clip,scale=0.35]{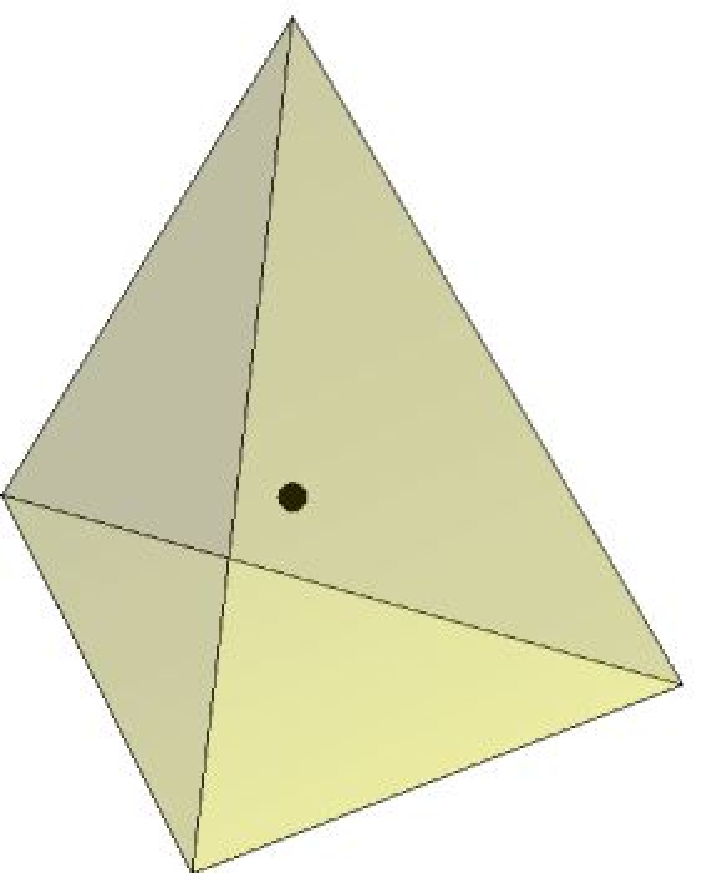} &
  \includegraphics[viewport=0 0 201 253,clip,scale=0.35]{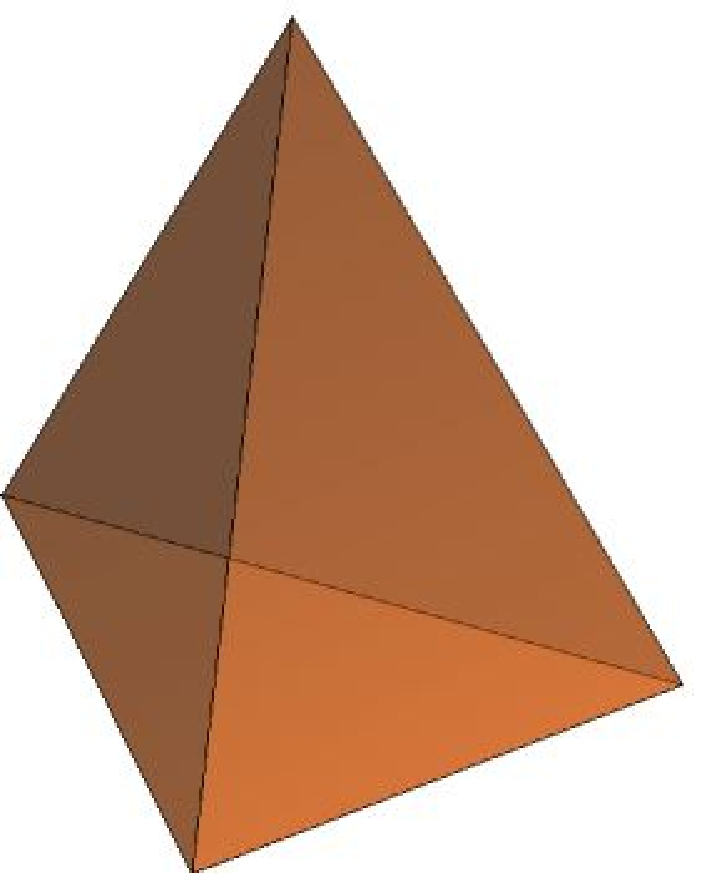} \\ $\sigma^3$,
  $3$-simplex & $\star\sigma^3$, $0$-cell & $V_{\sigma^3}=V_{\star\sigma^3}$\\
  \hline
\end{tabular}
\end{center}
\end{table}

In our subsequent discussion, we will assume that we are given a
simplicial complex $K$ of dimension $n$ in $\RR^N$. Thus, the
highest-dimensional simplex in the complex is of dimension $n$ and
each $0$-simplex (vertex) is in $\RR^N$. One can obtain this, for
example, by starting from $0$-simplices, i.e., vertices, and then
constructing a Delaunay triangulation, using the vertices as sites. Often,
our examples will be for two-dimensional discrete surfaces in
$\RR^3$ made up of triangles (here $n = 2$ and $N=3$) or three-dimensional manifolds made of tetrahedra, possibly embedded in a
higher-dimensional space.

\paragraph{Cell Complexes.}
The circumcentric dual of a primal simplicial complex is an
example of a cell complex. The definition of a cell complex
follows.

\begin{definition}
A \textbf{cell complex}\index{complex!cell|see{cell, complex}}\index{cell!complex} $\star K$ in $\RR^N$ is a collection of
cells in $\RR^N$ such that,
\begin{enumerate}
\item There is a partial ordering of cells in $\star K$,
$\hat\sigma^k\prec\hat\sigma^l$, which is read as $\hat\sigma^k$
is a face of $\hat\sigma^l$. \item The intersection of any two
cells in $\star K$, is either a face of each of them, or it is
empty. \item The boundary of a cell is expressible as a sum of its
proper faces.
\end{enumerate}
\end{definition}

We will see in the next section that the notion of boundary in the
circumcentric dual has to be modified slightly from the geometric
notion of a boundary in order for the circumcentric dual to be
made into a cell complex.

\section{Local and Global Embeddings}\label{dec:sec:local_global_embeddings}
While it is computationally more convenient to have a global
embedding of the simplicial complex into a higher-dimensional
ambient space to account for non-flat manifolds\index{manifolds!non-flat} it suffices to
have an abstract simplicial complex along with a local metric on
vertices. The metric is local in the sense that distances between two
vertices are only defined if they are part of a common $n$-simplex
in the abstract simplicial complex. Then, the local metric is a map
$d:\{(v_0,v_1)\mid v_0,v_1\in K^{(0)}, [v_0,v_1]\prec\sigma^n\in
K\}\rightarrow\RR$.

The axioms for a {\bfi local metric}\index{metric!local} are as follows,
\begin{description}
\item[Positive.] $d(v_0,v_1)\geq 0$, and $d(v_0,v_0)=0$, $\forall
[v_0,v_1]\prec \sigma^n \in K$. \item[Strictly Positive.] If
$d(v_0,v_1)=0$, then $v_0=v_1$, $\forall [v_0,v_1]\prec\sigma^n\in
K$. \item[Symmetry.] $d(v_0,v_1)=d(v_1,v_0)$, $\forall
[v_0,v_1]\prec\sigma^n\in K$. \item[Triangle
Inequality.]$d(v_0,v_2)\leq d(v_0,v_1)+d(v_1,v_2)$, $\forall
[v_0,v_1,v_2]\prec\sigma^n\in K$.
\end{description}
This allows us to embed each $n$-simplex locally into $\RR^n$, and
thereby compute all the necessary metric dependent quantities in
our formulation. For example, the volume of a $k$-dual cell will
be computed as the sum of the $k$-volumes of the dual cell
restricted to each $n$-simplex in its local embedding into
$\RR^n$.

This notion of local metrics and local embeddings is consistent
with the point of view that exterior calculus is a local theory
with operators that operate on objects in the tangent and
cotangent space of a fixed point. The issue of comparing objects
in different tangent spaces is addressed in the discrete theory of
connections on principal bundles in \cite{LeMaWe2003}.

This also provides us with a criterion for evaluating a global
embedding. The embedding should be such that the metric of the
ambient space $\RR^N$ restricted to the vertices of the complex,
thought of as points in $\RR^N$, agrees with the local metric
imposed on the abstract simplicial complex. A global embedding
that satisfies this condition will produce the same numerical
results in discrete exterior calculus as that obtained using the
local embedding method.

It is essential that the metric condition we impose is local,
since the notion of distances between points in a manifold which
are far away is not a well-defined concept, nor is it particularly
useful for embeddings. As the simple example below illustrates,
there may not exist any global embeddings into Euclidean space
that satisfies a metric constraint imposed for all possible pairs
of vertices.

\begin{example} \label{example:global_nonexistence}
Consider a circle, with the distance between two points given by
the minimal arc length. Consider a discretization given by 4
equidistant points on the circle, labelled $v_0,\ldots, v_3$, with
the metric distances as follows,
\begin{align*}
d(v_i,v_{i+1})=1, d(v_i,v_{i+2})=2,
\end{align*}
where the indices are evaluated modulo 4, and this distance
function is extended to a metric on all pairs of vertices by
symmetry. It is easy to verify that this distance function is
indeed a metric on vertices.
\newcounter{node}
\newcommand{\Letter}{\setcounter{node}{\xypolynode}\addtocounter{node}{-1}\arabic{node}}
\[\begin{xy} /r15mm/:
  \xypolygon4"I"{
    ~={90}
    ~*{\xybox{*{\Letter}*\cir<2mm>{}}}
    }
  \end{xy}
  \qquad\qquad
  \begin{xy} /r15mm/:
    \xypolygon4"I"{~={90}~>{}}
    \ar@{-}_{1}"I1";"I2"
    \ar@{.}"I2";"I3"
    \ar@{.}"I3";"I4"
    \ar@{.}"I4";"I1"
    \ar@{-}_{2}"I1";"I3"
    \ar@{.}"I2";"I4"
  \end{xy}\]
If we only use the local metric constraint, then we only require
that adjacent vertices are separated by $1$, and the following is
an embedding of the simplicial complex into $\RR^2$,
\[\begin{xy} /r15mm/:
  \xypolygon4"I"{
    ~={90}~>{}
    }
    \ar@{-}_{1}"I1";"I2"
    \ar@{-}_{1}"I2";"I3"
    \ar@{-}_{1}"I3";"I4"
    \ar@{-}_{1}"I4";"I1"
  \end{xy}\]
If, however, we use the metric defined on all possible pairs of vertices, by considering $v_0, v_1, v_2$, we have that $d(v_0,v_1)+d(v_1,v_2)=d(v_0,v_2)$. Since we are embedding these points into a Euclidean space, it follows that $v_0, v_1, v_2$ are collinear.

Similarly, by considering $v_0, v_2, v_3$, we conclude that they
are collinear as well, and that $v_1, v_3$ are coincident, which
contradicts $d(v_1,v_3)=2$. Thus, we find that there does not exist
a global embedding of the circle into Euclidean space if we
require that the embedding is consistent with the metric on
vertices defined for all possible pairs of vertices.
\end{example}

\section{Differential Forms and Exterior Derivative}
\label{dec:sec:Differential}

We will now define discrete differential forms. We will use some
terms (which we will define) from algebraic topology, but it will
become clear by looking at the examples that one can gain a clear
and working notion of what a discrete form is without any
algebraic topology. We start with a few definitions for which more
details can be found on page 26 and 27 of \cite{Mu1984}.

\begin{definition}
Let $K$ be a simplicial complex. We denote the free abelian group
generated by a basis consisting of oriented $k$-simplices by
$C_k\left(K;\mathbb{Z}\right).$ This is the space of finite formal
sums of the $k$-simplices, with coefficients in $\mathbb{Z}$.
Elements of $C_k(K;\ZZ)$ are called \textbf{$k$-chains}\index{chain}.
\end{definition}

\begin{example}\index{chain!example} Figure~\ref{dec:fig:chains} shows examples of $1$-chains and $2$-chains.
\begin{figure}[H]
\[\begin{xy} <1cm,0cm>:
(0.3,0.3)*[o]+{\circlearrowleft}*+!L{1},
(0.7,0.7)*[o]+{\circlearrowleft}*+!L{2},
(1.4,0.7)*[o]+{\circlearrowleft}*+!L{1},
(2,0.3)*[o]+{\circlearrowleft}*+!L{3},
(1.5,-0.75)*+{2\mathrm{-chain}}, (-3.5,-0.75)*+{1\mathrm{-chain}},
(0,0)*+@{*}="a2", (0,1)*+@{*}="b2", (1,0)*+@{*}="c2",
(1,1)*+@{*}="d2", (2.5,1.2)*+@{*}="e2", (3,-0.2)*+@{*}="f2",
(-2,1.3)*+@{*}="e1", (-2.4,0.2)*+@{*}="d1", (-3,1)*+@{*}="c1",
(-3.7,0.5)*+@{*}="b1", (-5,1.2)*+@{*}="a1", \ar@{-} "a2";"b2",
\ar@{-} "b2";"c2", \ar@{-} "c2";"a2", \ar@{-} "b2";"d2", \ar@{-}
"c2";"d2", \ar@{-} "c2";"e2", \ar@{-} "d2";"e2", \ar@{-}
"c2";"f2", \ar@{-} "e2";"f2", \ar^{1} "a1";"b1", \ar^{3}
"b1";"c1", \ar^{2} "c1";"d1", \ar_{5} "d1";"e1"
\end{xy}
\]
\caption{\label{dec:fig:chains}Examples of chains.}
\end{figure}
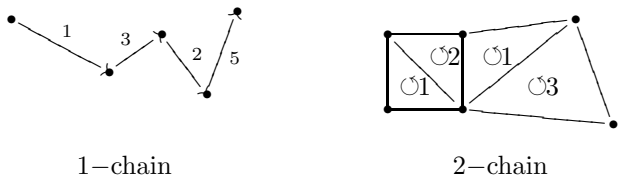
\end{example}

We view discrete $k$-forms as maps from the space of $k$-chains to
$\RR$. Recalling that the space of $k$-chains is a group, we
require that these maps be homomorphisms into
the additive group $\RR$. Thus, discrete forms are what are called
cochains in algebraic topology. We will define cochains below in
the definition of forms but for more context and more details
readers can refer to any algebraic topology text, for example, page
251 of \cite{Mu1984}.

This point of view of forms as cochains is not new.  The idea of
defining forms as cochains appears, for example, in the works of
\cite{Ad1996}, \cite{De1995}, \cite{Hi1999}, and
\cite{SeSeSeAd2000}. Our point of departure is that the
other authors go on to develop a theory of discrete exterior
calculus of forms only by introducing interpolation of forms, which
we will be able to avoid. The formal definition of discrete forms
follows.

\begin{definition}
A \textbf{primal discrete $k$-form}\index{form!primal} $\alpha$ is a homomorphism
from the chain group $C_k(K;\ZZ)$ to the additive group $\RR$.
Thus, a discrete $k$-form is an element of
$\operatorname{Hom}(C_k(K),\RR)$, the space of \textbf{cochains}\index{cochain}.
This space becomes an abelian group if we add two homomorphisms by
adding their values in $\RR$.  The standard notation for
$\operatorname{Hom}(C_k(K), \RR)$ in algebraic topology is $C^k(K;
\RR)$. But we will often use the notation $\Omega_d^k(K)$\index{$\Omega_d^k(K)$|see{form, primal}}\index{$\Omega_d^k(\star K)$|see{form, dual}} for this
space as a reminder that this is the space of discrete (hence the
$d$ subscript) $k$-forms on the simplicial complex $K$. Thus,
\[
    \Omega_d^k(K) := C^k(K;\RR) = \operatorname{Hom}(C_k(K),\RR) \, .
\]
\end{definition}

Note that, by the above definition, given a $k$-chain $\sum_i a_i c^k_i$
(where $a_i \in \ZZ$) and a discrete $k$-form $\alpha$, we have that
\[
    \alpha\left(\sum_i a_i c^k_i\right) = \sum_i a_i \alpha(c^k_i)\, ,
\]
and for two discrete $k$-forms $\alpha, \beta \in \Omega_d^k(K)$
and a $k$-chain $c \in C_k(K;\ZZ)$,
\[
    (\alpha + \beta)(c) = \alpha(c) + \beta(c) \, .
\]

In the usual exterior calculus on smooth manifolds integration of
$k$-forms on a $k$-dimensional manifold is defined in terms of the
familiar integration in $\RR^k$. This is done roughly speaking by
doing the integration in local coordinates, and showing that the
value is independent of the choice of coordinates, due to the
change of variables theorem in $\RR^k$. For details on this, see
the first few pages of Chapter 7 of \cite{AbMaRa1988}. We will not try to introduce the notion of
integration of discrete forms on a simplicial complex. Instead the
fundamental quantity that we will work with is the natural
bilinear pairing of cochains and chains, defined by evaluation.
More formally, we have the following definition.

\begin{definition}
The \textbf{natural pairing}\index{pairing!natural} of a $k$-form $\alpha$ and a
$k$-chain $c$ is defined as the bilinear pairing
\[
    \langle \alpha, c \rangle = \alpha(c).
\]
\end{definition}

As mentioned above, in discrete exterior calculus, this natural
pairing plays the role that integration of forms on chains plays in
the usual exterior calculus on smooth manifolds. The two are
related by a procedure done at the time of discretization. Indeed,
consider a simplicial triangulation $K$ of a polyhedron in
$\RR^n$, i.e., consider a ``flat''\index{manifolds!flat} discrete
manifold. If we are discretizing a continuous problem, we will have
some smooth forms defined in the space $|K| \subset \RR^n$.
Consider such a smooth $k$-form $\alpha^k$. In order to define the
discrete form $\alpha_d^k$ corresponding to $\alpha^k$, one would
integrate $\alpha^k$ on all the $k$-simplices in $K$. Then, the
evaluation of $\alpha_d^k$ on a $k$-simplex $\sigma^k$ is defined
by $\alpha_d^k(\sigma^k) := \int_{\sigma^k} \alpha^k$. Thus,
discretization is the only place where integration plays a role in
our discrete exterior calculus.

In the case of a non-flat manifold\index{manifolds!non-flat}, the situation is somewhat
complicated by the fact that the smooth manifold, and the
simplicial complex, as geometric sets embedded in the ambient space
do not coincide. A smooth differential form on the manifold
can be discretized into the cochain representation by identifying
the vertices of the simplicial complex with points on the
manifold, and then using a local chart to identify $k$-simplices
with $k$-volumes on the manifold.

There is the possibility of $k$-volumes overlapping even when
their corresponding $k$-simplices do not intersect, and this
introduces a discretization error that scales like the mesh size.
One can alternatively construct geodesic boundary surfaces in an
inductive fashion, which yields a partition of the manifold, but
this can be computationally prohibitive to compute.

Now we can define the discrete exterior derivative which we will
call $\d$, as in the usual exterior calculus. The discrete exterior
derivative will be defined as the dual, with respect to the natural
pairing defined above, of the boundary operator, which is defined
below.

\begin{definition}
The \textbf{boundary}\index{boundary} operator
$\partial_k:C_k( K;\mathbb{Z}) \rightarrow
C_{k-1}(K;\mathbb{Z})  $ is a homomorphism defined by
its action on a simplex $\sigma^k = [v_0, \ldots, v_k]$,
\begin{align*}
\partial_k\sigma^k=\partial_k( [ v_0, \ldots
,v_k] )  =\sum_{i=0}^{k}(-1)^i [ v_0,\ldots,\hat{v}_{i},\ldots,v_{k}]\, ,
\end{align*}
where $[v_0, \ldots, \hat{v}_i, \ldots, v_k]$ is the
$(k-1)$-simplex obtained by omitting the vertex $v_i$. Note that
$\partial_k \circ
\partial_{k+1} = 0$.
\end{definition}

\begin{example}\index{boundary!example}
Given an oriented triangle $[v_0, v_1, v_2]$ the boundary, by the
above definition, is the chain $[v_1, v_2] - [v_0, v_2] + [v_0,
v_1]$, which are the three boundary edges of the triangle.
\end{example}

\begin{definition}
On a simplicial complex of dimension $n$, a \textbf{chain complex}\index{chain!complex}\index{complex!chain|see{chain, complex}}
is a collection of chain groups and homomorphisms $\partial_k$,
such that,
\[
\xymatrix{0\quad \ar[r]& \quad C_n(K)\quad
\ar[r]^(0.57){\partial_n} & \quad\ldots\quad
\ar[r]^(0.52){\partial_{k+1}} & \quad C_k (K)\quad
\ar[r]^(0.57){\partial_k} & \quad\ldots\quad
\ar[r]^(0.45){\partial_1} & \quad C_0(K)\quad \ar[r]
 & \quad 0 }\! ,
\]
and $\partial_k \circ \partial_{k+1} = 0$.
\end{definition}

\begin{definition}
The \textbf{coboundary operator},\index{coboundary} $\delta^k:C^k(K) \rightarrow C^{K+1}(K) $, is defined by duality to the boundary operator, with respect to the natural bilinear pairing between discrete forms and chains. Specifically, for a discrete form
$\alpha^k \in \Omega_d^k(K)$, and a chain $c_{k+1} \in
C_{k+1}(K;\ZZ)$, we define $\delta^k$ by
\begin{equation}
\left\langle \delta^k\alpha^{k},c_{k+1}\right\rangle =
    \left\langle \alpha^{k}, \partial_{k+1}c_{k+1}\right\rangle \, .
\label{E:coboundary}
\end{equation}
That is to say
\[
    \delta^k(\alpha^k) = \alpha^k \circ \partial_{k+1} \, .
\]
This definition of the coboundary operator induces the
\textbf{cochain complex}\index{complex!cochain|see{cochain, complex}}\index{cochain!complex},
\[
\xymatrix{0\quad & \quad C^n(K)\quad \ar[l] & \quad\ldots\quad
\ar[l]_(0.35){\delta^{n-1}} & \quad C^k (K)\quad
\ar[l]_(0.53){\delta^k} & \quad\ldots\quad
\ar[l]_(0.35){\delta^{k-1}} & \quad C^0(K)\quad
\ar[l]_(0.53){\delta^0} & \quad 0 \ar[l]}\! ,
\]
where it is easy to see that $\delta^{k+1} \circ \delta^k = 0$.
\end{definition}

\begin{definition}
The \textbf{discrete exterior derivative}\index{exterior derivative} denoted by $\d :
\Omega_d^k(K) \rightarrow \Omega_d^{k+1}(K)$ is defined to be
the coboundary operator $\delta^k$.
\end{definition}

\begin{remark}
With the above definition of the exterior derivative, $\d:
\Omega^k_d(K) \rightarrow \Omega^{k+1}_d(K)$, and the relationship
between the natural pairing and integration, one can regard
equation \ref{E:coboundary} as a discrete {\bfi generalized
Stokes' theorem}\index{Stokes' theorem!generalized}. Thus, given a
$k$-chain $c$, and a discrete $k$-form $\alpha$, the discrete
Stokes' theorem, which is true by definition, states that
\[
    \left\langle \d \alpha, c \right\rangle = \left\langle \alpha,
        \partial c \right \rangle \, .
\]
Furthermore, it also follows immediately that $\d^{k+1}\d^k=0$.
\end{remark}

\paragraph{Dual Discrete Forms.}\index{form!dual}
Everything we have said above in terms of simplices and the
simplicial complex $K$ can be said in terms of the cells that are
duals of simplices and elements of the dual complex $\star K$. One
just has to be a little more careful in the definition of the
boundary operator, and the definition we construct below is well-defined on the
dual cell complex. This gives us the notion of cochains of cells
in the dual complex and these are the {\bfi dual discrete
forms}.

\begin{definition}\label{dec:def:dual_boundary}
The \textbf{dual boundary}\index{boundary!dual} operator,
$\partial_{k}:C_{k}\left( \star K;\mathbb{Z}\right) \rightarrow
C_{k-1}\left(  \star K;\mathbb{Z}\right)  $, is a homomorphism
defined by its action on a dual cell $\hat\sigma^k=\star
\sigma^{n-k} = \star [v_0, \ldots, v_{n-k}]$,
\begin{align*}
\partial \hat \sigma^k &= \partial \star [v_0,...,v_{n-k}]\\
&= \sum_{\sigma^{n-k+1}\succ \sigma^{n-k}} \star\sigma^{n-k+1}\, ,
\end{align*}
where $\sigma^{n-k+1}$ is oriented so that it is consistent with the induced orientation
on $\sigma^{n-k}$.
\end{definition}

\section{Hodge Star and Codifferential} \label{dec:sec:hodge}

In the exterior calculus for smooth manifolds, the Hodge star,
denoted $\ast$, is an isomorphism between the space of $k$-forms
and $(n-k)$-forms. The Hodge star is useful in defining the
adjoint of the exterior derivative and this is adjoint is called
the codifferential. The Hodge star, $\ast:\Omega^k(M)\rightarrow \Omega^{n-k}(M)$, is in the smooth case uniquely defined by the identity,
\[ \langle\!\langle \alpha^k, \beta^k \rangle\!\rangle \mathbf{v} = \alpha^k\wedge\ast\beta^k\, , \]
where $\langle\!\langle\, , \, \rangle\!\rangle$ is a metric on differential forms, and $\mathbf{v}$ is the volume-form. For a more in-depth discussion, see, for example, page 411 of \cite{AbMaRa1988}.

The appearance of $k$ and $(n-k)$ in the definition of Hodge star
may be taken to be a hint that primal and dual meshes will play
some role in the definition of a discrete Hodge star, since the
dual of a $k$-simplex is an $(n-k)$-cell. Indeed, this is the case.

\begin{definition}\label{dec:defn:hodge_star}
The \textbf{discrete Hodge Star}\index{Hodge star} is a map $\ast:\Omega^k_d(K)\rightarrow
\Omega_d^{n-k}(\star K)$, defined by its action on simplices. For a $k$-simplex $\sigma^k$, and a discrete $k$-form $\alpha^k$,
\[
\frac{1}{|\sigma^k|}\langle \alpha^k, \sigma^k
\rangle=\frac{1}{|\star\sigma^k|}\langle
\ast\alpha^k,\star\sigma^k\rangle.
\]
\end{definition}

The idea that the discrete Hodge star maps primal discrete forms
to dual forms, and vice versa, is well-known. See, for example, \cite{SeSeSeAd2000}. However, notice we now make use of the
volume of these primal and dual meshes. But the definition we have
given above does appear in the work of \cite{Hi2002}.

The definition implies that the primal and dual \emph{averages}
must be equal. This idea has already been introduced, not in the
context of exterior calculus, but in an attempt at defining
discrete differential geometry operators, see \cite{MeDeScBa2002}.

\begin{remark}
Although we have defined the discrete Hodge star above, we will show in Remark~\ref{dec:rem:hodge_star_from_inner_product} of \S\ref{dec:sec:commute} that if an appropriate discrete wedge product and metric on discrete $k$-forms is defined, then the expression for the discrete Hodge star operator follows from the smooth definition.
\end{remark}

\begin{lemma} \label{lemma:starstar}
For a $k$-form $\alpha^k$,
\[
    \ast \ast \alpha^k = (-1)^{k(n-k)} \alpha^k \, .
\]
\begin{proof}
The proof is a simple calculation using the property that for a
simplex or a cell $\sigma^k$, $\star\star(\sigma^k) =
(-1)^{k(n-k)} \sigma^k$ (Equation \ref{E:starstarsigma}).
\end{proof}
\end{lemma}

\begin{definition} Given a simplicial or a dual cell complex $K$
the \textbf {discrete codifferential operator},\index{codifferential} $\boldsymbol{\delta} :
\Omega_d^{k+1}(K) \rightarrow \Omega_d^k(K)$, is defined by
$\boldsymbol{\delta}(\Omega_d^0(K)) = 0$ and on discrete $(k+1)$-forms to be
\[
    \boldsymbol{\delta} \beta = (-1)^{nk+1} \ast\d\ast\beta \, .
\]
\end{definition}
With the discrete forms, Hodge star, $\d$ and $\boldsymbol{\delta}$ defined so
far, we already have enough to do an interesting calculation
involving the Laplace--Beltrami operator. But, we will show this
calculation in \S\ref{dec:sec:Divergence} after we have
introduced discrete divergence operator.

\section{Maps between $1$-Forms and Vector Fields}
\label{dec:sec:Maps}

\begin{figure}[htbp]
\begin{center}
\subfigure[Primal vector field]{\includegraphics[scale=0.85]{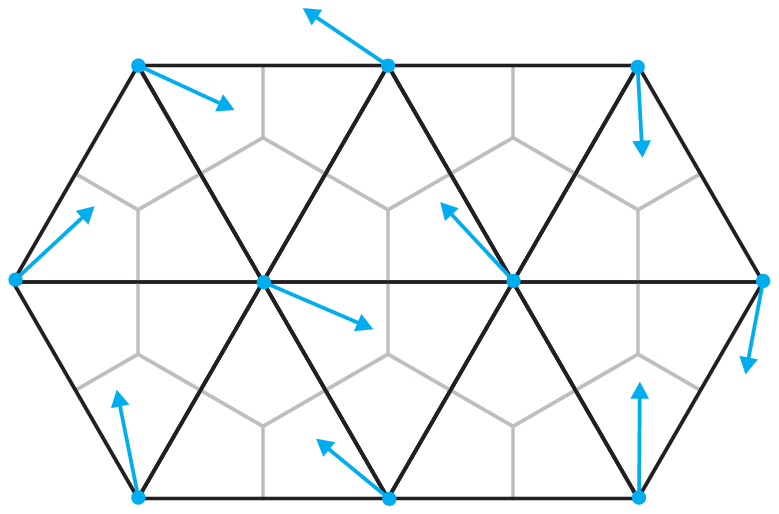}}\qquad
\subfigure[Dual vector field]{\includegraphics[scale=0.85]{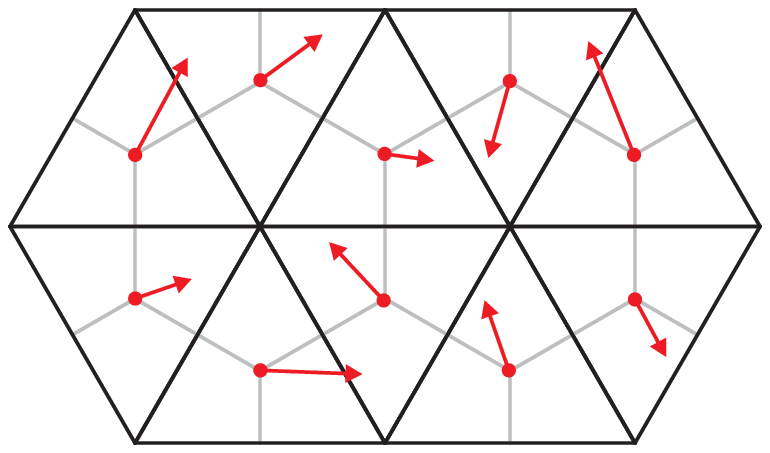}}
\caption{\label{dec:fig:discrete_vf}Discrete vector fields.}
\end{center}
\end{figure}

Just as discrete forms come in two flavors, primal and dual
(being linear functionals on primal chains or chains made up of dual cells), discrete vector fields also come in two flavors.
Before formally defining primal and dual discrete vector fields,
consider the examples illustrated in
Figure~\ref{dec:fig:discrete_vf}.
The distinction lies in the choice of basepoints, be they primal or dual vertices, to which we assign vectors.

\begin{definition}
Let $K$ be a flat simplicial complex, that is, the dimension of $K$ is
the same as that of the embedding space. A \textbf{primal discrete
vector field}\index{vector field!primal} $X$ on a flat simplicial complex
$K$ is a map from the zero-dimensional primal subcomplex $K^{(0)}$
(i.e., the primal vertices) to $\RR^N$.  We will denote the space of
such vector fields by ${\mathfrak X}_d(K)$\index{$X_d(K)$@$\mathfrak{X}_d(K)$|see{vector field, primal}}.  The value of such a
vector field is piecewise constant on the dual $n$-cells of $\star
K$. Thus, we could just as well have called such vector fields dual
and defined them as functions on the $n$-cells of $\star K$.
\end{definition}

\begin{definition}
A \textbf{dual discrete vector field}\index{vector field!dual} $X$ on a
simplicial complex $K$ is a map from the zero-dimensional dual
subcomplex $(\star K)^{(0)}$ (i.e, the circumcenters of the primal
$n$ simplices) to $\RR^N$ such that its value on each dual vertex
is tangential to the corresponding primal $n$-simplex. We will
denote the space of such vector fields by ${\mathfrak X}_d(\star
K)$\index{$X_d(\star K)$@$\mathfrak{X}_d(\star K)$|see{vector field, dual}}. The value of such a vector field is piecewise constant on the
$n$-simplices of $K$. Thus, we could just as well have called such
vector fields primal and defined them as functions on the
$n$-simplices of $K$.
\end{definition}

\begin{remark}
In this paper we have defined the primal vector fields only for
flat meshes. We will address the issue of non-flat meshes in
separate work.
\end{remark}

As in the smooth exterior calculus, we want to define the flat
($\flat$) and sharp ($\sharp$) operators that relate forms to vector
fields. This allows one to write various vector calculus
identities in terms of exterior calculus.

\begin{definition}
Given a simplicial complex $K$ of dimension $n$, the
\textbf{discrete flat operator on a dual vector field}\index{flat},\index{$\flat$|see{flat}} $\flat :
\mathfrak{X}_d(\star K) \rightarrow \Omega^d(K)$, is defined by
its evaluation on a primal 1 simplex $\sigma^1$,
\[
    \langle X^\flat, \sigma^1 \rangle = \sum_{\sigma^n \succ \sigma^1}
    \frac{|\star \sigma^1 \cap \sigma^n|}{|\star \sigma^1|}
    X \cdot \vec{\sigma}^1\, ,
\]
where $X \cdot \vec{\sigma}^1$ is the usual dot product of vectors
in $\RR^N$, and $\vec{\sigma}^1$ stands for the vector
corresponding to $\sigma^1$, and with the same orientation. The sum
is over all $\sigma^n$ containing the edge $\sigma^1$. The volume factors are in dimension $n$.
\end{definition}

\begin{definition}
Given a simplicial complex $K$ of dimension $n$, the
\textbf{discrete sharp operator on a primal $1$-form},\index{sharp}\index{$\sharp$|see{sharp}} $\sharp :
\Omega^d(K)\rightarrow\mathfrak{X}_d(\star K) $, is defined by
its evaluation on a given vertex $v$,
\[
   \alpha^\sharp(v) = \sum_{[v,\sigma^0]}\langle\alpha,[v,\sigma^0]\rangle \sum_{\sigma^n\succ[v,\sigma^0]}\frac{|\star v \cap \sigma^n|}{|\sigma^n|}\hat{n}_{[v,\sigma^0]}\, ,
\]
where the outer sum is over all $1$-simplices containing the vertex $v$, and the inner sum is over all $n$-simplices containing the $1$-simplex $[v,\sigma^0]$. The volume factors are in dimension $n$, and the vector $\hat{n}_{[v,\sigma^0]}$ is the normal vector to the simplex $[v,\sigma^0]$, pointing into the $n$-simplex $\sigma^n$.
\end{definition}

For a discussion of the proliferation of discrete sharp and flat operators that arise from considering the interpolation of differential forms and vector fields, please see \cite{Hirani2003}.
\section{Wedge Product} \label{dec:sec:wedge}
As in the smooth case, the wedge product we will construct is a
way to build higher degree forms from lower degree ones. For
information about the smooth case, see the first few pages of
Chapter 6 of \cite{AbMaRa1988}.

\begin{definition}\label{dec:def:wedge}
Given a primal discrete $k$-form $\alpha^k \in \Omega^k_d(K)$, and
a primal discrete $l$-form $\beta^l \in \Omega^l_d(K)$, the
\textbf{discrete primal-primal wedge product},\index{wedge product}\index{wedge product!primal-primal} $\wedge :
\Omega^k_d(K) \times \Omega^l_d(K) \rightarrow \Omega^{k+l}_d(K)$,
defined by the evaluation on a $(k+l)$-simplex $\sigma^{k+l} =
[v_0, \ldots, v_{k+l}]$ is given by
\[
\langle \alpha^k \wedge \beta^l, \sigma^{k+l} \rangle =
\frac{1}{(k+l)!} \sum_{\tau \in S_{k+l+1}}
\operatorname{sign}(\tau)
    \frac{|\sigma^{k+l} \cap \star v_{\tau(k)}|}{|\sigma^{k+l}|}
    \alpha \smile \beta (\tau(\sigma^{k+l})) \, ,
\]
where $S_{k+l+1}$ is the permutation group, and its elements are
thought of as permutations of the numbers $0, \ldots, k+l+1$. The
notation $\tau(\sigma^{k+l})$ stands for the simplex
$[v_{\tau(0)}, \ldots, v_{\tau(k+l)}]$. Finally, the notation
$\alpha \smile \beta (\tau(\sigma^{k+l}))$ is borrowed from
algebraic topology (see, for example, page 206 of \cite{Ha2001}) and is defined as
\[
  \alpha \smile \beta (\tau(\sigma^{k+l})) :=
    \langle \alpha, \left[v_{\tau(0)},\ldots,v_{\tau(k)}\right] \rangle
    \langle \beta, \left[v_{\tau(k)},\ldots,v_{\tau(k+l)}\right] \rangle\, .
\]
\end{definition}

\begin{example}\index{wedge product!example}
When we take the wedge product of two discrete $1$-forms, we obtain terms in the sum that are graphically represented in Figure~\ref{dec:fig:wedge}.
\begin{figure}[H]
\begin{center}
\includegraphics[scale=0.5]{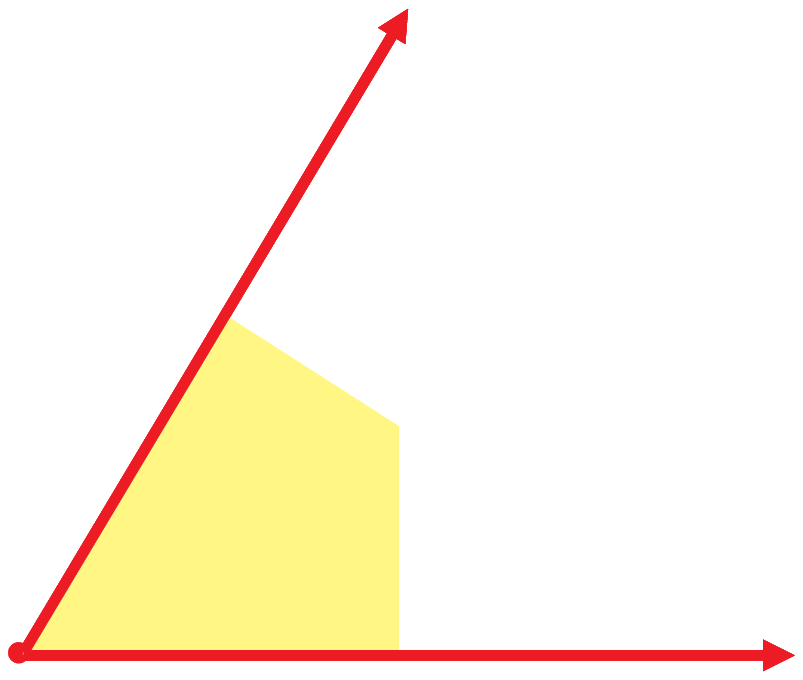}\qquad
\includegraphics[scale=0.5]{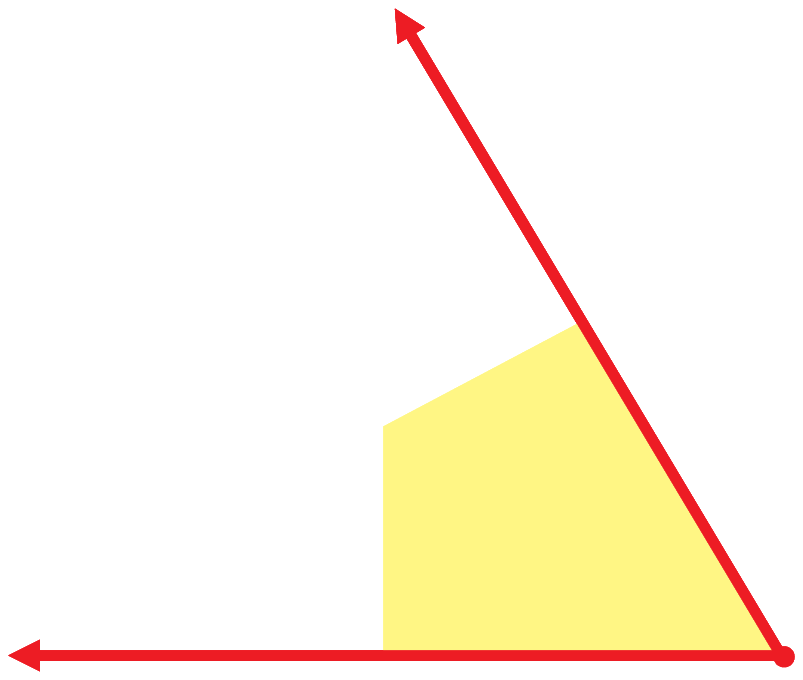}\qquad
\includegraphics[scale=0.5]{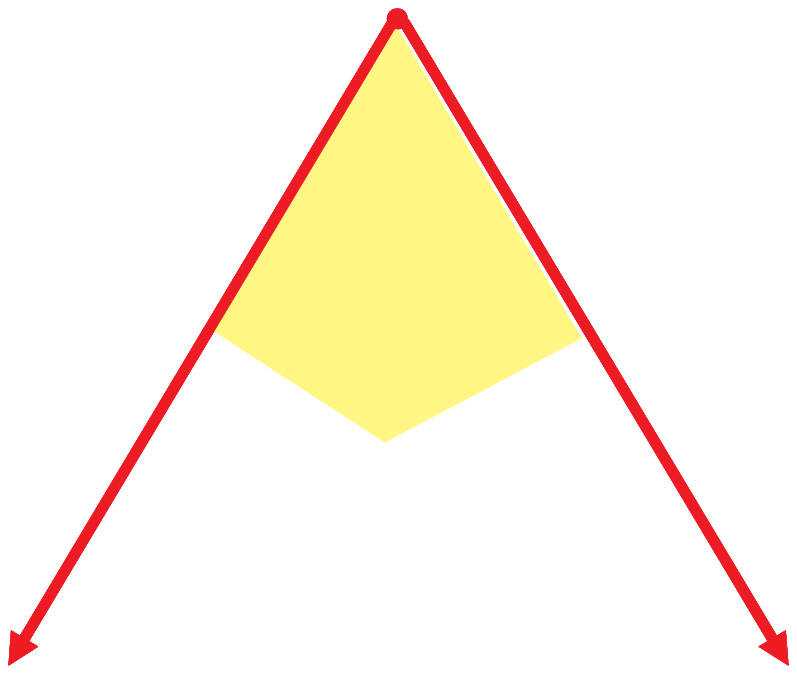}
\end{center}
\caption{\label{dec:fig:wedge}Terms in the wedge product of two discrete $1$-forms.}
\end{figure}
\end{example}

\begin{definition}
 Given a dual discrete $k$-form $\hat\alpha^k
\in \Omega^k_d(\star K)$, and a primal discrete $l$-form $\hat
\beta^l \in \Omega^l_d(\star K)$, the \textbf{discrete dual-dual
wedge product},\index{wedge product! dual-dual} $\wedge : \Omega^k_d(\star K) \times
\Omega^l_d(\star K) \rightarrow \Omega^{k+l}_d(\star K)$, defined
by the evaluation on a $(k+l)$-cell $\hat\sigma^{k+l}=\star
\sigma^{n-k-l}$, is given by
\begin{align*}
\langle \hat\alpha^k\wedge\hat\beta^l,\hat\sigma^{k+l}\rangle =&
\langle \hat\alpha^k\wedge\hat\beta^l,\star\sigma^{n-k-l}\rangle\\
=& \sum_{\sigma^n\succ\sigma^{n-k-l}}
\operatorname{sign}(\sigma^{n-k-l},[v_{k+l},\ldots,v_n])
\sum_{\tau\in S_{k+l}} \operatorname{sign}(\tau)\\
&\qquad\cdot \langle \hat\alpha^k,
\star[v_{\tau(0)},\ldots,v_{\tau(l-1)},v_{k+l},\ldots,v_n]\rangle
\langle
\hat\beta^l,\star[v_{\tau(l)},\ldots,v_{\tau(k+l-1)},v_{k+l},\ldots,v_n]\rangle
\end{align*}
where $\sigma^n=[v_0,\ldots,v_n]$, and, without loss of
generality, assumed that $\sigma^{n-k-l}=\pm [v_{k+l},\ldots,v_n]$.
\end{definition}

\paragraph{Anti-Commutativity of the Wedge Product.}

\begin{lemma}\index{wedge product!anti-commutative}
The discrete wedge product, $\wedge:C^k(K)\times C^l(K)\rightarrow
C^{k+l}(K)$, is anti-commutative, i.e.,
\[\alpha^k\wedge\beta^k=(-1)^{kl}\beta^l\wedge\alpha^k\, .\]
\end{lemma}
\begin{proof}
We first rewrite the expression for the discrete wedge product
using the following computation,
\begin{align*}
&\sum_{\bar\tau\in
S_{k+l+1}}\operatorname{sign}(\bar\tau)|\sigma^{k+l} \cap \star
v_{\bar\tau(k)}|
\langle\alpha^k,\bar\tau[v_0,\ldots,v_k]\rangle \beta^l,\bar\tau[v_k,\ldots,v_{k+l}]\rangle\\
&\qquad=\sum_{\bar\tau\in S_{k+l+1}} (-1)^{k-1}
\operatorname{sign}(\bar\tau) |\sigma^{k+l} \cap \star
v_{\bar\tau(k)}|\
\langle\alpha^k,\bar\tau[v_1,\ldots,v_0,v_k]\rangle\langle\beta^l,\bar\tau[v_k,\ldots,v_{k+l}]\rangle\\
&\qquad=\sum_{\bar\tau\in S_{k+l+1}} (-1)^{k-1}
\operatorname{sign}(\bar\tau) |\sigma^{k+l} \cap \star
v_{\bar\tau\rho(0)}|\langle\alpha^k,\bar\tau\rho[v_1,\ldots,v_k,
v_0]\rangle\langle\beta^l,\bar\tau\rho[v_0,
v_{k+1},\ldots,v_{k+l}]\rangle\\
&\qquad=\sum_{\bar\tau\in S_{k+l+1}} (-1)^{k-1} (-1)^k
\operatorname{sign}(\bar\tau) |\sigma^{k+l} \cap \star
v_{\bar\tau\rho(0)}|
\langle\alpha^k,\bar\tau\rho[v_0,\ldots,v_k]\rangle\langle\beta^l,\bar\tau\rho [v_0,v_{k+1},\ldots,v_{k+l}]\rangle\\
&\qquad=\sum_{\bar\tau\rho\in S_{k+l+1}\rho} (-1)^{k-1} (-1)^k
(-1) \operatorname{sign}(\bar\tau\rho) |\sigma^{k+l} \cap \star
v_{\bar\tau\rho(0)}|\\
&\qquad\qquad\qquad\qquad\cdot
\langle\alpha^k,\bar\tau\rho[v_0,\ldots,v_k]\rangle\langle\beta^l,\bar\tau\rho [v_0,v_{k+1},\ldots,v_{k+l}]\rangle\\
&\qquad=\sum_{\tau\in S_{k+l+1}} \operatorname{sign}(\tau)
|\sigma^{k+l} \cap \star v_{\tau(0)}|
\langle\alpha^k,\tau[v_0,\ldots,v_k]\rangle\langle\beta^l,\tau [v_0,v_{k+1},\ldots,v_{k+l}]\rangle\, .
\end{align*}
Here, we used the elementary fact, from permutation group theory,
that a $k+1$ cycle can be written as the product of $k$
transpositions, which accounts for the $(-1)^k$ factors. Also,
$\rho$ is a transposition of $0$ and $k$. Then, the discrete wedge
product can be rewritten as
\begin{align*}
\langle \alpha^k \wedge \beta^l, \sigma^{k+l} \rangle &=
\frac{1}{(k+l)!} \sum_{\tau \in S_{k+l+1}}
\operatorname{sign}(\tau) \frac{|\sigma^{k+l} \cap \star
v_{\tau(0)}|}{|\sigma^{k+l}|}
\\*
&\qquad\qquad\qquad\qquad\cdot    \langle\alpha^k,[v_{\tau(0)},\ldots,v_{\tau(k)}]\rangle
    \langle\beta^l,[v_{\tau(0),\tau(k+1)},\ldots,v_{\tau(k+l)}]\rangle.
\end{align*}
For ease of notation, we denote $[v_0,\ldots,v_k]$ by $\sigma^k$,
and $[v_0,v_{k+1},\ldots,v_{k+l}]$ by $\sigma^l$. Then, we have
\[
\langle \alpha^k \wedge \beta^l, \sigma^{k+l} \rangle =
\frac{1}{(k+l)!} \sum_{\tau \in S_{k+l+1}}
\operatorname{sign}(\tau) \frac{|\sigma^{k+l} \cap \star
v_{\tau(0)}|}{|\sigma^{k+l}|}
    \langle\alpha^k,\tau(\sigma^k)\rangle
    \langle\beta^l,\tau(\sigma^l)\rangle.
\]
Furthermore, we denote $[v_0,v_{l+1}\ldots,v_{k+l}]$ by
$\bar\sigma^k$, and $[v_0,v_1,\ldots,v_l]$ by $\bar\sigma^l$.
Then,
\begin{align*}
\langle \beta^l \wedge \alpha^k, \sigma^{k+l} \rangle =
\frac{1}{(k+l)!} \sum_{\bar\tau \in S_{k+l+1}}
\operatorname{sign}(\bar\tau) \frac{|\sigma^{k+l} \cap \star
v_{\bar\tau(0)}|}{|\sigma^{k+l}|}
    \langle\alpha^k,\bar\tau(\bar\sigma^k)\rangle
    \langle\beta^l,\bar\tau(\bar\sigma^l)\rangle.
\end{align*}
Consider the permutation $\theta\in S_{k+l+1}$, given by
\[
\theta=\begin{pmatrix}
  0 & 1 & \ldots & k & k+1 & \ldots & k+l \\
  0 & l+1 & \ldots & k+l & 1 & \ldots & l
\end{pmatrix},
\]
which has the property that
\begin{align*}
\bar\sigma^k &= \theta(\sigma^k),\\
\bar\sigma^l &= \theta(\sigma^l).
\end{align*}
Then, we have
\begin{align*}
\langle \beta^l \wedge \alpha^k, \sigma^{k+l} \rangle &=
\frac{1}{(k+l)!} \sum_{\bar\tau \in S_{k+l+1}}
\operatorname{sign}(\bar\tau) \frac{|\sigma^{k+l} \cap \star
v_{\bar\tau(0)}|}{|\sigma^{k+l}|}
    \langle\alpha^k,\bar\tau(\bar\sigma^k)\rangle
    \langle\beta^l,\bar\tau(\bar\sigma^l)\rangle\\
&= \frac{1}{(k+l)!} \sum_{\bar\tau \in S_{k+l+1}}
\operatorname{sign}(\bar\tau) \frac{|\sigma^{k+l} \cap \star
v_{\bar\tau\theta(0)}|}{|\sigma^{k+l}|}
    \langle\alpha^k,\bar\tau\theta(\sigma^k)\rangle
    \langle\beta^l\bar\tau\theta(\sigma^l)\rangle\\
&= \frac{1}{(k+l)!} \sum_{\bar\tau\theta \in S_{k+l+1}\theta}
\operatorname{sign}(\bar\tau\theta)\operatorname{sign}(\theta)
\frac{|\sigma^{k+l} \cap \star
v_{\bar\tau\theta(0)}|}{|\sigma^{k+l}|}
    \langle\alpha^k,\bar\tau\theta(\sigma^k)\rangle
    \langle\beta^l,\bar\tau\theta(\sigma^l)\rangle.
\end{align*}
By making the substitution, $\tau=\bar\tau\theta$, and noting that
$S_{k+l+1}\theta=S_{k+l+1}$, we obtain
\begin{align*}
\langle \beta^l \wedge \alpha^k, \sigma^{k+l} \rangle &=
\operatorname{sign}(\theta)\frac{1}{(k+l)!} \sum_{\tau \in
S_{k+l+1}} \operatorname{sign}(\tau) \frac{|\sigma^{k+l} \cap
\star v_{\tau(0)}|}{|\sigma^{k+l}|}
    \langle\alpha^k,\tau(\sigma^k)\rangle
    \langle\beta^l,\tau(\sigma^l)\rangle\\
&=\operatorname{sign}(\theta)\langle \alpha^k\wedge\beta^l,
\sigma^{k+l} \rangle\, .
\end{align*}
To obtain the desired result, we simply need to compute the sign
of $\theta$, which is given by
\[
\operatorname{sign}(\theta)=(-1)^{kl}.
\]
This follows from the observation that in order to move each of
the last $l$ vertices of $\sigma^{k+l}$ forward, we require $k$
transpositions with $v_1,\ldots,v_k$. Therefore, we obtain
\[
\langle \beta^l \wedge \alpha^k, \sigma^{k+l} \rangle
=\operatorname{sign}(\theta)\langle \alpha^k\wedge\beta^l,
\sigma^{k+l} \rangle = (-1)^{kl}\langle \alpha^k\wedge\beta^l,
\sigma^{k+l} \rangle,
\]
and
\[ \alpha^k \wedge \beta^l = (-1)^{kl} \beta^l \wedge \alpha^k. \qedhere\]
\end{proof}

\paragraph{Leibniz Rule for the Wedge Product.}

\begin{lemma}\index{wedge product!Leibniz rule}
The discrete wedge product satisfies the Leibniz rule,
\[
\d(\alpha^k\wedge\beta^l)=(\d\alpha^k)\wedge\beta^l+(-1)^k\alpha^k\wedge
(\d\beta^l).
\]
\end{lemma}\begin{proof}
The proof of the Leibniz rule for discrete wedge products is
directly analogous to the proof of the coboundary formula for the simplicial
cup product on cochains, which can be found on page 206 of
\cite{Ha2001}. This is because the discrete exterior
derivative is precisely the coboundary operator, and the wedge product is constructed out of weighted sums of cup products.

The cup product satisfies the Leibniz rule for an given partial
ordering of the vertices, and the permutations in the signed sum
in the discrete wedge product correspond to different choices of
partial ordering. We then obtain the Leibniz rule for the discrete
wedge product by applying it term-wise for each choice of
permutation.

Consider
\begin{align*}
\langle (\d\alpha^k)\wedge\beta^l,\sigma^{k+l+1}\rangle &=
\sum_{i=0}^{k+1}(-1)^i \frac{1}{(k+l)!} \sum_{\tau \in S_{k+l+1}}
\operatorname{sign}(\tau) \frac{|\sigma^{k+l} \cap \star
v_{\tau(0)}|}{|\sigma^{k+l}|}\\
&\qquad\cdot
    \langle\alpha^k,[v_{\tau(0)},\ldots,\hat v_i,\ldots, v_{\tau(k+1)}]\rangle
    \langle\beta^l,[v_{\tau(k+1)},\ldots,v_{\tau(k+l+1)}]\rangle,
\end{align*}
and
\begin{align*}
(-1)^k\langle \alpha^k\wedge(\d\beta^l),\sigma^{k+l+1}\rangle &=
(-1)^k\sum_{i=k}^{k+l+1}(-1)^{i-k} \frac{1}{(k+l)!} \sum_{\tau \in
S_{k+l+1}} \operatorname{sign}(\tau) \frac{|\sigma^{k+l} \cap
\star v_{\tau(0)}|}{|\sigma^{k+l}|}\\
&\qquad\cdot
    \langle\alpha^k,[v_{\tau(0)},\ldots, v_{\tau(k)}]\rangle
    \langle\beta^l,[v_{\tau(k)},\ldots,\hat
    v_i,\ldots,v_{\tau(k+l+1)}]\rangle.
\end{align*}
The last set of terms, $i=k+1$, of the first expression cancels
the first set of terms, $i=k$, of the second expression, and what
remains is simply $\langle
\alpha^k\wedge\beta^l,\partial\sigma^{k+l+1}\rangle$. Therefore,
we can conclude that
\[
\langle (\d\alpha^k)\wedge\beta^l,\sigma^{k+l+1}\rangle +
(-1)^k\langle \alpha^k\wedge(\d\beta^l),\sigma^{k+l+1}\rangle =
\langle \alpha^k\wedge\beta^l,\partial\sigma^{k+l+1}\rangle =
\langle \d(\alpha^k\wedge\beta^l),\sigma^{k+l+1}\rangle,
\]
or simply that the Leibniz rule for discrete differential forms
holds,
\[
\d(\alpha^k\wedge\beta^l)=(\d\alpha^k)\wedge\beta^l+(-1)^k\alpha^k\wedge
(\d\beta^l).\qedhere
\]
\end{proof}

\paragraph{Associativity for the Wedge Product.}\index{wedge product!associativity}
The discrete wedge product which we have introduced is not associative in general. This is a consequence of the fact that the stencil for the two possible triple wedge products are not the same. In the expression for $\langle \alpha^k\wedge(\beta^l\wedge\gamma^m), \sigma^{k+l+m} \rangle$, each term in the double summation consists of a geometric factor multiplied by $\langle\alpha^k,\sigma^k\rangle\langle\beta^l,\sigma^l\rangle\langle\gamma^m,\sigma^m\rangle$ for some $k,l,m$ simplices $\sigma^k,\sigma^l,\sigma^m$.

\begin{figure}[htbp]
\begin{center}
\begin{minipage}{0.3\textwidth}{
\begin{center}
\includegraphics{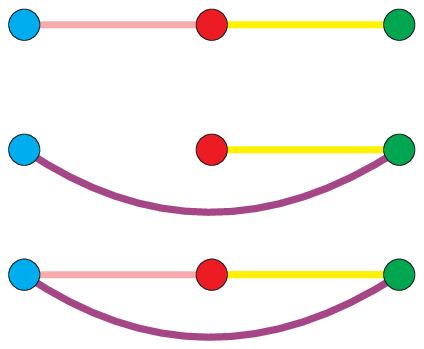}
$\alpha\wedge(\beta\wedge\gamma)$
\end{center}}
\end{minipage}
\qquad
\begin{minipage}{0.3\textwidth}{
\begin{center}
\includegraphics{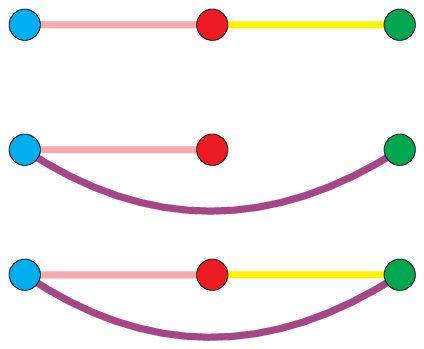}
$(\alpha\wedge\beta)\wedge\gamma$
\end{center}}
\end{minipage}
\end{center}
\caption{\label{dec:fig:associative45}Stencils arising in the double summation for the two triple wedge products.}
\end{figure}

Since $\beta^l$ and $\gamma^m$ are wedged together first, $\sigma^l$ and $\sigma^m$ will always share a common vertex, but $\sigma^k$ could have a vertex in common with only $\sigma^l$, or only $\sigma^m$, or both. We can represent this in a graph, where the nodes denote the three simplices, which are connected by an edge if, and only if, they share a common vertex. The graphical representation of the terms which arise in the two possible triple wedge products are given in Figure~\ref{dec:fig:associative45}.

\begin{figure}[htbp]
\begin{center}
\includegraphics{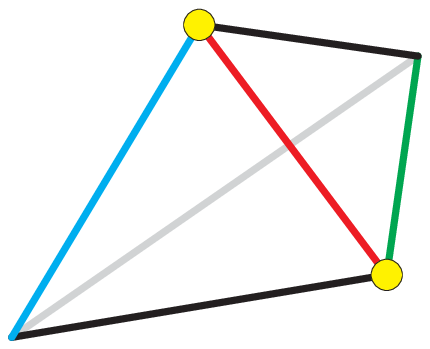}\qquad
\includegraphics{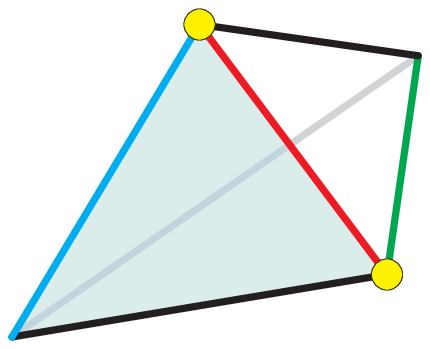}\qquad
\includegraphics{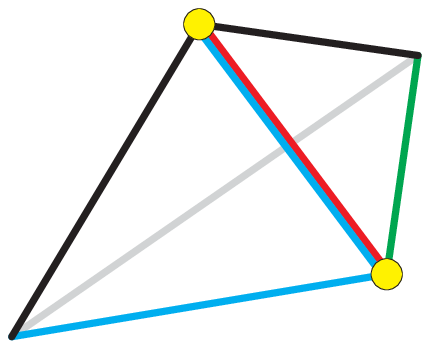}
\end{center}
\caption{\label{dec:fig:associative123}Associativity for closed forms.}
\end{figure}

For the wedge product to be associative for all forms, the two stencils must agree. Since the stencils for the two possible triple wedge products differ, the wedge product is not associative in general. However, in the case of closed forms, we can rewrite the terms in the sum so that all the discrete forms are evaluated on triples of simplices that share a common vertex. This is illustrated graphically in Figure~\ref{dec:fig:associative123}.

This result is proved rigorously in the follow lemma.

\begin{lemma}
The discrete wedge product is associative for closed forms. That
is to say, for $\alpha^k\in C^k(K)$, $\beta^l\in C^l(K)$,
$\gamma^m\in C^m(K)$, such that $\d\alpha^k=0$, $\d\beta^l=0$,
$\d\gamma^m=0$, we have that
\[ (\alpha^k\wedge\beta^l)\wedge\gamma^m =
\alpha^k\wedge(\beta^l\wedge\gamma^m).\]
\end{lemma}
\begin{proof}
\begin{align*}
&\langle
(\alpha^k\wedge\beta^l)\wedge\gamma^m,\sigma^{k+l+m}\rangle\\
&\qquad= \sum_{\tau\in S_{k+l+m+1}} \operatorname{sign}(\tau)
\langle
\alpha^k\wedge\beta^l,\tau[v_0,\ldots,v_{k+l}]\rangle\langle\gamma^m,\tau[v_{k+l},\ldots,v_{k+l+m}]\rangle\\
&\qquad= \sum_{\tau\in S_{k+l+m+1}} \sum_{\rho\in S_{k+l+1}}
\operatorname{sign}(\tau)\operatorname{sign}(\rho) \langle
\alpha^k,\rho\tau[v_0,\ldots,v_k]\rangle\\
&\hspace{1.8in}\cdot\langle\beta^l,\rho\tau[v_k,\ldots,v_{k+l}]\rangle
\langle\gamma^m,\tau[v_{k+l},\ldots,v_{k+l+m}]\rangle
\end{align*}
Here, either $\rho\tau(k)=\tau(k+l)$, in which case all three
permuted simplices share $v_{\tau(k+l)}$ as a common vertex, or we
need to rewrite either $\langle
\alpha^k,\rho\tau[v_0,\ldots,v_k]\rangle$ or
$\langle\beta^l,\rho\tau[v_k,\ldots,v_{k+l}]\rangle$, using the
fact that $\alpha^k$ and $\beta^l$ are closed forms.

If $v_{\tau(k+l)}\notin \rho\tau[v_0,\ldots,v_k]$, then we need to
rewrite $\langle \alpha^k,\rho\tau[v_0,\ldots,v_k]\rangle$ by
considering the simplex obtained by adding the vertex
$v_{\tau(k+l)}$ to $\rho\tau[v_0,\ldots,v_k]$, which is
$[v_{\tau(k+l)},v_{\rho\tau(0)},\ldots,v_{\rho\tau(k)}]$. Then,
since $\alpha^k$ is closed, we have that
\begin{align*}
0 &= \langle \d \alpha^k,
[v_{\tau(k+l)},v_{\rho\tau(0)},\ldots,v_{\rho\tau(k)}]\rangle\\
&= \langle \alpha^k,\partial
[v_{\tau(k+l)},v_{\rho\tau(0)},\ldots,v_{\rho\tau(k)}]\rangle\\
&= \langle\alpha^k,
[v_{\rho\tau(0)},\ldots,v_{\rho\tau(k)}]\rangle - \sum_{i=0}^k
(-1)^i \langle\alpha^k, [v_{\tau(k+l)},v_{\rho\tau(0)},\ldots,\hat
v_{\rho\tau(i)},\ldots,v_{\rho\tau(k)}]\rangle
\end{align*}
or equivalently,
\[
\langle\alpha^k,
[v_{\rho\tau(0)},\ldots,v_{\rho\tau(k)}]\rangle = \sum_{i=0}^k
(-1)^i \langle\alpha^k, [v_{\tau(k+l)},v_{\rho\tau(0)},\ldots,\hat
v_{\rho\tau(i)},\ldots,v_{\rho\tau(k)}]\rangle.
\]
Notice that all the simplices in the sum, with the exception of the last
one, will share two vertices, $v_{\tau(k+l)}$ and
$v_{\rho\tau(k)}$ with $\rho\tau[v_{k},\ldots,v_{k+l}]$, and so
their contribution in the triple wedge product will vanish due to
the anti-symmetrized sum.

Similarly, if $v_{\tau(k+l)}\notin \rho\tau[v_k,\ldots,v_{k+l}]$, using the fact
that $\beta^l$ is closed yields
\[\langle\alpha^k,
[v_{\rho\tau(k)},\ldots,v_{\rho\tau(k+l)}]\rangle =
\sum_{i=k}^{k+l} (-1)^{(i-k)} \langle\alpha^k,
[v_{\tau(k+l)},v_{\rho\tau(k)},\ldots,\hat
v_{\rho\tau(i)},\ldots,v_{\rho\tau(k+l)}]\rangle.\]
As before, all the simplices in the sum, with the exception of the
last one, will share two vertices, $v_{\tau(k+l)}$ and
$v_{\rho\tau(k)}$ with $\rho\tau[v_{0},\ldots,v_{k}]$, and so
their contribution in the triple wedge product will vanish due to
the anti-symmetrized sum.

This allows us to rewrite the triple wedge product in the case of
closed forms as
\begin{align*} \langle
(\alpha^k\wedge\beta^l)\wedge\gamma^m,\sigma^{k+l+m}\rangle =
\sum_{i=0}^{k+l+m}\sum_{\tau\in S_{k+l+m}} &
\operatorname{sign}(\rho_i\tau)\langle \alpha^k,\rho_i\tau
[v_0,\ldots,v_k]\rangle\langle \beta^l,\rho_i\tau
[v_0,v_{k+1},\ldots,v_{k+l}]\rangle\\
&\quad\cdot \langle \gamma^m, \rho_i\tau
[v_0,v_{k+l+1},\ldots,v_{k+l+m}]\rangle\, ,
\end{align*}
where $\tau\in S_{k+l+m}$ is thought of as acting on the set
$\{1,\ldots,k+l+m\}$, and $\rho_i$ is a transposition of $0$ and
$i$. A similar argument allows us to write
$\alpha^k\wedge(\beta^l\wedge\gamma^m)$ in the same form, and
therefore, the wedge product is associative for closed forms.
\end{proof}
\begin{remark}\index{wedge product!convergence}
This lemma is significant, since if we think of a constant
smooth differential form, and discretize it to obtain a
discrete differential form, this discrete form will be closed. As
such, this lemma states that in the infinitesimal limit, the
discrete wedge product we have defined will be associative.

In practice, if we have a mesh with characteristic length $\Delta
x$, then we will have that
\[ \frac{1}{|\sigma^{k+l+m}|}\langle \alpha^k\wedge(\beta^l\wedge\gamma^m) -
(\alpha^k\wedge\beta^l)\wedge\gamma^m,\sigma^{k+l+m}\rangle =
\mathcal{O}(\Delta x),\] which is to say that the average of the
associativity defect is of the order of the mesh size, and
therefore vanishes in the infinitesimal limit.
\end{remark}

\section{Divergence and Laplace--Beltrami}
\label{dec:sec:Divergence}

In this section, we will illustrate the application of some of the DEC operations we have previously defined to the construction of new discrete operators such as the divergence and Laplace--Beltrami operators.

\paragraph{Divergence.}\index{divergence}\index{$\Div$|see{divergence}}
The divergence of a vector field is given in terms of the Lie derivative of the volume-form, by the expression, $(\Div(\X)\mu = \pounds_\X \mu$. Physically, this corresponds to the net flow per unit volume of an infinitesimal volume about a point.

We will define the discrete divergence by using the
formulas defining them in the smooth exterior calculus. The
divergence definition will be valid for arbitrary dimensions. The resulting expressions involve operators that
we have already defined and so we can actually perform some
calculations to express these quantities in terms of geometric
quantities. We will show that the resulting expression in terms of
geometric quantities is the same as that derived by variational
means in \cite{ToLoHiDe2003}.

\begin{definition}
For a discrete dual vector field $X$ the divergence $\Div(X)$ is
defined to be
\[
    \Div(\X) = -\boldsymbol{\delta} X^\flat \, .
\]
\end{definition}

\begin{remark}
The above definition is a theorem in smooth exterior
calculus. See, for example, page 458 of \cite{AbMaRa1988}.
\end{remark}

\begin{figure}[H]
\begin{center}
\includegraphics[scale=0.75]{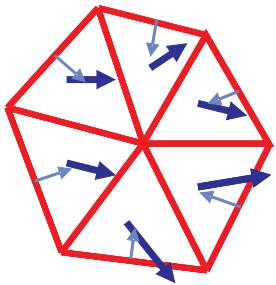}
\end{center}
\caption{\label{dec:fig:divergence}Divergence of a discrete dual vector field.}
\end{figure}

As an example, we will now compute the divergence of a discrete
dual vector field on a two-dimensional simplicial complex $K$, as illustrated in Figure~\ref{dec:fig:divergence}.

A similar derivation works in higher dimensions, where one needs to be mindful of the sign that arises from applying the Hodge star twice, $\ast\ast\alpha^k=(-1)^{k(n-k)}\alpha^k$. Since $\Div(X) =
-\boldsymbol{\delta} X^\flat$, it follows that $\Div(X) = \ast \d \ast X^\flat$.
Since this is a primal $0$-form it can be evaluated on a $0$-simplex
$\sigma^0$, and we have that
\[ \langle \Div(x),\sigma^0\rangle = \langle \ast \d \ast X^\flat, \sigma^0 \rangle\, .\]
Using the definition of discrete Hodge star, and
the discrete generalized Stokes' theorem, we get
\begin{align*}
    \frac{1}{|\sigma^0|} \langle \Div(X) , \sigma^0 \rangle & =
        \frac{1}{|\star \sigma^0|} \langle \ast \ast \d \ast X^\flat,
            \star \sigma^0 \rangle \\
        &= \frac{1}{|\star \sigma^0|} \langle \d \ast X^\flat,
            \star \sigma^0 \rangle \\
        &= \frac{1}{|\star \sigma^0|} \langle \ast X^\flat,
            \partial(\star \sigma^0) \rangle \, .
\end{align*}
The second equality is obtained by applying the definition of the Hodge star,
and the last equality is obtained by applying the discrete generalized Stokes'
theorem. But,
\[
    \partial(\star \sigma^0) = \sum_{\sigma^1 \succ \sigma^0}
    \star \sigma^1\, ,
\]
as given by the expression for the boundary of a dual cell in Equation~\ref{dec:def:dual_boundary}. Thus,
\begin{align*}
    \frac{1}{|\sigma^0|} \langle \Div(X) , \sigma^0 \rangle & =
        \frac{1}{|\star \sigma^0|} \langle \ast X^\flat,
        \sum_{\sigma^1 \succ \sigma^0}  \star \sigma^1 \rangle \\
    &= \frac{1}{|\star \sigma^0|} \sum_{\sigma^1 \succ \sigma^0}
        \langle \ast X^\flat,\star \sigma^1 \rangle \\
    &= \frac{1}{|\star \sigma^0|}\sum_{\sigma^1 \succ \sigma^0}
    \frac{|\star \sigma^1|}{|\sigma^1|} \langle X^\flat,
        \sigma^1 \rangle \\
    &= \frac{1}{|\star \sigma^0|}\sum_{\sigma^1 \succ \sigma^0}
    \frac{|\star \sigma^1|}{|\sigma^1|}
    \sum_{\sigma^2 \succ \sigma^1}
    \frac{|\star \sigma^1 \cap \sigma^2|}{|\star \sigma^1|} X \cdot
        \vec{\sigma}^1 \\
    &= \frac{1}{|\star \sigma^0|}\sum_{\sigma^1 \succ \sigma^0}
    \sum_{\sigma^2 \succ \sigma^1}
    \frac{|\star \sigma^1 \cap \sigma^2|}{|\sigma^1|}
    X \cdot\vec{\sigma}^1 \\
    &=  \frac{1}{|\star \sigma^0|}\sum_{\sigma^1 \succ \sigma^0}
    |\star \sigma^1\cap \sigma^2| \; (X \cdot \frac{\vec{\sigma}^1}{|\sigma^1|}) \, .
\end{align*}

This expression has the nice property that the divergence theorem holds on any dual $n$-chain, which, as a set, is a simply connected subset of $|K|$. Furthermore, the coefficients we computed for the discrete divergence operator are the unique ones for which a discrete divergence theorem holds.

\paragraph{Laplace--Beltrami.}\index{Laplace--Beltrami}\index{$\Delta$|see{Laplace--Beltrami}} The
Laplace--Beltrami operator is the generalization of the Laplacian
to curved spaces. In the smooth case the Laplace--Beltrami
operator on smooth functions is defined to be $\nabla^2 = \Div
\circ \curl = \boldsymbol{\delta} d$. See, for example, page 459
of \cite{AbMaRa1988}. Thus, in the smooth case, the
Laplace--Beltrami on functions is a special case of the more
general Laplace--deRham operator\index{Laplace--deRham}, $\Delta : \Omega^k(M) \rightarrow
\Omega^k(M)$, defined by $\Delta = \d \boldsymbol{\delta} +
\boldsymbol{\delta} \d$.

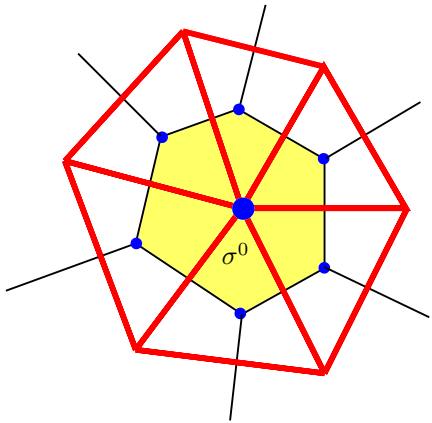
\begin{figure}[htbp]
\WARMprocessMoEPS{laplacebeltrami_new}{eps}{bb}
\renewcommand{\xyWARMinclude}[1]{\scaledfig{.7}{#1}}
\begin{center}
\leavevmode
\begin{xy}
\xyMarkedImport{}
\xyMarkedMathPoints{1}
\end{xy}
\end{center}
\caption{\label{dec:fig:laplace_beltrami}Laplace--Beltrami of a discrete function.}
\renewcommand{\xyWARMinclude}[1]{\includegraphics{#1}}
\end{figure}

As an example, we compute $\Delta f$ on a primal vertex
$\sigma^0$, where $f \in \Omega^0_d(K)$, and $K$ is a (not
necessarily flat) triangle mesh in $\RR^3$, as illustrated in Figure~\ref{dec:fig:laplace_beltrami}.

This calculation is
done below.
\begin{align*}
    \frac{1}{|\sigma^0|}\langle \Delta f, \sigma^0 \rangle & =
     \langle \boldsymbol{\delta} \d f, \sigma^0 \rangle \\
&= - \langle \ast\d \ast\d f, \sigma^0 \rangle \\
&= -\frac{1}{|\star \sigma^0|}\langle\d\ast\d f, \star \sigma^0 \rangle \\
&= -\frac{1}{|\star \sigma^0|} \langle \ast\d f, \partial(\star
\sigma^0)
        \rangle \\
&= -\frac{1}{|\star \sigma^0|} \langle \ast\d f,
    \sum_{\sigma^1 \succ \sigma^0} \star \sigma^1 \rangle \\
&= -\frac{1}{|\star \sigma^0|} \sum_{\sigma^1 \succ \sigma^0}
    \langle \ast\d f, \star \sigma^1 \rangle \\
&= -\frac{1}{|\star \sigma^0|} \sum_{\sigma^1 \succ \sigma^0}
    \frac{|\star \sigma^1|}{|\sigma^1|} \langle \d f, \sigma^1 \rangle \\
&=  -\frac{1}{|\star \sigma^0|} \sum_{\sigma^1 \succ \sigma^0}
    \frac{|\star \sigma^1|}{|\sigma^1|} (f(v) - f(\sigma^0))\, ,
\end{align*}
where $\partial \sigma^1 = v - \sigma^0$. But, the above is the
same as the formula involving cotangents found by \cite{MeDeScBa2002} without using discrete exterior calculus.

Another interesting aspect, which will be discussed in \S\ref{dec:sec:commute}, is that the characterization of harmonic functions as those functions which vanish when the Laplace--Beltrami operator is applied is equivalent to that obtained from a discrete variational principle using DEC as the means of discretizing the Lagrangian.

\section{Contraction and Lie Derivative} \label{dec:sec:Contraction}
In this section we will discuss some more operators that involve
vector fields, namely contraction, and Lie derivatives.

For contraction, we will first define the usual smooth contraction
algebraically, by relating it to Hodge star and wedge products.
This yields one potential approach to defining discrete
contraction. However, since in the discrete theory we are only
concerned with integrals of forms, we can use the interesting
notion of extrusion of a manifold by the flow of a vector field to
define the integral of a contracted discrete differential form.

We learned about this definition
of contraction via extrusion from \cite{Bo2002b}, who goes
on to define discrete extrusion in his paper. Thus, he is able to
obtain a definition of discrete contraction. Extrusion turns out
to be a very nice way to define integrals of operators involving
vector fields, and we will show how to define integrals of Lie
derivatives via extrusion, which will yield discrete Lie
derivatives.

\begin{definition}
Given a manifold $M$, and $S$, a $k$-dimensional submanifold of $M$,
and a vector field $X \in \mathfrak{X}(M)$, we call the manifold
obtained by sweeping $S$ along the flow of $X$ for time $t$ as the
\textbf{extrusion}\index{extrusion} of $S$ by $X$ for time $t$, and denote it by
$E_X^t(S)$. The manifold $S$ carried by the \textbf{flow}\index{flow} for time $t$ will
be denoted $\varphi_X^t(S)$.
\end{definition}

\begin{example}\index{extrusion!example}
Figure~\ref{dec:fig:extrusion} illustrates the $2$-simplex that arises from the extrusion of a $1$-simplex by a discrete vector field that is interpolated using a linear shape function.
\begin{figure}[htbp]
\begin{center}
\includegraphics[scale=0.5]{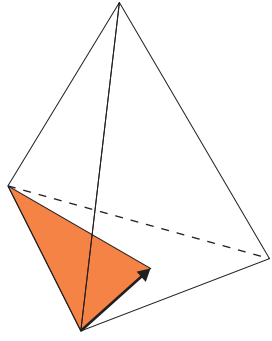}
\end{center}
\caption{\label{dec:fig:extrusion}Extrusion of $1$-simplex by a discrete vector field.}
\end{figure}
\end{example}

\paragraph{Contraction (Extrusion).}\index{contraction!extrusion}\index{extrusion!contraction|see{contraction, extrusion}} We first establish an integral property of the contraction operator.

\begin{lemma}
\[
    \int_S \mathbf{i}_X \beta = \left.\frac{d}{dt}\right|_{t=0}
        \int_{E_X^t(S)} \beta
\]
\begin{proof}
Prove instead that
\[
    \int_0^t \left[ \int_{S_\tau} \mathbf{i}_X \beta \right] d\tau=
        \int_{E_X^t(S)} \beta \, .
\]
Then, by first fundamental theorem of calculus, the desired result
will follow. To prove the above, simply take coordinates on $S$ and
carry them along with the flow and define the transversal
coordinate to be the flow of $X$. This proof is sketched in \cite{Bo2002b}.
\end{proof}
\end{lemma}

This lemma allows us to interpret contraction as being the dual, under the integration pairing between $k$-forms and $k$-volumes, to the geometric operation of extrusion. The discrete contraction operator is then given by
\[ \langle \mathbf{i}_X\alpha^{k+1}, \sigma^k \rangle = \left.\frac{d}{dt}\right|_{t=0}\langle \alpha^{k+1}, E_X^t(\sigma^k) \rangle,\]
where the evaluation of the RHS will typically require that the discrete differential form and the discrete vector field are appropriately interpolated.

\begin{remark}
Since the dynamic definition of the contraction operator only depends on the derivative of pairing of the differential form with the extruded region, it will only depend on the vector field in the region $S$, and not on its extension into the rest of the domain.

In addition, if the interpolation for the discrete vector field satisfies a superposition principle, then the discrete contraction operator will satisfy a corresponding superposition principle.
\end{remark}

\paragraph{Contraction (Algebraic).}\index{contraction!algebraic} Contraction is an operator
that allows one to combine vector fields and forms. For a smooth
manifold $M$, the contraction of a vector field $X \in
\mathfrak{X}(M)$ with a $(k+1)$-form $\alpha \in \Omega^{k+1}(M)$
is written as $\mathbf{i}_X \alpha$, and for vector fields $X_1,
\ldots, X_k \in \mathfrak{X}(M)$, the contraction in smooth
exterior calculus is defined by
\[
    \mathbf{i}_X \alpha(X_1,\ldots,X_k) =
    \alpha(X,X_1,\ldots,X_k) \, .
\]
We define contraction by using an identity that is true in smooth
exterior calculus. This identity originally appeared in \cite{Hirani2003}, and we state
it here with proof.

\begin{lemma}[\cite{Hirani2003}]\label{lemma:contraction}
Given a smooth manifold $M$ of dimension $n$, a vector field
$X\in \mathfrak{X}(M)$, and a $k$-form $\alpha \in \Omega^k(M)$, we
have that
\[
 \mathbf{i}_X \alpha = (-1)^{k(n-k)}\ast(\ast \alpha \wedge X^\flat) \, .
\]
\begin{proof}
Recall that for a smooth function $f\in \Omega^0(M)$, we have that
$\mathbf{i} _X \alpha = f \mathbf{i}_X \alpha$. This, and the
multilinearity of $\alpha$, implies that it is enough to show the
result in terms of basis elements. In particular, let $\tau \in
S_n$ be a permutation of the numbers $1, \ldots n$, such that
$\tau(1) < \ldots < \tau(k)$, and  $\tau(k+1) < \ldots < \tau(n)$.
Let $X = e_{\tau(j)}$, for some $j \in {1, \ldots, n}$. Then, we
have to show that
\[
    \mathbf{i}_{e_{\tau(j)}} e^{\tau(1)} \wedge \ldots \wedge e^{\tau(k)}
    = (-1)^{k(n-k)}
    \ast(\ast (e^{\tau(1)} \wedge \ldots \wedge e^{\tau(k)})
    \wedge e^{\tau(j)}) \, .\label{dec:eqn:algebraic_contraction}
\]
It is easy to see that the LHS is 0 if $j > k$, and it is
\[
    (-1)^{j-1} (e^{\tau(1)} \wedge \ldots \wedge \widehat{e^{\tau(j)}}
    \ldots \wedge e^{\sigma(k)})\, ,
\]
otherwise, where $\widehat{e^{\tau(j)}}$ means that $e^{\tau(j)}$ is omitted from the wedge product. Now, on the RHS of Equation~\ref{dec:eqn:algebraic_contraction}, we have that
\[
    \ast(e^{\tau(1)} \wedge \ldots \wedge e^{\tau(k)}) =
    \operatorname{sign}(\tau)
    (e^{\tau(k+1)}\wedge \ldots \wedge e^{\tau(n)}) \, .
\]
Thus, the RHS is equal to
\[
(-1)^{k(n-k)} \operatorname{sign}(\tau) \ast(e^{\tau(k+1)}\wedge
\ldots \wedge e^{\tau(n)} \wedge e^{\tau(j)}) \, ,
\]
which is 0 as required if $j > k$. So, assume that $1 \le j \le k$. We need to compute
\[
\ast (e^{\tau(k+1)}\wedge\ldots\wedge e^{\tau(n)}\wedge e^{\tau(j)})\, ,
\]
which is given by
\[
s\, e^{\tau(1)}\wedge\ldots\widehat{e^{\tau(j)}}\wedge\ldots\wedge e^{\tau(k)}\, ,
\]
where the sign $s=\pm 1$, such that the equation,
\[
s\, e^{\tau(k+1)} \wedge \ldots \wedge e^{\tau(n)}\wedge e^{\tau(j)}\wedge e^{\tau(1)}\wedge \ldots \wedge \widehat{e^{\tau(j)}}\wedge \ldots \wedge e^{\tau(k)}=\mu \, ,
\]
holds for the standard volume-form, $\mu=e^1\wedge\ldots\wedge e^n$. This implies that
\[s=(-1)^{j-1}(-1)^{k(n-k)}\operatorname{sign}(\tau)\, .\]
Then, RHS = LHS as required.
\end{proof}
\end{lemma}

Since we have expressions for the discrete Hodge star ($\ast$), wedge product ($\wedge$), and flat ($\flat$), we have the necessary ingredients to use the algebraic expression proved in the above lemma to construct a discrete contraction operator.

One has to note, however, that the wedge product is only associative for closed forms, and as a consequence, the Leibniz rule for the resulting contraction operator will only hold for closed forms as well. This is, however, sufficient to establish that the Leibniz rule for the discrete contraction will hold in the limit as the mesh is refined.

\paragraph{Lie Derivative (Extrusion).}\index{Lie derivative!extrusion}\index{extrusion!Lie derivative|see{Lie derivative, extrusion}} As was the case with contraction, we will establish a integral identity that allows the Lie derivative to be interpreted as the dual of a geometric operation on a volume. This involves the flow of a volume by a vector field, and it is illustrated in the following example.
\begin{example}\index{flow!example}
Figure~\ref{dec:fig:flow} illustrates the flow of a $1$-simplex by a discrete vector field interpolated using a linear shape function.
\begin{figure}[htbp]
\begin{center}
\includegraphics[scale=0.5]{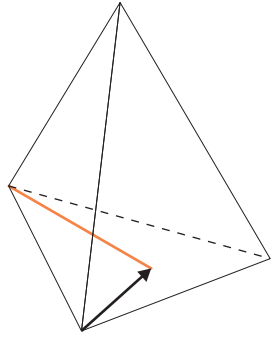}
\end{center}
\caption{\label{dec:fig:flow}Flow of a $1$-simplex by a discrete vector field.}
\end{figure}
\end{example}
\begin{lemma}
\[
\int_S \pounds_X\beta = \left.\frac{d}{dt}\right|_{t=0} \int_{\varphi_X^t(S)}
\beta\, .
\]
\begin{proof}
\begin{align*}
    F_t^\ast(\pounds_X\beta) &= \frac{d}{dt} F_t^\ast \beta \\
    \int_0^t F_\tau^\ast(\pounds_X \beta) d\tau &= F_t^\ast \beta - \beta\\
    \int_S \int_0^t F_\tau^\ast(\pounds_X\beta)d\tau &=
        \int_S F_t^\ast\beta - \int_S \beta \\
    \int_0^t \int_{\varphi_X^\tau(S)}\pounds_X\beta d\tau &= \int_{\varphi_X^t(S)} \beta -
        \int_S\beta \, .\qedhere
\end{align*}
\end{proof}
\end{lemma}
This lemma allows us to define a discrete Lie derivative as follows,
\[ \langle \pounds_X \beta^k, \sigma^k \rangle = \left.\frac{d}{dt}\right|_{t=0} \langle \beta^k,\varphi_X^t(\sigma^k) \rangle\, ,\]
where, as before, evaluating the RHS will require the discrete differential form and discrete vector field to be appropriately interpolated.

\paragraph{Lie Derivative (Algebraic).}\index{Lie derivative!algebraic} Alternatively, as we have expressions for the discrete contraction operator ($\mathbf{i}_X$), and exterior derivative ($\d$), we can construct a
discrete Lie derivative using the Cartan magic formula,
\[\pounds_X \omega =
\mathbf{i}_X\d\omega+\d\mathbf{i}_X\omega.\]
As is the case with the algebraic definition of the discrete contraction, the discrete Lie derivative will only satisfy a Leibniz rule for closed forms. As before, this is sufficient to establish that the Leibniz rule will hold in the limit as the mesh is refined.

\section{Discrete Poincar\'e Lemma}\index{Poincar\'e lemma}
In this section, we will prove the discrete Poincar\'e lemma
by constructing a homotopy operator though a generalized cocone
construction. This section is based on the work in \cite{DeHiLeMa2003b}.

The standard cocone construction fails at the discrete level, since
the cone of a simplex is not, in general, expressible as a chain
in the simplicial complex. As such, the standard cocone does not
necessarily map $k$-cochains to $(k-1)$-cochains.

\begin{figure}[H]
\begin{center}
\includegraphics[scale=0.65]{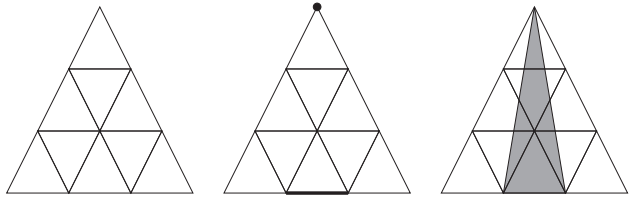}
\end{center}
\caption{\label{dec:fig:nocone}The cone of a simplex is, in general, not expressible as a chain.}
\end{figure}

An example of how the standard cone construction fails to map
chains to chains is illustrated in Figure~\ref{dec:fig:nocone}. Given the
simplicial complex on the left, consisting of triangles, edges and
nodes, we wish, in the center figure, to consider the cone of the
bold edge with respect to the top most node. Clearly, the
resulting cone in the right figure, which is shaded grey, cannot
be expressed as a combination of the triangles in the original
complex.

In this subsection, a generalized cone operator that is valid for
chains is developed which has the essential homotopy properties to
yield a discrete analogue of the Poincar\'e lemma.

We will first consider the case of
trivially star-shaped complexes, followed by logically star-shaped
complexes, before generalizing the result to contractible
complexes.
\begin{definition}
Given a $k$-simplex $\sigma^k=[v_0,\ldots,v_k]$ we construct the
\textbf{cone}\index{cone}\index{cone!geometric} with vertex $w$ and base $\sigma^k$, as follows,
\[w\diamond\sigma^k=[w,v_0,\ldots,v_k].\]
\end{definition}
\begin{lemma}
The geometric cone operator satisfies the following property,
\[\partial(w\diamond\sigma^k)+w\diamond(\partial
\sigma^k)=\sigma^k.\]
\end{lemma}
\begin{proof}
This is a standard result from simplicial algebraic topology.
\end{proof}
\paragraph{Trivially Star-Shaped Complexes.}
\begin{definition}
A complex $K$ is called \textbf{trivially star-shaped}\index{star-shaped!trivially} if there
exists a vertex $w\in K^{(0)}$, such that for all $\sigma^k\in K$,
the cone with vertex $w$ and base $\sigma^k$ is expressible as a
chain in $K$. That is to say, \[ \exists w\in K^{(0)}\mid \forall
\sigma^k\in K, w\diamond\sigma^k\in C_{k+1}(K).\] We can then
denote the cone operation with respect to $w$ as
$p:C_k(K)\rightarrow C_{k+1}(K)$.
\end{definition}
\begin{lemma}
In trivially star-shaped complexes, the cone operator,
$p:C_k(K)\rightarrow C_{k+1}(K)$, satisfies the following identity,
\[
p\partial+\partial p = I,
\]
at the level of chains.
\end{lemma}
\begin{proof}
Follows immediately from the identity for cones, and noting that
the cone is well-defined at the level of chains on trivially
star-shaped complexes.
\end{proof}
\begin{definition}
The \textbf{cocone}\index{cocone} operator, $H:C^k(K)\rightarrow C^{k-1}(K)$, is
defined by
\[
\langle H \alpha^k, \sigma^{k-1} \rangle = \langle \alpha^k,
p(\sigma^{k-1}) \rangle.
\]
This operator is well-defined on trivially star-shaped simplicial
complexes.
\end{definition}
\begin{lemma}
The cocone operator, $H:C^k(K)\rightarrow C^{k-1}(K)$, satisfies the
following identity,
\[ H \d + \d H = I,\]
at the level of cochains.
\end{lemma}
\begin{proof}
A simple duality argument applied to the cone identity,
\[ p \partial + \partial p = I,\]
yields the following,
\begin{align*}
\langle \alpha^k,\sigma^k \rangle &=\langle \alpha^k, (p\partial +
\partial p)\sigma^k\rangle\\
&=\langle \alpha^k, p\partial \sigma^k \rangle + \langle \alpha^k,
\partial p \sigma^k \rangle\\
&=\langle H \alpha^k, \partial \sigma^k \rangle + \langle
\d\alpha^k, p \sigma^k \rangle\\
&=\langle (\d H \alpha^k, \sigma^k \rangle+\langle H \d\alpha^k,
\sigma^k \rangle\\
&=\langle (\d H + H \d)\alpha^k, \sigma^k \rangle.
\end{align*}
Therefore,
\[ H \d + \d H = I,\]
at the level of cochains.
\end{proof}
\begin{corollary}[Discrete Poincar\'{e} Lemma for Trivially
Star-shaped Complexes] Given a closed cochain $\alpha^k$, that is
to say, $\d\alpha^k=0$, there exists a cochain $\beta^{k-1}$, such
that, $\d\beta^{k-1}=\alpha^k$.
\end{corollary}
\begin{proof}
Applying the identity for cochains,
\[ H\d + \d H=I,\]
we have,
\begin{align*}
\langle \alpha^k, \sigma^k \rangle &= \langle (H\d+\d H)\alpha^k,
\sigma^k \rangle\, ,\\
\intertext{but, $\d\alpha^k=0$, so,}
\langle \alpha^k, \sigma^k \rangle &=\langle \d (H\alpha^k),
\sigma^k \rangle.
\end{align*}
Therefore, $\beta^{k-1}=H\alpha^k$ is such that
$\d\beta^{k-1}=\alpha^k$ at the level of cochains.
\end{proof}
\begin{example}
We demonstrate the construction of the tetrahedralization of the
cone of a $(n-1)$-simplex over the origin.

If we denote by $v_i^k$, the projection of the $v_i$ vertex to the
$k$-th concentric sphere, where the $0$-th concentric sphere is simply
the central point, then we fill up the cone $[c,v_1,...v_n]$ with
simplices as follows,
\[ [v_1^0,v_1^1,\ldots, v_n^1], [v_1^2, v_1^1,\ldots, v_n^1],
[v_1^2, v_2^2, v_2^1,\ldots, v_n^1],\ldots, [v_1^2,\ldots, v_n^2,
v_n^1].
\]
Since $S^{n-1}$ is orientable, we can use a consistent
triangulation of $S^{n-1}$ and these $n$-cones to consistently
triangulate $B^n$ such that the resulting triangulation is
star-shaped.

This fills up the region to the 1st concentric sphere, and we
repeat the process by leapfrogging at the last vertex to add
$[v_1^2, ..., v_n^2, v_n^3]$, and continuing the construction, to
fill up the annulus between the 1st and 2nd concentric sphere.
Thus, we can keep adding concentric shells to create an
arbitrarily dense triangulation of a $n$-ball about the origin.

In three dimensions, these simplices are given by
\begin{center}
\begin{minipage}{0.22\textwidth}
\begin{center}
\includegraphics[scale=0.35]{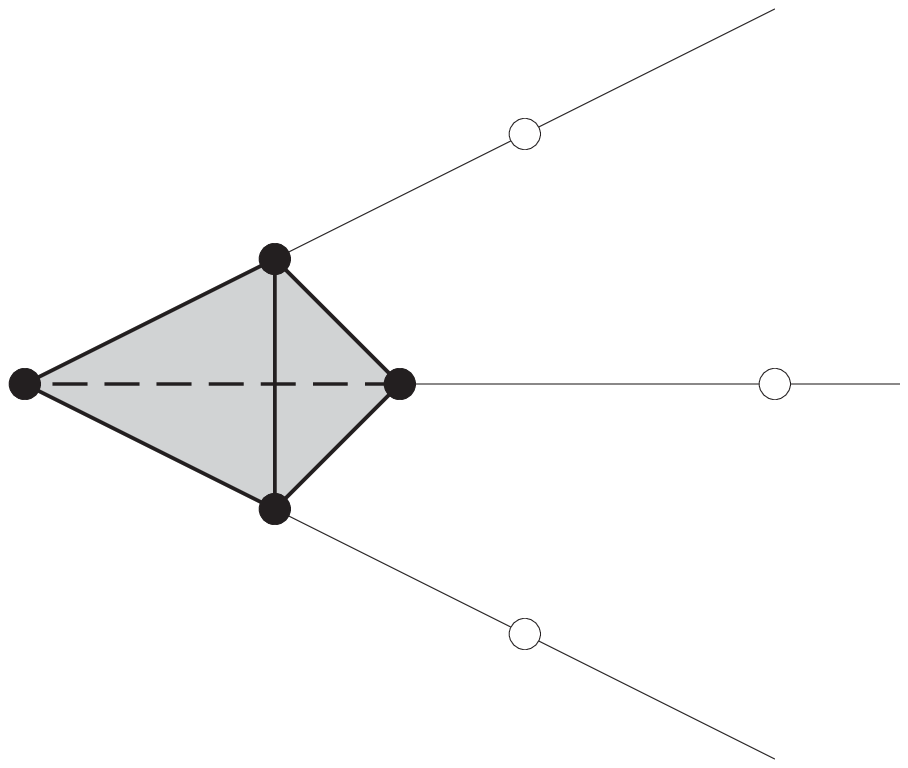}\\
$[c,v_1^1, v_2^1, v_3^1],$
\end{center}
\end{minipage}
\begin{minipage}{0.22\textwidth}
\begin{center}
\includegraphics[scale=0.35]{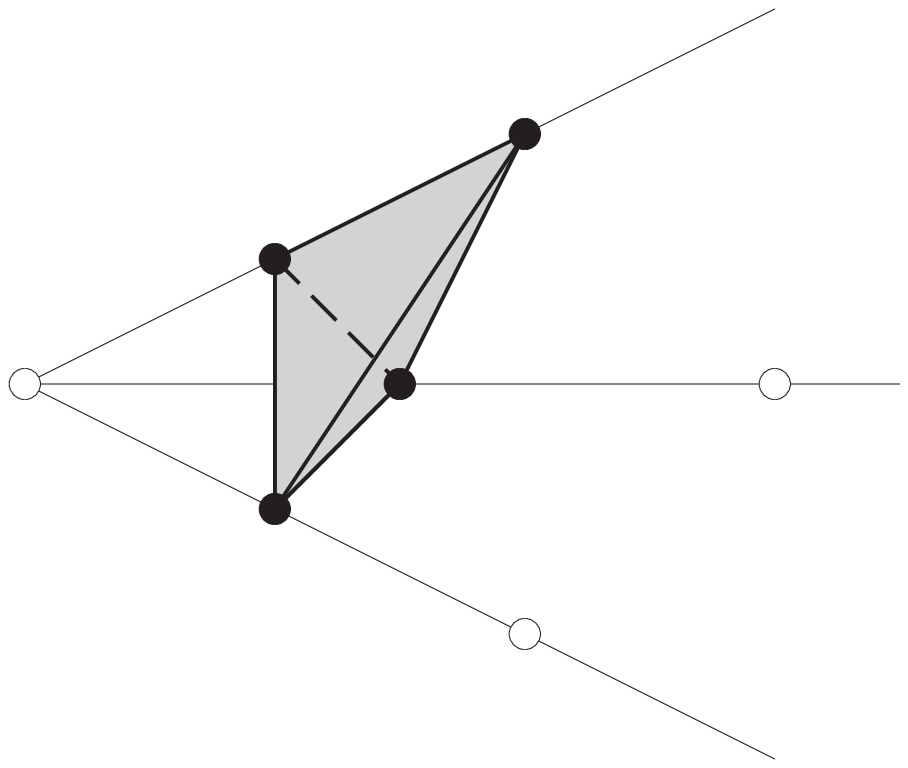}\\
$[v_1^2,v_1^1, v_2^1, v_3^1],$
\end{center}
\end{minipage}
\begin{minipage}{0.22\textwidth}
\begin{center}
\includegraphics[scale=0.35]{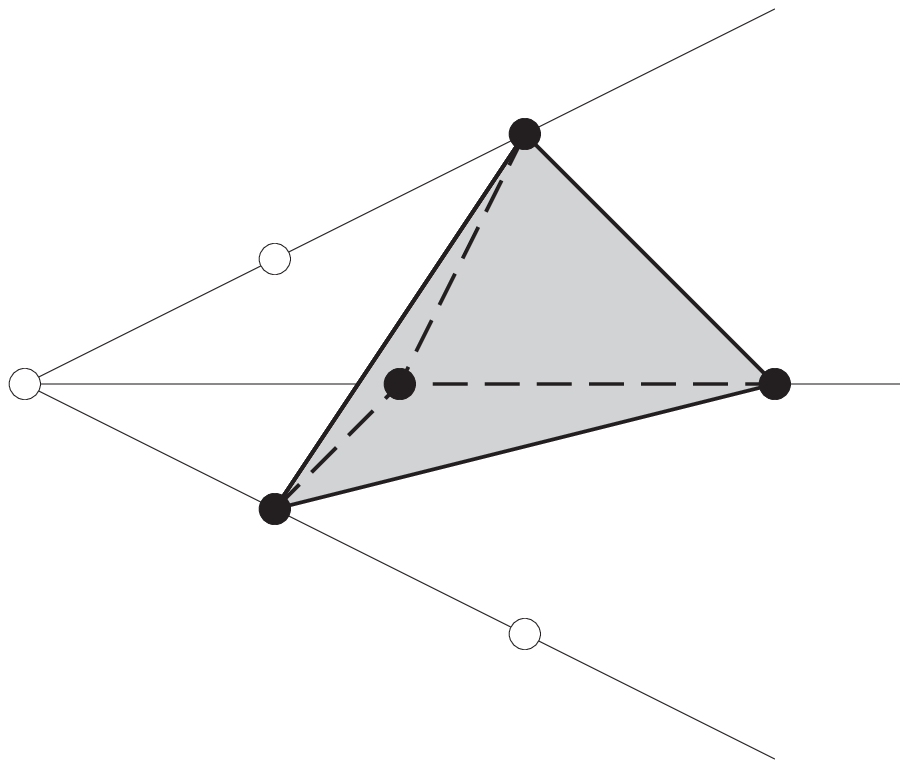}\\
$[v_1^2,v_2^2, v_2^1, v_3^1],$
\end{center}
\end{minipage}
\begin{minipage}{0.22\textwidth}
\begin{center}
\includegraphics[scale=0.35]{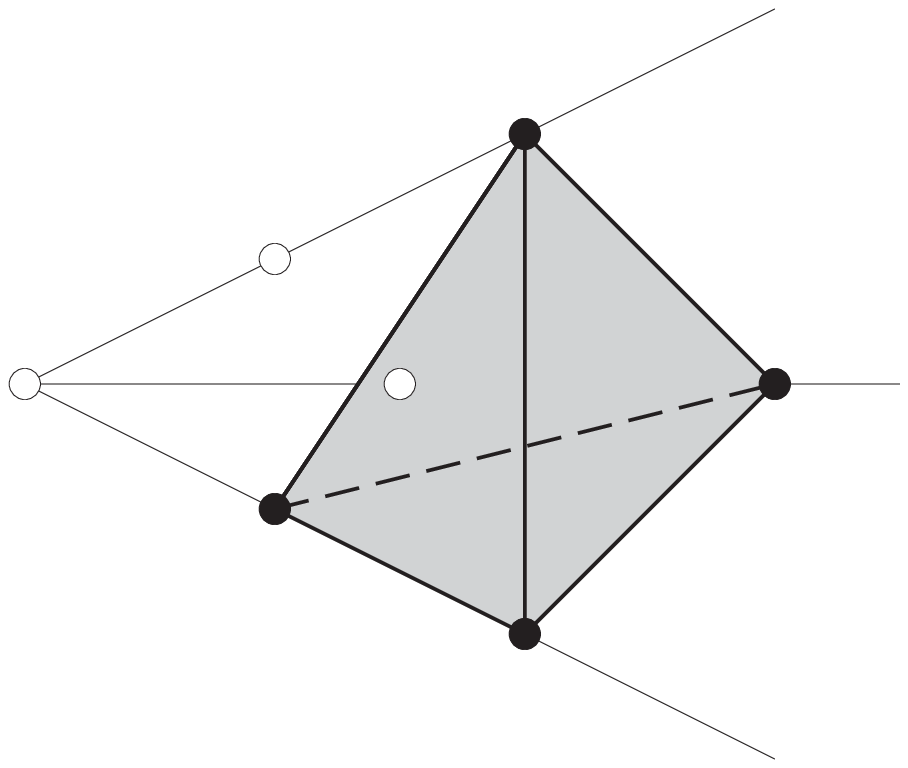}\\
$[v_1^2,v_2^2, v_3^2, v_3^1].$
\end{center}
\end{minipage}
\end{center}
Putting them together, we obtain Figure~\ref{dec:fig:triangulation_3_cone}.

\begin{figure}[H]
\begin{center}
\includegraphics[scale=0.70]{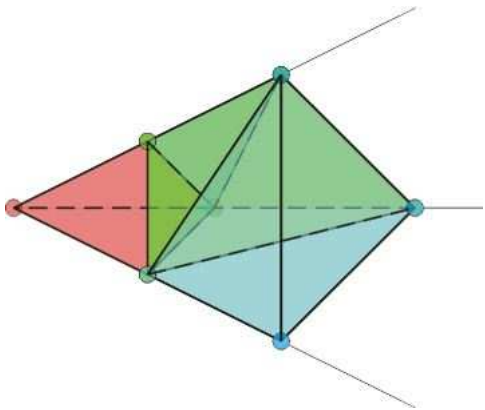}
\end{center}
\caption{\label{dec:fig:triangulation_3_cone}Triangulation of a three-dimensional cone.}
\end{figure}

This example is significant, since we have demonstrated that for
any $n$-dimensional ball about a point, we can construct a
trivially star-shaped triangulation of the ball, with arbitrarily
high resolution. This allows us to recover the smooth
Poincar\'{e} lemma in the limit of an infinitely fine mesh, using
the discrete Poincar\'{e} lemma for trivially star-shaped
complexes.
\end{example}
\paragraph{Logically Star-Shaped Complexes.}
\begin{definition}
A simplicial complex $L$ is \textbf{logically star-shaped}\index{star-shaped!logically} if it
is isomorphic, at the level of an abstract simplicial complex, to a
trivially star-shaped complex $K$.
\end{definition}
\begin{example}
We see two simplicial complexes, in Figure~\ref{dec:fig:isomorphic_abstract}, which are clearly isomorphic as abstract simplicial complexes.

\begin{figure}[htbp]
\begin{center}
\begin{minipage}{0.4\textwidth}
\begin{center}
\includegraphics[scale=0.45]{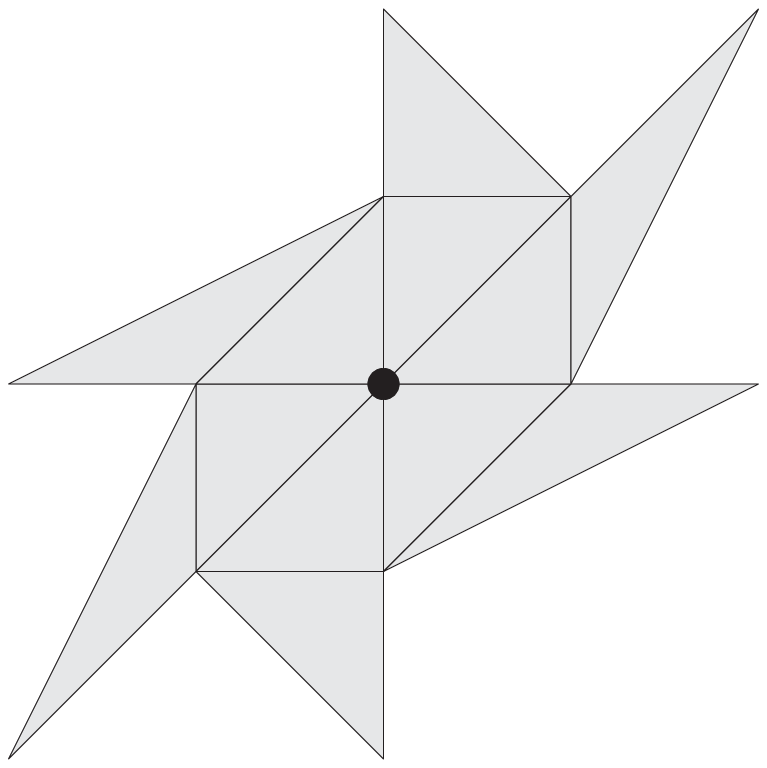}
\end{center}
\end{minipage}$\cong$
\begin{minipage}{0.4\textwidth}
\begin{center}
\includegraphics[scale=0.45]{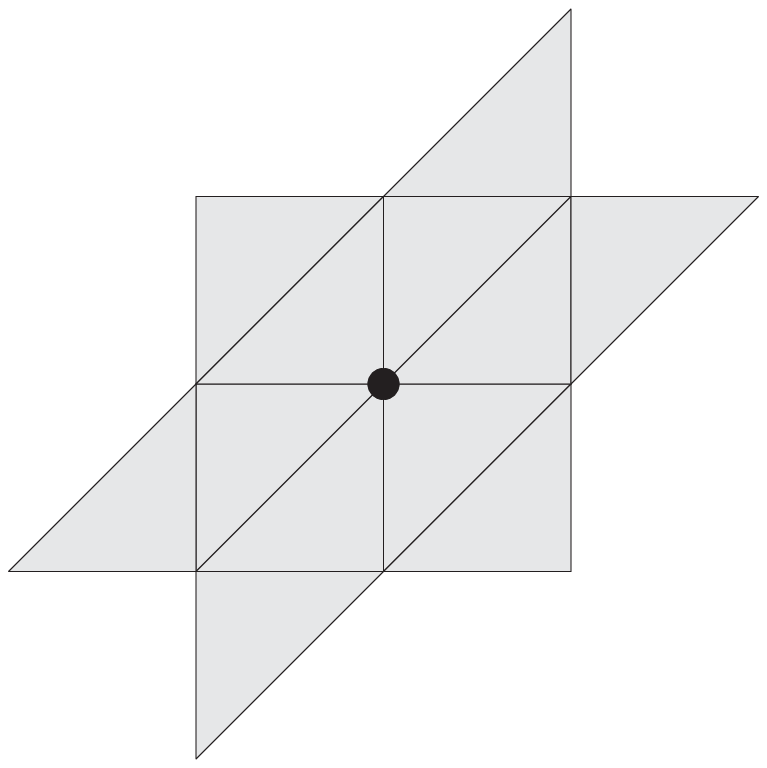}
\end{center}
\end{minipage}
\end{center}
\caption{\label{dec:fig:isomorphic_abstract}Trivially star-shaped complex (left); Logically star-shaped complex (right).}
\end{figure}
\end{example}
\begin{definition}
The \textbf{logical cone}\index{cone!logical} operator $p:C^k(L)\rightarrow
C^{k+1}(L)$ is defined by making the following diagram commute,
\[\xymatrix{
C^k (K) \ar[r]^{p_K} \ar@{=}[d] &
C^{k+1}(K) \ar@{=}[d]\\
C^k (L) \ar[r]^{p_L} & C^{k+1}(L)}\] Which is to say that, given
the isomorphism $\varphi:K\rightarrow L$, we define
\[p_L=\varphi\circ p_K\circ\varphi^{-1}.\]
\end{definition}
\begin{example}\index{cone!example}
We show an example of the construction of the logical cone
operator.
\[\xymatrix@!0@R=4cm@C=4cm{
\includegraphics[scale=0.3]{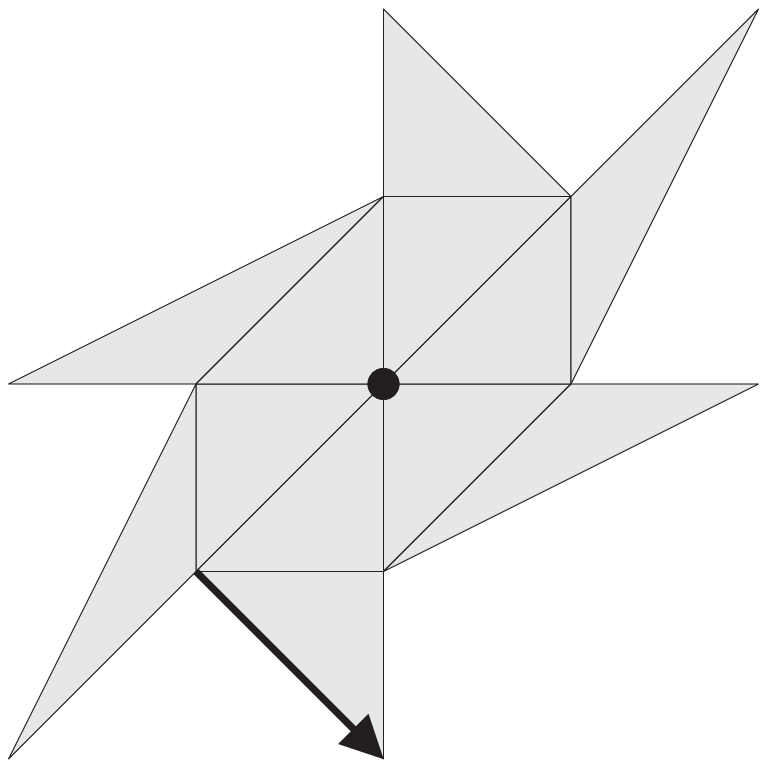} \ar[r]^{p_K} \ar@{=}[d] &
\includegraphics[scale=0.3]{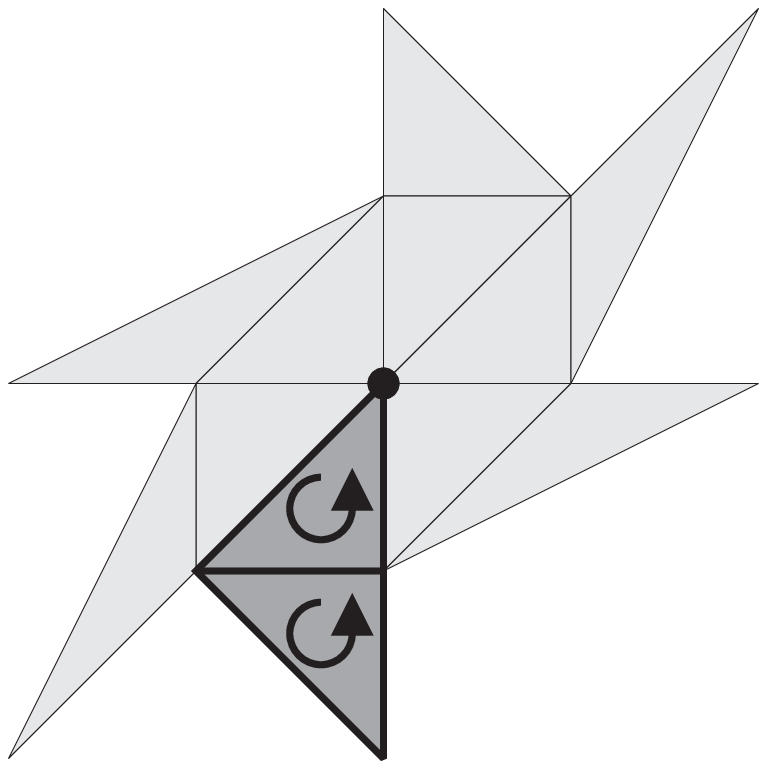} \ar@{=}[d]\\
\includegraphics[scale=0.3]{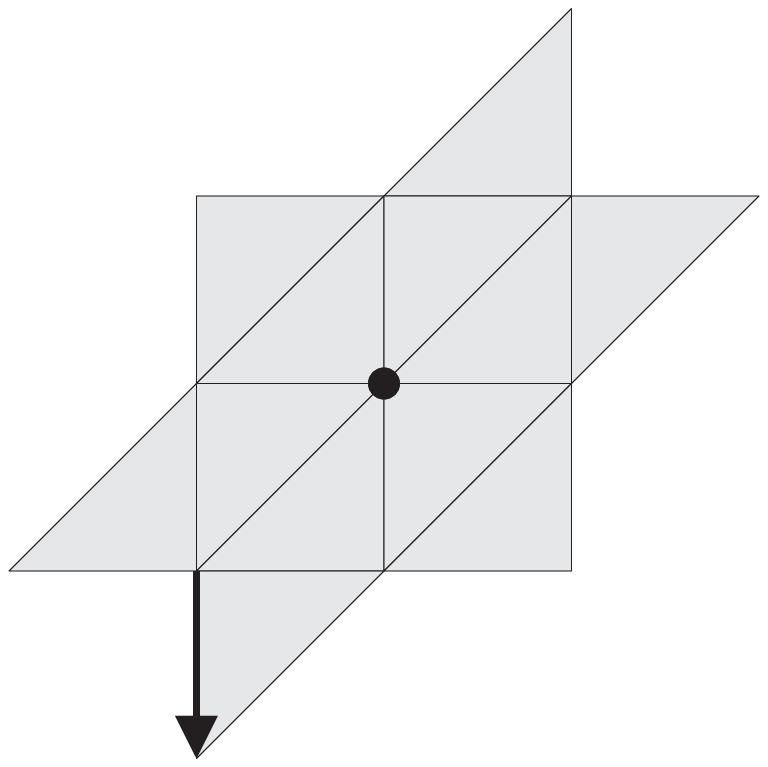} \ar[r]^{p_L} &
\includegraphics[scale=0.3]{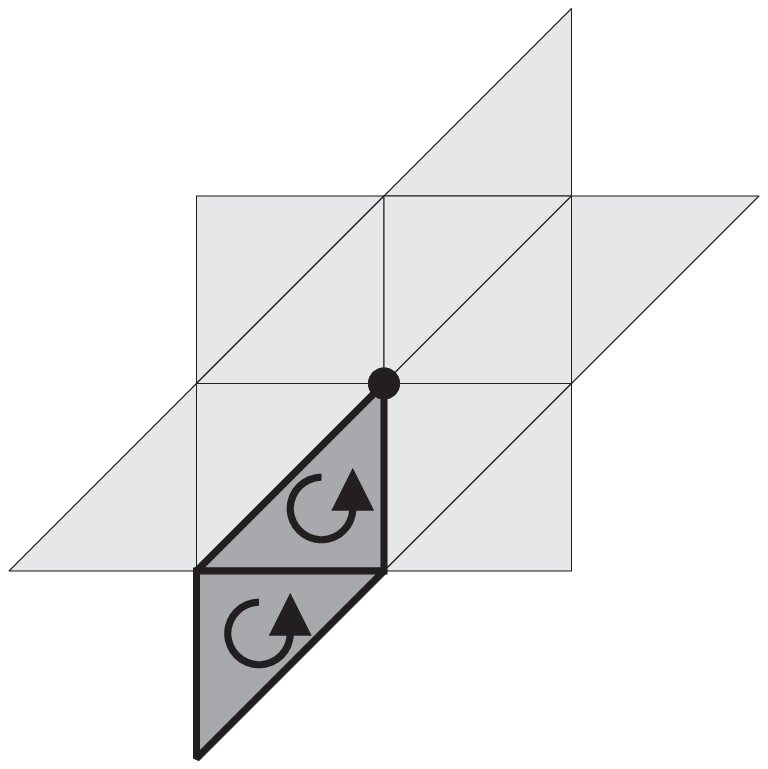}}\]
\end{example}
This definition of the logical cone operator results in
identities for the cone and cocone operator that follow from the
trivially star-shaped case, and we record the results as follows.
\begin{lemma}
In logically star-shaped complexes, the logical cone operator
satisfies the following identity,
\[
p\partial+\partial p = I,
\]
at the level of chains.
\end{lemma}
\begin{proof}
Follows immediately by pushing forward the result for trivially
star-shaped complexes using the isomorphism.
\end{proof}
\begin{lemma}
In logically star-shaped complexes, the logical cocone operator
satisfies the following identity,
\[
H \d + \d H = I,
\]
at the level of cochains.
\end{lemma}
\begin{proof}
Follows immediately by pushing forward the result for trivially
star-shaped complexes using the isomorphism.
\end{proof}
Similarly, we have a Discrete Poincar\'{e} Lemma for logically
star-shaped complexes.
\begin{corollary}[Discrete Poincar\'{e} Lemma for Logically
Star-shaped Complexes] Given a closed cochain $\alpha^k$, that is
to say, $\d\alpha^k=0$, there exists a cochain $\beta^{k-1}$, such
that, $\d\beta^{k-1}=\alpha^k$.
\end{corollary}
\begin{proof}
Follows from the above lemma using the proof for the trivially
star-shaped case.
\end{proof}
\paragraph{Contractible Complexes.}
For arbitrary contractible complexes, we construct a generalized
cone operator such that it satisfies the identity,
\[ p \partial + \partial p = I,\]
which is the crucial property of the cone operator, from the point
of view of proving the discrete Poincar\'{e} lemma.

The trivial cone construction gives a clue as to how to proceed in
the construction of a {\bfi generalized cone
operator}\index{cone!generalized}. Notice that if a $\sigma^{k+1}$
is a term in $p(\sigma^k)$, then $p(\sigma^{k+1})=\varnothing$.
This suggests how we can use the cone identity to inductively
construct the generalized cone operator.

To define $p(\sigma^k)$, we consider $\sigma^{k+1}\succ\sigma^k$,
such that, $\sigma^{k+1}$ and $\sigma^k$ are consistently oriented.
We apply $p\partial+\partial p$ to $\sigma^{k+1}$. Then, we have
\begin{align*} \sigma^{k+1} &=
p(\sigma^k)+p(\partial\sigma^{k+1}-\sigma^k) + \partial
p(\sigma^{k+1}).\\ \intertext{If we set
$p(\sigma^{k+1})=\varnothing$,} \sigma^{k+1} &=
p(\sigma^k)+p(\partial\sigma^{k+1}-\sigma^k) +
\partial (\varnothing)\\ &=
p(\sigma^k)+p(\partial\sigma^{k+1}-\sigma^k).\end{align*}
Rearranging, we have
\[p(\sigma^k) = \sigma^{k+1} -
p(\partial\sigma^{k+1}-\sigma^k),\] and
\[p(\sigma^{k+1})=\varnothing.\]We are done, so long as
the simplices in the chain $\partial\sigma^{k+1}-\sigma^k$ already
have $p$ defined on it. This then reduces to enumerating the
simplices in such a way that in the right hand side of the
equation, we never evoke terms that are undefined.

We now introduce a method of augmenting a complex so that the
enumeration condition is always satisfied.

\begin{definition} Given a $n$-complex
$K$, consider a $(n-1)$-chain $c_{n-1}$ that is contained on the
boundary of $K$, and is included in the one-ring of some vertex on
$\partial K$. Then, the \textbf{one-ring cone augmentation}\index{augmentation!one-ring cone} of $K$
is the complex obtained by adding the $n$-cone $w\diamond
c_{n-1}$, and all its faces to the complex.
\end{definition}

\begin{definition}
A complex is \textbf{generalized star-shaped}\index{star-shaped!generalized} if it can be
constructed by repeatedly applying the one-ring augmentation
procedure.
\end{definition}

We will explicitly show in Examples~\ref{oneringaug2}, and
\ref{oneringaug3}, how to enumerate the vertices in two and
three dimensions. And in Examples~\ref{2regtetra}, and \ref{3regtetra}, we
will introduce regular triangulations of $\mathbb{R}^2$ and
$\mathbb{R}^3$ that can be constructed by inductive one-ring cone
augmentation.

\begin{remark}
Notice that a non-contractible complex cannot be constructed by
inductive one-ring cone augmentation, since it will involve adding
a cone to a vertex that has two disjoint base chains. This
prevents us from enumerating the simplices in such as way that all
the terms in $\partial\sigma^{k+1}-\sigma^k$ have had $p$ defined
on them, and we see in Example~\ref{cex} how this causes the cone
identity, and hence the discrete Poincar\'{e} lemma to break.
\end{remark}

\begin{example}\label{oneringaug2}\index{augmentation!example}
In two dimensions, the one-ring condition implies that the base of the
cone consists of either one or two $1$-simplices. To aid in
visualization, consider Figure~\ref{dec:fig:onering2}.

\begin{figure}[htbp]
\begin{center}
\includegraphics[scale=1]{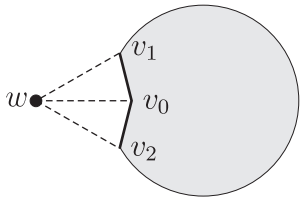}
\end{center}
\caption{\label{dec:fig:onering2}One-ring cone augmentation of a complex in two dimensions.}
\end{figure}
In the case of one $1$-simplex, $[v_0,v_1]$, when we augment using
the cone construction with the new vertex $w$, we define,
\begin{align*}
p([w]) &= [v_0,w]+p([v_0]), &p([v_0,w]) &=\varnothing,\\
p([v_1,w]) &=[v_0,v_1,w]-p([v_0,v_1]),
&p([v_0,v_1,w])&=\varnothing.
\end{align*}
In the case of two $1$-simplices, $[v_0,v_1]$, $[v_0,v_2]$, we have,
\begin{align*}
p([w]) &= [v_0,w]+p([v_0]), &p([v_0,w]) &=\varnothing,\\
p([v_1,w]) &=[v_0,v_1,w]-p([v_0,v_1]), &p([v_0,v_1,w])&=\varnothing,\\
p([v_2,w]) &=[v_0,v_2,w]-p([v_0,v_2]),
&p([v_0,v_2,w])&=\varnothing.
\end{align*}
\end{example}
\begin{example}\index{cone!example}
We will now explicitly utilize the one-ring cone augmentation
procedure to compute the generalized cone operator for part of a
regular two-dimensional triangulation that is not logically
star-shaped.

As a preliminary, we shall consider a logically star-shaped
complex, and augment with a new vertex, as seen in Figure~\ref{dec:fig:star_shaped_augmented}.

\begin{figure}[htbp]
\begin{center}
\includegraphics[scale=0.4]{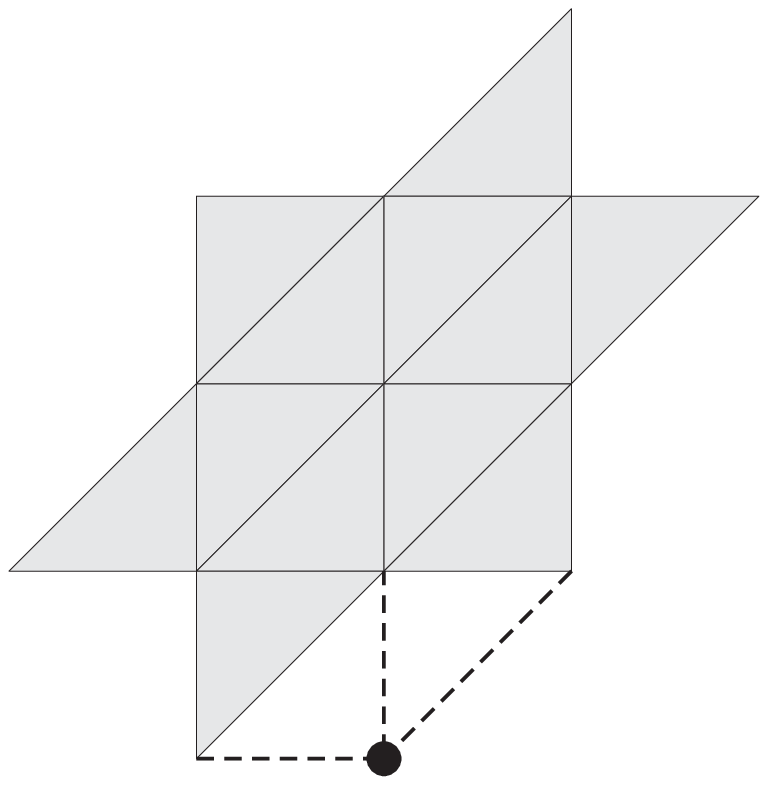}
\end{center}
\caption{\label{dec:fig:star_shaped_augmented}Logically star-shaped complex augmented by cone.}
\end{figure}

We use the logical cone operator for the subcomplex that is
logically star-shaped, and the augmentation rules in the example
above for the newly introduced simplices. This yields,
\begin{align*}
p\begin{pmatrix}
\raisebox{-6.9ex}{\includegraphics[scale=0.25]{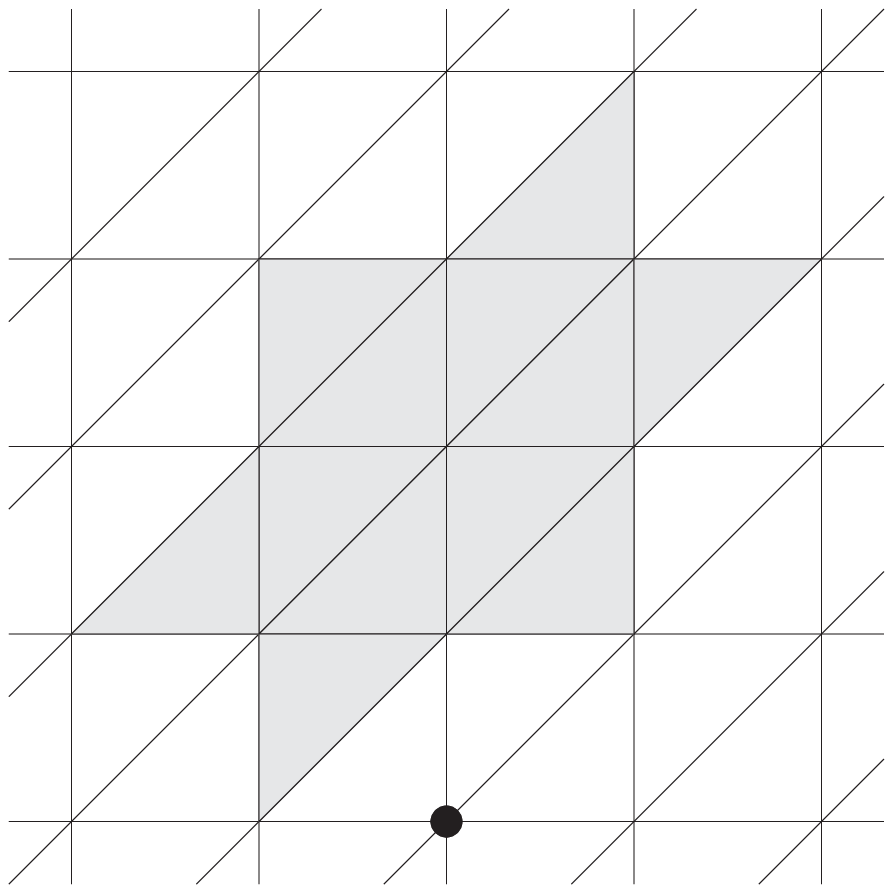}}
\end{pmatrix}
&=\raisebox{-6.9ex}{\includegraphics[scale=0.25]{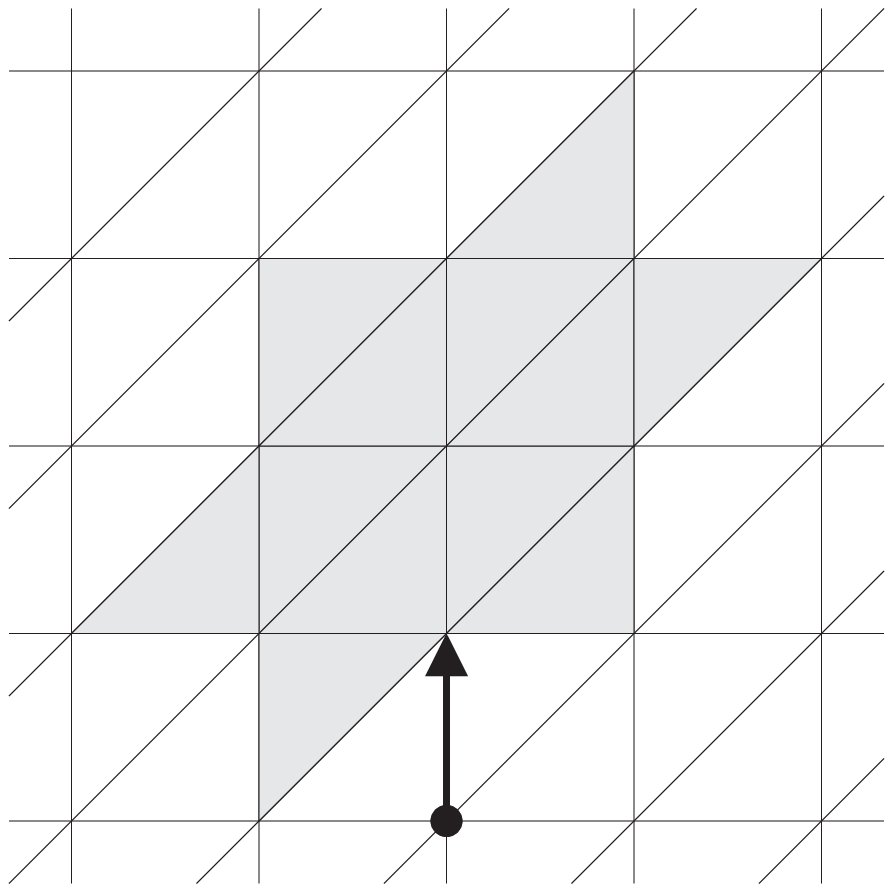}} +
p
\begin{pmatrix}\raisebox{-6.9ex}{\includegraphics[scale=0.25]{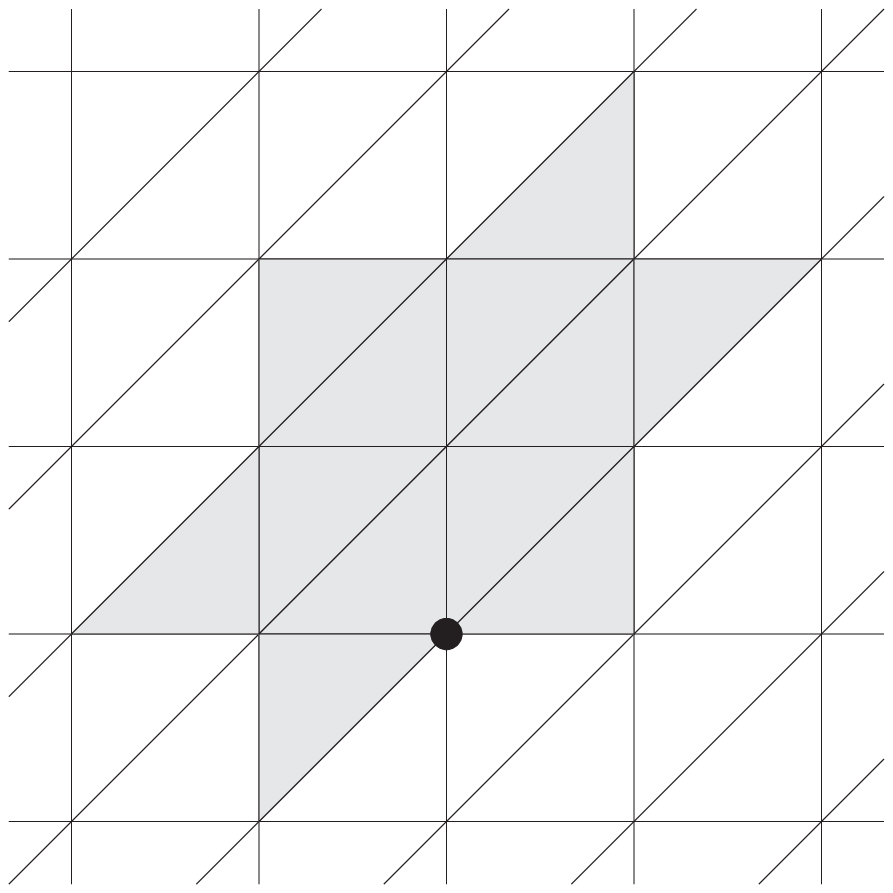}}\end{pmatrix}
=\raisebox{-6.9ex}{\includegraphics[scale=0.25]{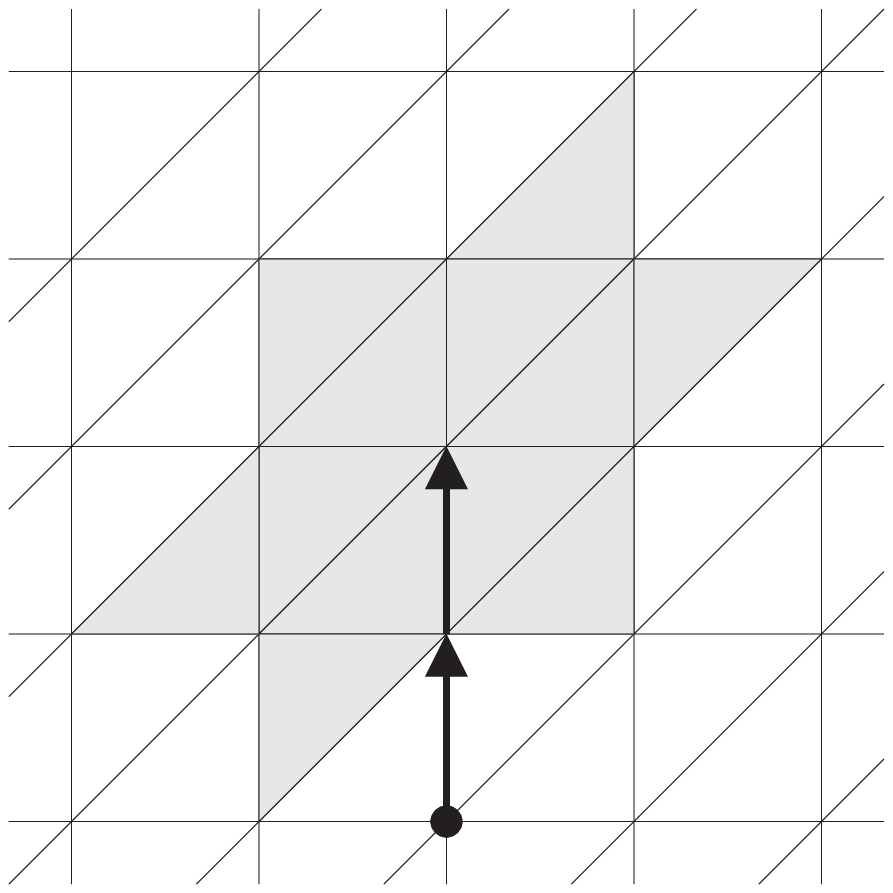}},\\
p\begin{pmatrix}\raisebox{-6.9ex}{\includegraphics[scale=0.25]{p01}}\end{pmatrix}&=\varnothing,\\
p\begin{pmatrix}\raisebox{-6.9ex}{\includegraphics[scale=0.25]{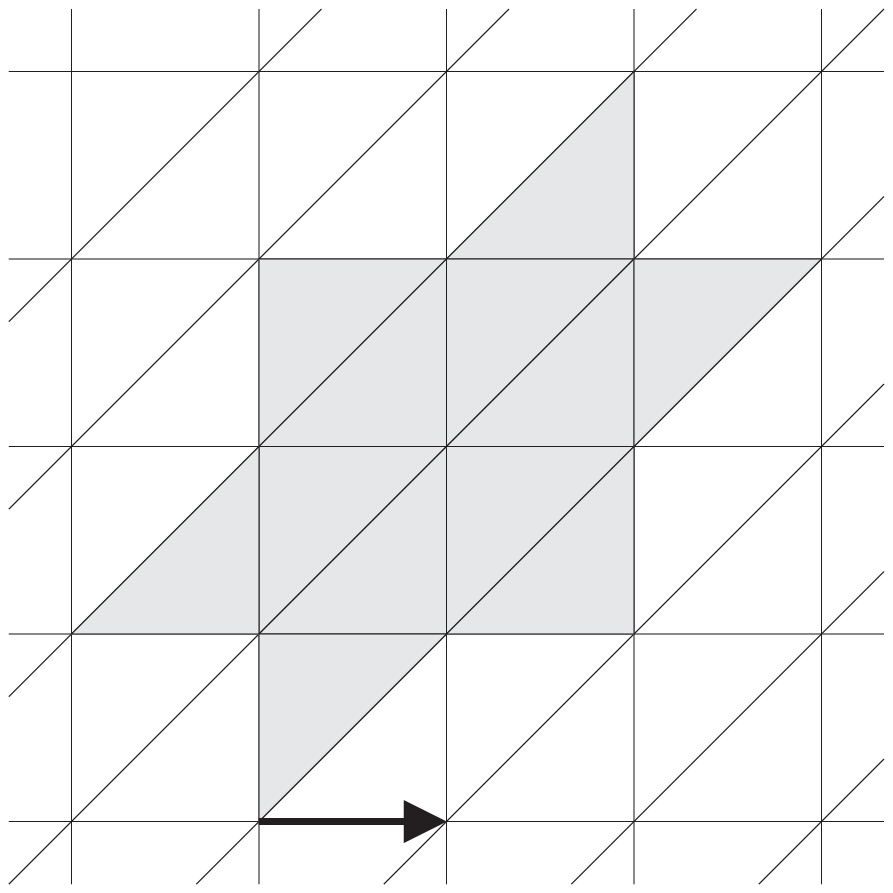}}\end{pmatrix}&=
\raisebox{-6.9ex}{\includegraphics[scale=0.25]{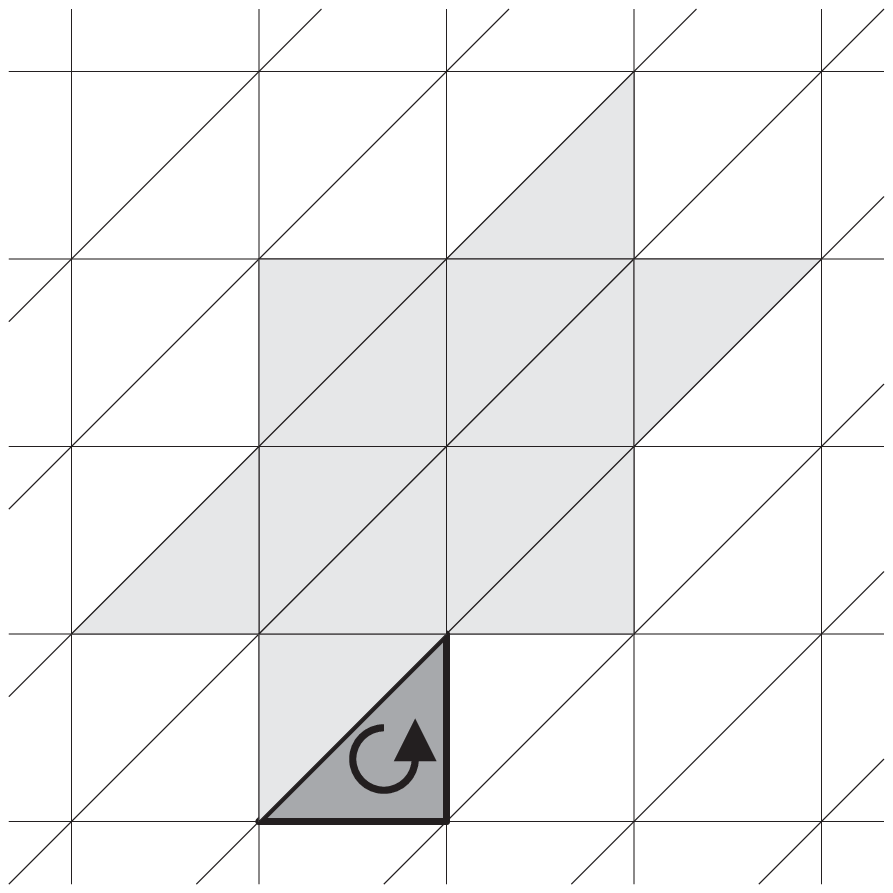}}+p\begin{pmatrix}\raisebox{-6.9ex}{\includegraphics[scale=0.25]{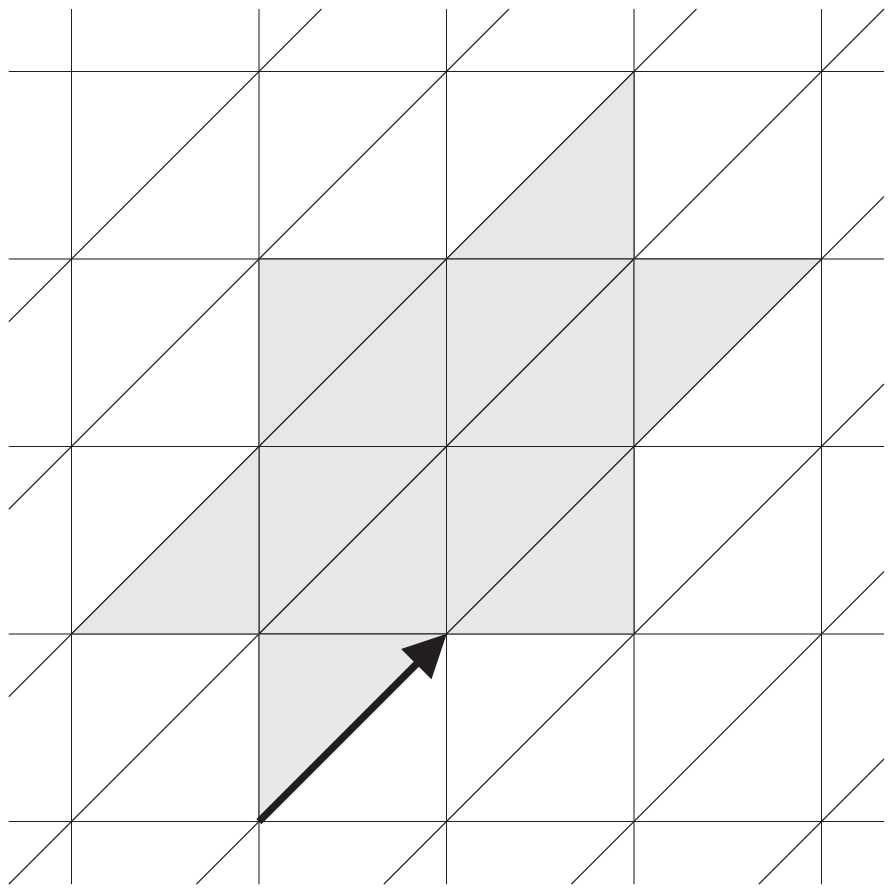}}\end{pmatrix}
=\raisebox{-6.9ex}{\includegraphics[scale=0.25]{p11}}+\varnothing\\
&=\raisebox{-6.9ex}{\includegraphics[scale=0.25]{p11}},\\
p\begin{pmatrix}\raisebox{-6.9ex}{\includegraphics[scale=0.25]{p11}}\end{pmatrix}&=\varnothing,\\
p\begin{pmatrix}\raisebox{-6.9ex}{\includegraphics[scale=0.25]{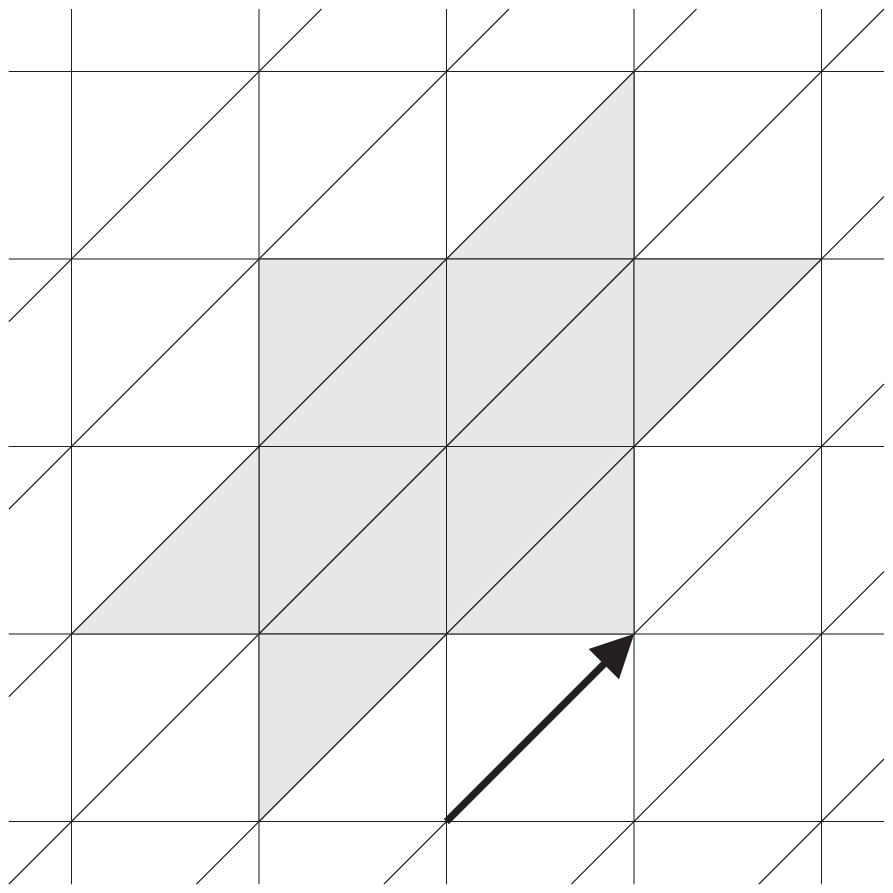}}\end{pmatrix}&=
\raisebox{-6.9ex}{\includegraphics[scale=0.25]{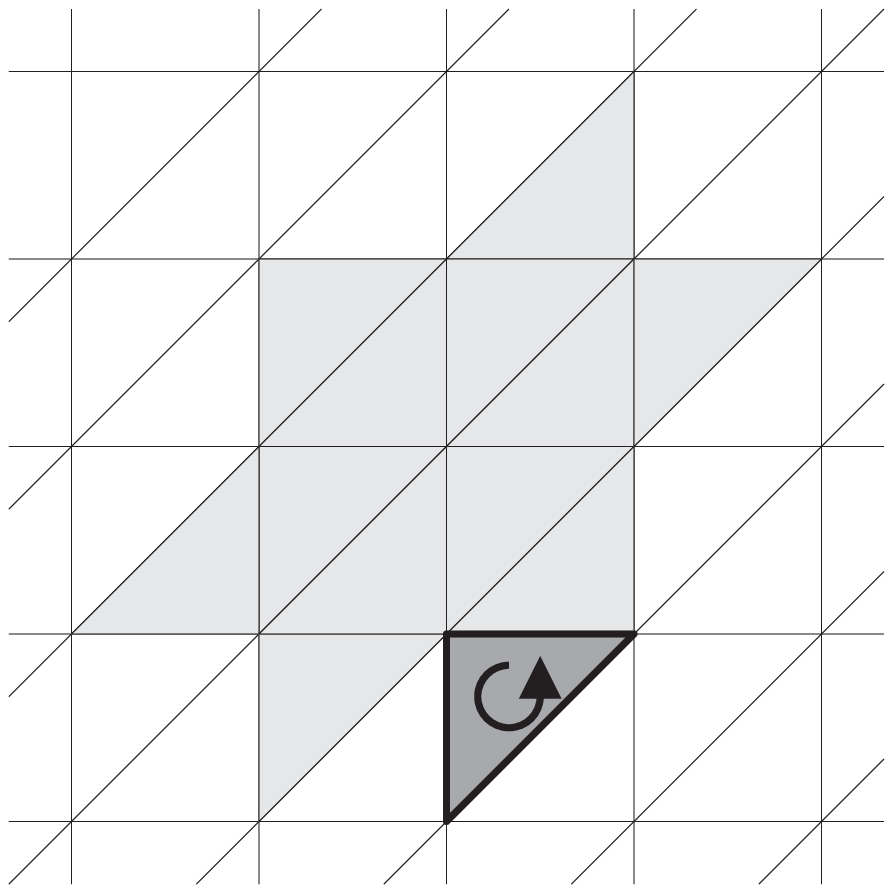}}+p\begin{pmatrix}\raisebox{-6.9ex}{\includegraphics[scale=0.25]{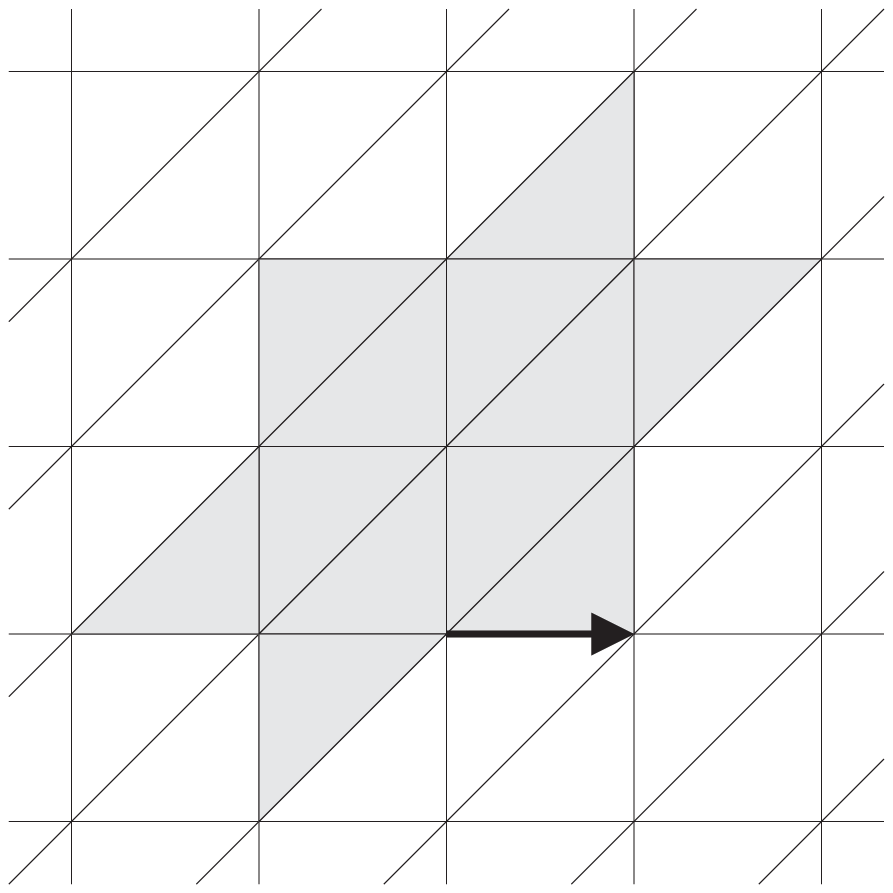}}\end{pmatrix}\\
&=\raisebox{-6.9ex}{\includegraphics[scale=0.25]{p21}}+\raisebox{-6.9ex}{\includegraphics[scale=0.25]{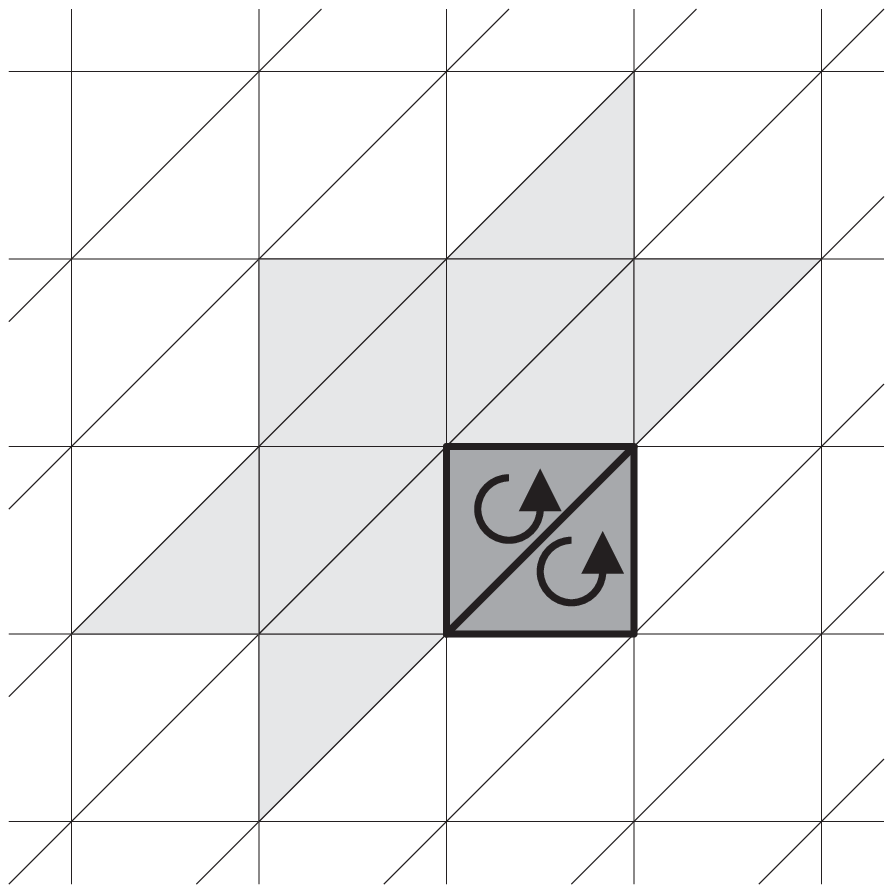}}
=\raisebox{-6.9ex}{\includegraphics[scale=0.25]{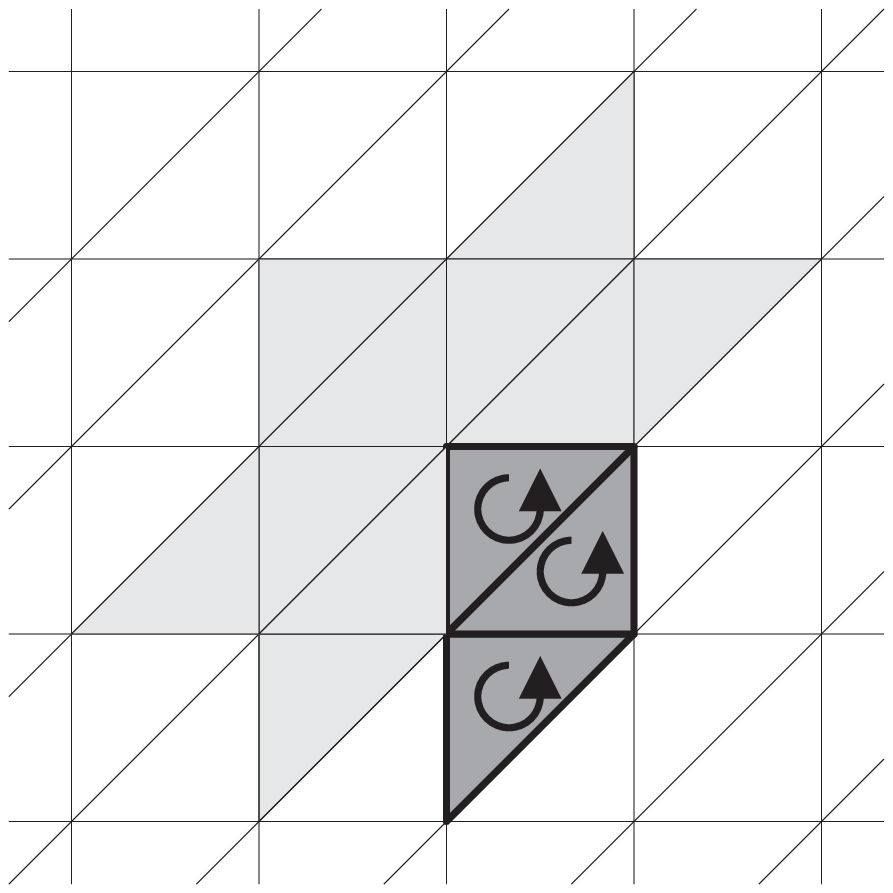}},\\
p\begin{pmatrix}\raisebox{-6.9ex}{\includegraphics[scale=0.25]{p21}}\end{pmatrix}&=\varnothing.
\end{align*}
\end{example}
\begin{example}\label{2regtetra}
Clearly, the regular two-dimensional triangulation can be
obtained by the successive application of the one-ring cone
augmentation procedure, as the following sequence illustrates,
\begin{center}
\raisebox{-6.9ex}{\includegraphics[scale=0.25]{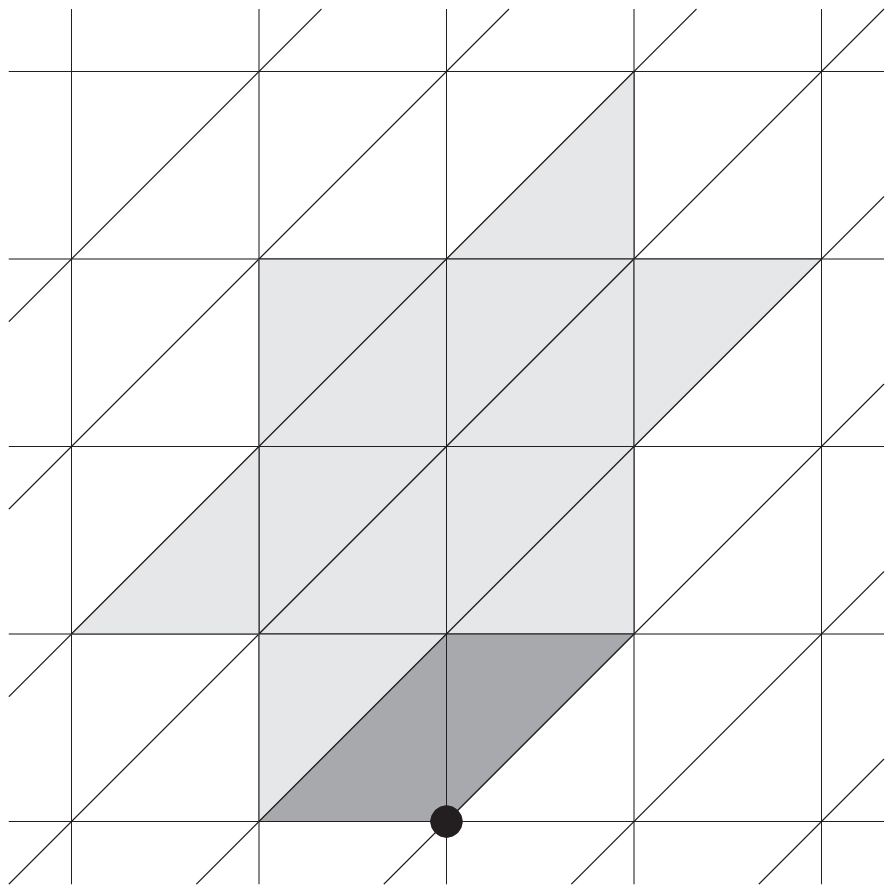}}\quad$\mapsto$\quad
\raisebox{-6.9ex}{\includegraphics[scale=0.25]{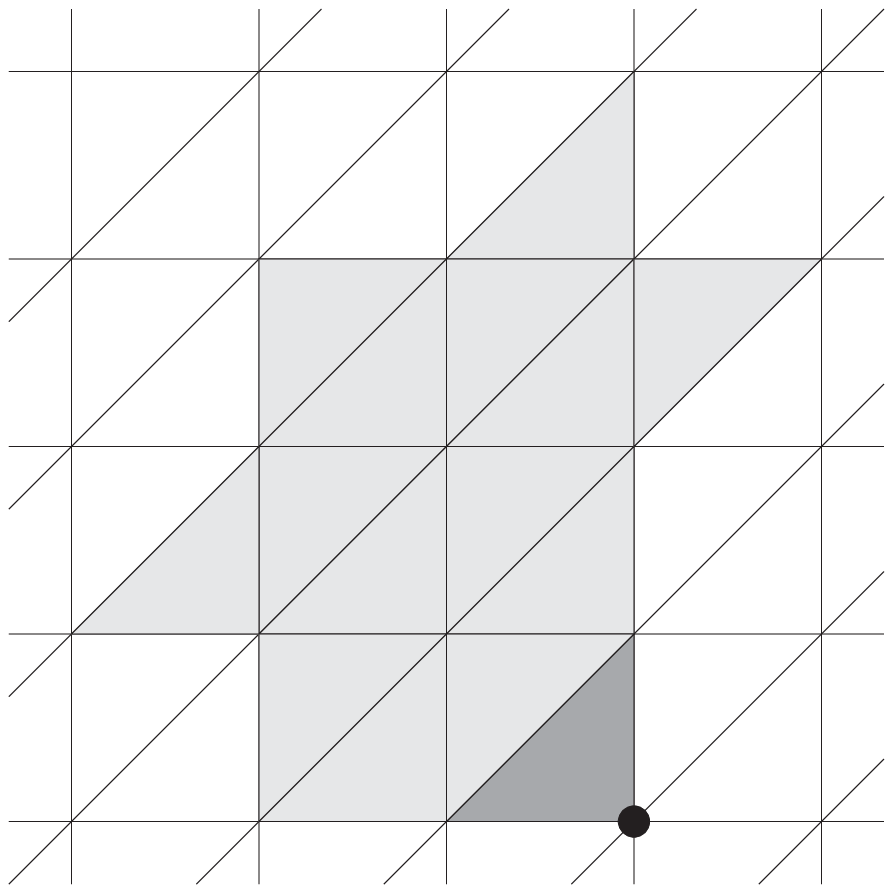}}\quad$\mapsto$\quad
\raisebox{-6.9ex}{\includegraphics[scale=0.25]{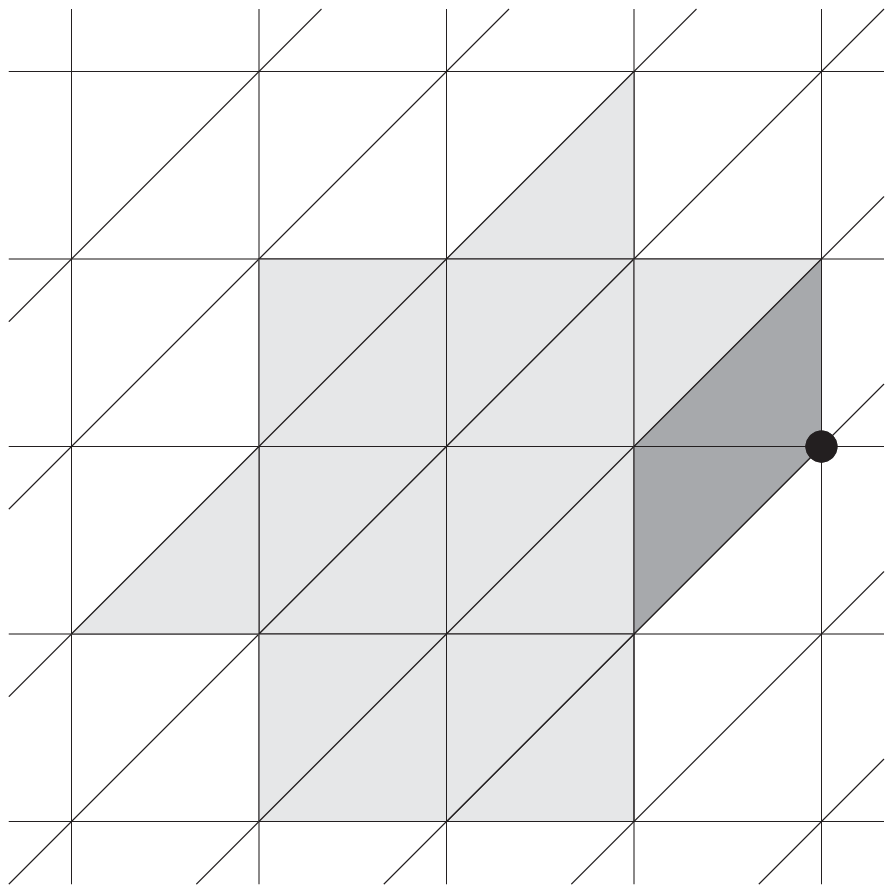}}\quad$\mapsto$\quad
\raisebox{-6.9ex}{\includegraphics[scale=0.25]{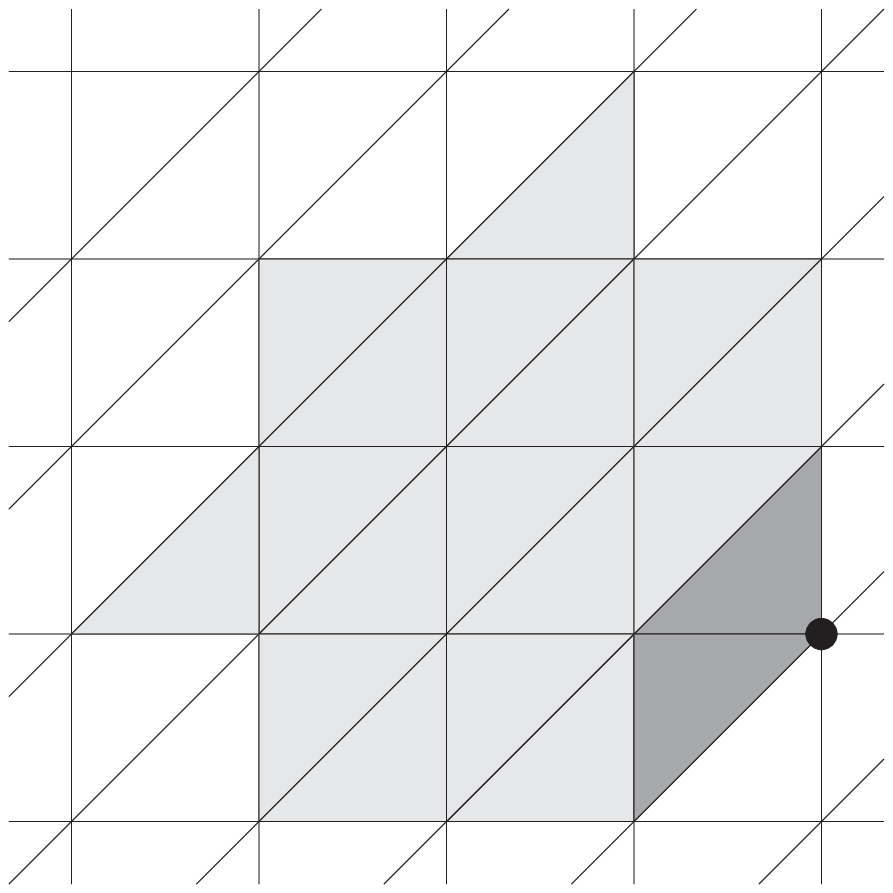}}\quad$\mapsto\quad\ldots\quad$,
\end{center}
which means that the discrete Poincar\'{e} lemma can be extended
to the entire regular triangulation of the plane.
\end{example}
\begin{example}\label{oneringaug3}\index{augmentation!example} We consider the case of augmentation in three dimensions. Denote by $v_0$, the center of the one-ring on the
two-surface, to which we are augmenting the new vertex $w$. The
other vertices of the one-ring are enumerated in order, $v_1,\ldots,
v_m$. To aid in visualization, consider Figure~\ref{dec:fig:onering3}.

\begin{figure}[htbp]
\begin{center}
\includegraphics[scale=1]{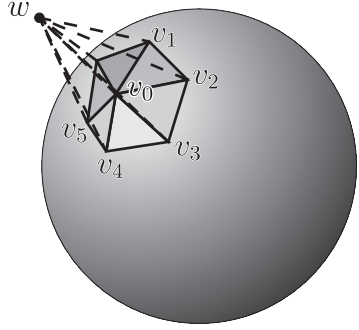}
\end{center}
\caption{\label{dec:fig:onering3}One-ring cone augmentation of a complex in
three dimensions.}
\end{figure}

If the one-ring does not go completely around $v_0$, we shall denote
the missing term by $[v_0,v_1,v_m]$. The generalized cone operators are given as follows.
\begin{align*}
\intertext{k=0,}
p([w])&= [v_0,w]+p([v_0]), &p([v_0,w]) &= \varnothing,\\
\intertext{k=1,} p([v_1,w])&= [v_0,v_1,w] - p([v_0,v_1]), &p([v_0,v_1,w]) &= \varnothing,\\
p([v_m,w])&= [v_0,v_m,w] - p([v_0,v_m]), &p([v_0,v_m,w]) &= \varnothing,\\
\intertext{k=2,} p([v_1,v_2,w]) &= [v_0,v_1,v_2,w] +
p([v_0,v_1,v_2]), &p([v_1,v_2,w]) &=
\varnothing,\\
p([v_{m-1},v_m,w]) &= [v_0,v_{m-1},v_m,w] + p([v_0,v_{m-1},v_m]),   &p([v_{m-1},v_m,w]) &= \varnothing.\\
\intertext{If it does go around completely,} p([v_m,v_1,w]) &=
[v_0,v_m,v_1,w] + p([v_0,v_m,v_1]), &p([v_0,v_m,v_1,w]) &=
\varnothing.
\end{align*}
\end{example}

\begin{example}\label{3regtetra}
We provide a tetrahedralization of the unit cube that can be tiled
to yield a regular tetrahedralization of $\mathbb{R}^3$. The
$3$-simplices are as follows,
\begin{align*}
[v_{000},v_{001},v_{010},v_{10}],
[v_{001},v_{010},v_{100},v_{101}],[v_{001},v_{010},v_{011},v_{101}],\\
[v_{010},v_{100},v_{101},v_{110}],[v_{010},v_{011},v_{101},v_{110}],
[v_{011},v_{101},v_{110},v_{111}].
\end{align*}
The tetrahedralization of the unit cube can be seen in Figure~\ref{dec:fig:tetrahedralization}.

\begin{figure}[htbp]
\begin{center}
\subfigure[Tileable tetrahedralization of the unit cube]{\quad\includegraphics[scale=0.5]{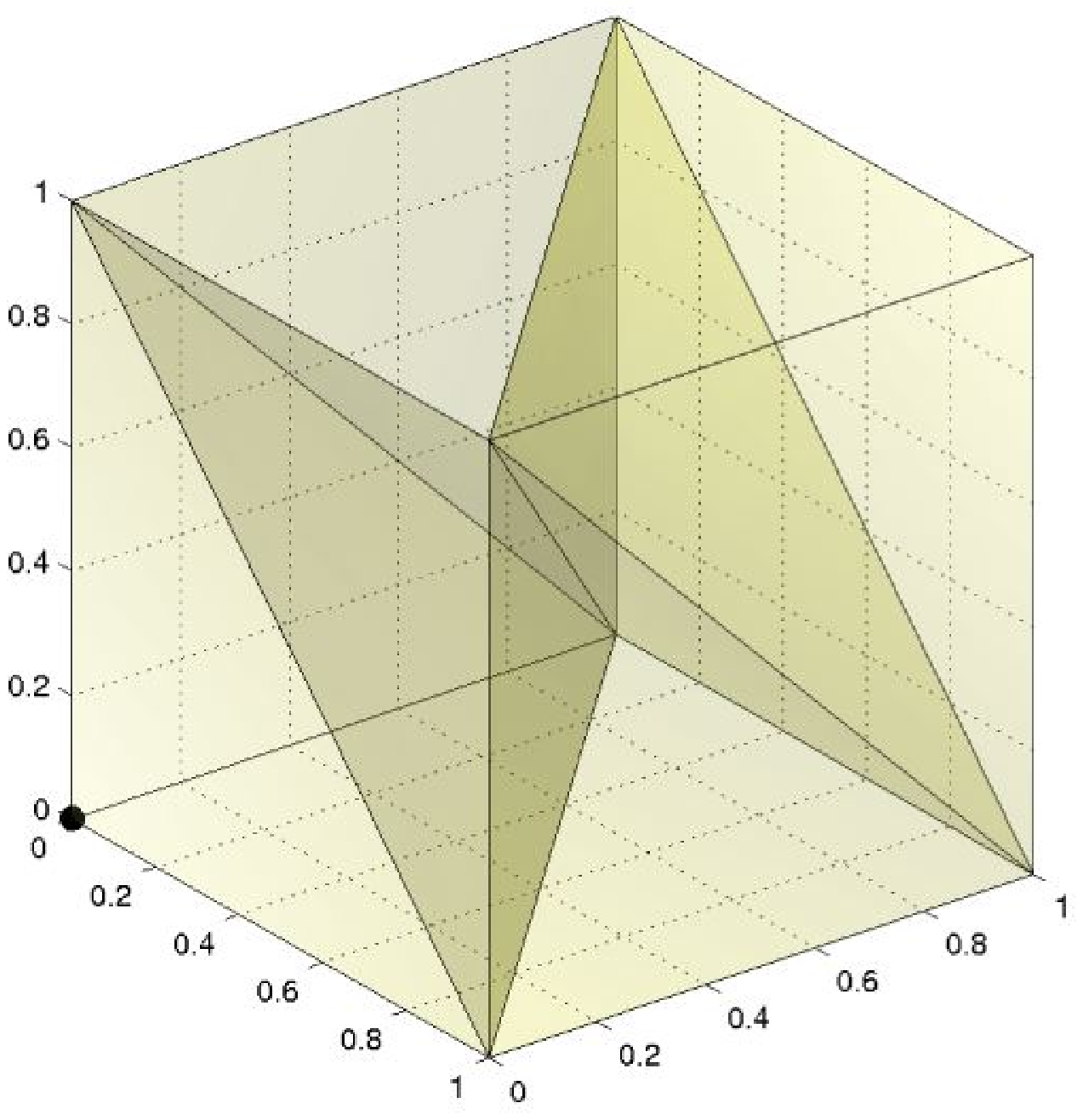}\quad}
\subfigure[Partial tiling of $\mathbb{R}^3$]{\quad\includegraphics[scale=0.5]{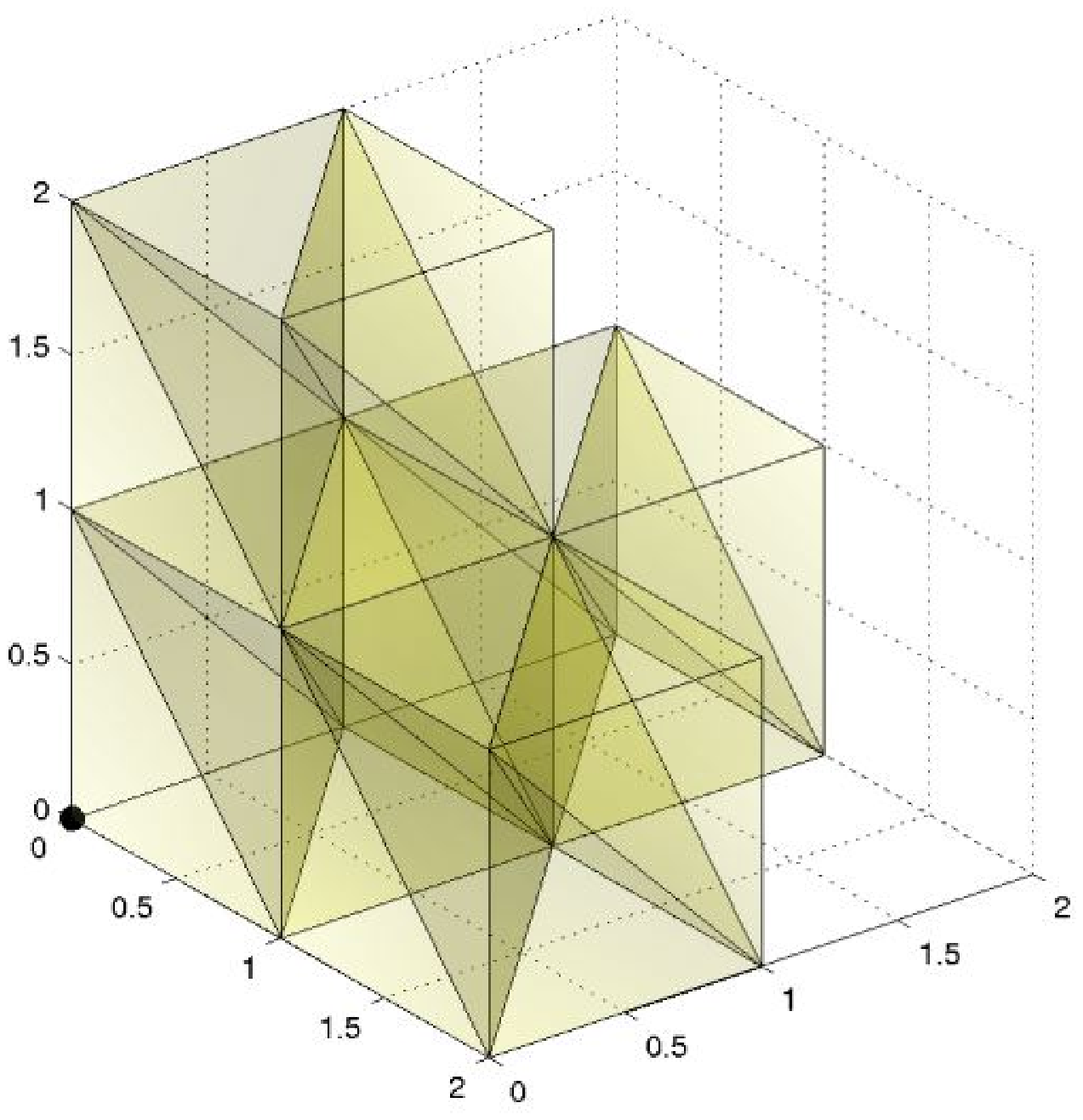}\quad}
\caption{\label{dec:fig:tetrahedralization}Regular tiling of $\mathbb{R}^3$ that admits a generalized cone operator.}
\end{center}
\end{figure}

Since this regular tetrahedralization can be constructed by the
successive application of the one-ring cone augmentation
procedure, the Discrete Poincar\'{e} lemma can be extended to the
entire regular tetrahedralization of $\mathbb{R}^3$.
\end{example}
In higher dimensions, we can extend the construction of the
generalized cone operator inductively using the one-ring cone
augmentation by choosing an appropriate enumeration of the base
chain. Topologically, the base chain will be the cone of
$S^{n-2}$ (with possibly an open $(n-2)$-ball removed) with respect
to the central point.

\begin{figure}[htbp]
\begin{center}
\includegraphics[scale=0.5]{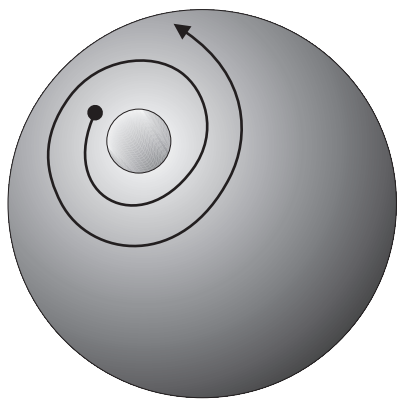}
\end{center}
\caption{\label{dec:fig:spiral_s_n-2}Spiral enumeration of $S^{n-2}$, $n=4$.}
\end{figure}

By spiraling around $S^{n-2}$, starting from around the boundary
of the $n-2$ ball, and covering the rest of $S^{n-2}$, as in Figure~\ref{dec:fig:spiral_s_n-2}, we obtain
the higher-dimensional generalization of the procedure we have
taken in Examples~\ref{oneringaug2}, and \ref{oneringaug3}.

Notice that $n=2$ is distinguished, since $S^{2-2}=S^0$ is
disjoint, which is why in the two-dimensional case, we were not able
to use the spiraling technique to enumerate the simplices.

Since we have constructed the generalized cone operator such that
the cone identity holds, we have,
\begin{lemma} In generalized star-shaped complexes, the generalized cone
operator satisfies the following identity,
\[
p\partial+\partial p = I,
\]
at the level of chains.
\end{lemma}
\begin{proof}
By construction.
\end{proof}
\begin{lemma}
In generalized star-shaped complexes, the generalized cocone
operator satisfies the following identity,
\[
H \d + \d H = I,
\]
at the level of cochains.
\end{lemma}
\begin{proof}
Follows immediately from applying the proof in the trivially
star-shaped case, and using the identity in the previous lemma.
\end{proof}
Similarly, we have a discrete Poincar\'{e} lemma for generalized
star-shaped complexes.
\begin{corollary}[Discrete Poincar\'{e} Lemma for Generalized
Star-shaped Complexes] Given a closed cochain $\alpha^k$, that is
to say, $\d\alpha^k=0$, there exists a cochain $\beta^{k-1}$, such
that, $\d\beta^{k-1}=\alpha^k$.
\end{corollary}
\begin{proof}
Follows from the above lemma using the proof for the trivially
star-shaped case.
\end{proof}
\begin{example}\label{cex}\index{Poincar\'e lemma!counterexample}
We will consider an example of how the Poincar\'{e} lemma fails in
the case when the complex is not contractible. Consider the
following trivially star-shaped complex, and augment by one vertex
so as to make the region non-contractible, as show in Figure~\ref{dec:fig:non_contract}.

\begin{figure}[htbp]
\begin{center}
\subfigure[Trivially star-shaped complex]{\qquad\includegraphics[scale=0.8]{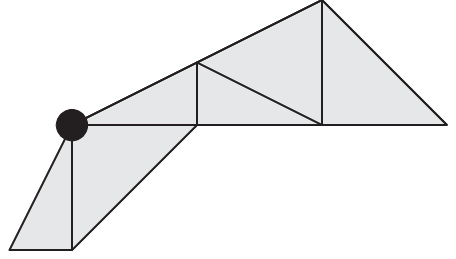}\qquad}
\subfigure[Non-contractible complex]{\qquad\includegraphics[scale=0.8]{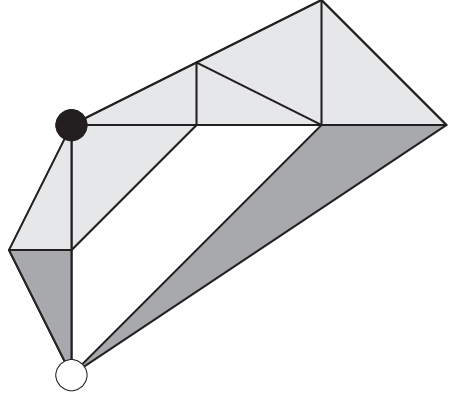}\qquad}
\end{center}
\caption{\label{dec:fig:non_contract}Counter-example for the discrete Poincar\'e lemma for a non-contractible complex.}
\end{figure}
Now we attempt to verify the identity,
\[ p\partial + \partial p = I,\]
and we will see how this is only true up to a chain that is
homotopic to the inner boundary.
\begin{align*}
(p\partial + \partial
p)\begin{pmatrix}\raisebox{-5ex}{\includegraphics[scale=0.5]{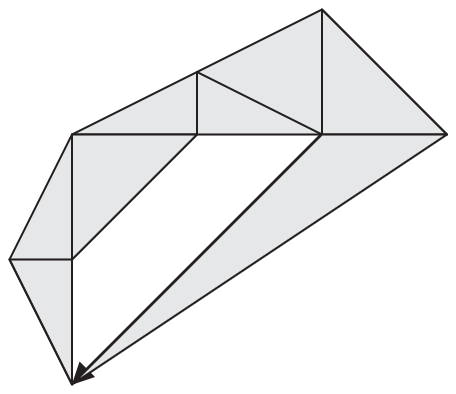}}\end{pmatrix}
&=p\begin{pmatrix}\raisebox{-6.9ex}{\includegraphics[scale=0.5]{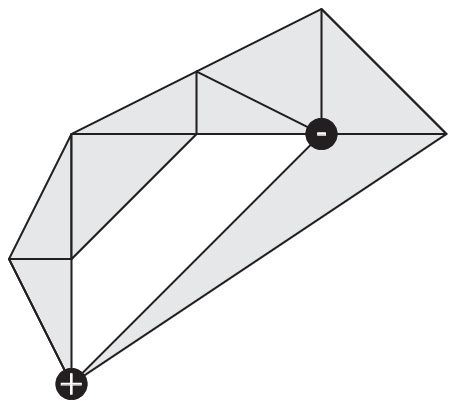}}\end{pmatrix}
+
\partial\begin{pmatrix}\raisebox{-5ex}{\includegraphics[scale=0.5]{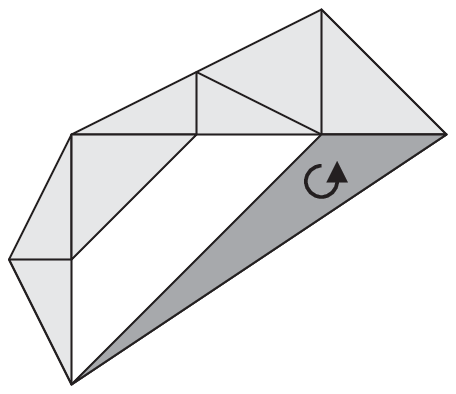}}\end{pmatrix}\\
&=\raisebox{-5ex}{\includegraphics[scale=0.5]{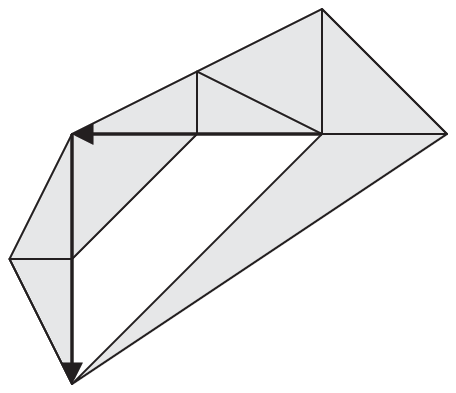}} + \raisebox{-5ex}{\includegraphics[scale=0.5]{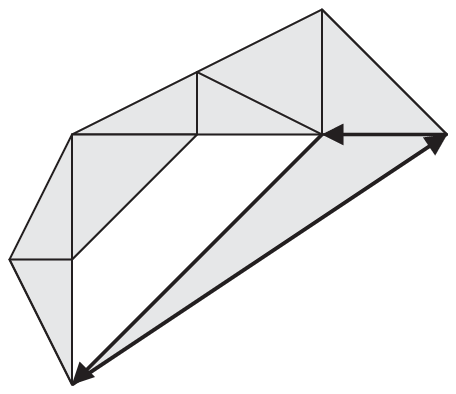}}\\
&=\raisebox{-5ex}{\includegraphics[scale=0.5]{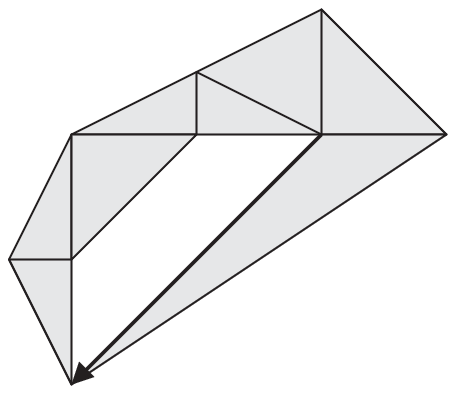}} +
\raisebox{-5ex}{\includegraphics[scale=0.5]{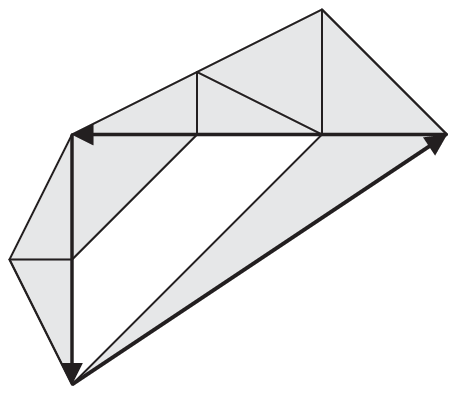}}
\end{align*}
Since the second term cannot be expressed as the boundary of a
$2$-chain, it will contribute a non-trivial effect, even on closed
discrete forms, and therefore the discrete Poincar\'{e} lemma does not hold for non-contractible complexes, as expected.
\end{example}

\section{Discrete Variational Mechanics and DEC}\index{discrete mechanics}
\label{dec:sec:commute}
We recall that discrete variational mechanics is based on a discrete analogue of Hamilton's principle, and they yield the discrete Euler--Lagrange equations. A particularly interesting property of DEC arises when it is used to construct the discrete Lagrangian for harmonic functions, and Maxwell's equations.

In particular, for these examples, the following diagram commutes,
\[
\xymatrix
@R=4pc@C=6pc{\txt{\textbf{Lagrangian}\\$L:TQ\rightarrow\mathbb{R}$}
\ar[r]^(0.42){DEC} \ar[d]&
\txt{\textbf{Discrete Lagrangian}\\$L_d:Q\times Q\rightarrow\mathbb{R}$} \ar[d]\\
\txt{\textbf{Euler--Lagrange}\\$\mathcal{EL}:T^2 Q\rightarrow T^*
Q$} \ar[r]^(0.42){DEC} & \txt{\textbf{Discrete
Euler--Lagrange}\\$\mathcal{EL}_d:Q^3\rightarrow T^* Q$}}\]

Which is to say that directly discretizing the differential
equations for harmonic functions, and Maxwell's equations using
DEC results in the same expressions as the discrete
Euler--Lagrange equations associated with a discrete Lagrangian
which is discretized from the corresponding continuous Lagrangian
by using DEC as the discretization scheme.

This is significant, since it implies that when DEC is used to discretize these equations, the corresponding numerical scheme which is obtained is variational, and consequently exhibits excellent structure-preserving properties.

In the variational principles for both harmonic functions and Maxwell's equations, we require the $L^2$ norm obtained from the $L^2$ inner product on $\Omega^k(M)$, which is given by
\[\langle\alpha^k, \beta^k \rangle = \int_M \alpha\wedge\ast\beta\, .\]
The discrete analogue of this requires a primal-dual wedge product, which is given below for forms of complementary dimension.

\begin{definition}
Given a primal discrete $k$-form $\alpha^k\in \Omega_d^k(K)$, and a dual discrete $(n-k)$-form $\hat\beta^{n-k}\in\Omega_d^{n-k}(\star K)$, the \textbf{discrete primal-dual wedge product}\index{wedge product!primal-dual} is defined as follows,
\begin{align*}
\langle \alpha^k \wedge \hat\beta^{n-k}, V_{\sigma^k}
\rangle &=
\frac{|V_{\sigma^k}|}{|\sigma^k||\star\sigma^k|}\langle\alpha^k,\sigma^k\rangle\langle\hat\beta^{n-k},\star\sigma^k\rangle\\
&=
\frac{1}{n}\langle\alpha^k,\sigma^k\rangle\langle\hat\beta^{n-k},\star\sigma^k\rangle\, ,
\end{align*}
where $V_{\sigma^k}$ is the $n$-dimensional support volume obtained by taking the convex hull of the simplex $\sigma^k$ and its dual cell $\star \sigma^k$.
\end{definition}

The corresponding $L^2$ inner product is as follows.
\begin{definition}
Given two primal discrete $k$-forms, $\alpha^k,\beta^k\in\Omega_d^k(K)$, their \textbf{discrete $L^2$ inner product}\index{inner product}, $\langle\alpha^k,\beta^k\rangle_d$, is given by
\begin{align*}
\langle \alpha^k,\beta^k\rangle_d\
&= \sum_{\sigma^k\in K}\frac{| V_{\sigma^k} |}{|\sigma^k||\star\sigma^k|}\langle \alpha^k,\sigma^k\rangle \langle\ast\beta,\star\sigma^k\rangle\,\\
&= \frac{1}{n} \sum_{\sigma^k\in K}\langle \alpha^k,\sigma^k\rangle \langle\ast\beta,\star\sigma^k\rangle
\, .
\end{align*}
\end{definition}

\begin{remark}\label{dec:rem:hodge_star_from_inner_product}
Notice that it would have been quite natural from the smooth theory to propose the following metric  tensor $\langle\!\langle\, ,\, \rangle\!\rangle$ for differential forms,
\[ \langle\, \langle\!\langle \alpha^k,\beta^k\rangle\!\rangle \mathbf{v},V_{\sigma^k} \rangle = |V_{\sigma^k}|\frac{\langle\alpha^k,\sigma^k\rangle}{|\sigma^k|}
\frac{\langle\beta^k,\sigma^k\rangle}{|\sigma^k|}\, ,\]
where the $|V_{\sigma^k}|$ is the factor arising from integrating the volume-form over $V_{\sigma^k}$, and
\[ \frac{\langle\alpha^k,\sigma^k\rangle}{|\sigma^k|}
\frac{\langle\beta^k,\sigma^k\rangle}{|\sigma^k|}\]
is what we would expect for $\langle\!\langle\alpha^k,\beta^k\rangle\!\rangle$, if the forms $\alpha^k$ and $\beta^k$ were constant on $\sigma^k$, which is the product of the average values of $\alpha^k$, and $\beta^k$.

If we adopt this as our definition of the metric tensor for forms, we can recover the definition we obtained in \S\ref{dec:sec:hodge} for the Hodge star operator. Starting from the definition from the smooth theory,
\[\int \langle\!\langle \alpha^k,\beta^k\rangle\!\rangle\mathbf{v} = \int \alpha^k\wedge\ast\beta^k\, , \]
and expanding this in terms of the metric tensor for discrete forms, and the primal-dual wedge operator, we obtain
\begin{align*}
\langle\, \langle\!\langle \alpha^k, \beta^k \rangle\!\rangle\mathbf{v}, V_{\sigma^k}\rangle
&=  \langle \alpha^k\wedge\ast\beta^k,V_{\sigma^k}\rangle\, ,\\
|V_{\sigma^k}|\frac{\langle\alpha^k,\sigma^k\rangle}{|\sigma^k|}
\frac{\langle\beta^k,\sigma^k\rangle}{|\sigma^k|}
&= \frac{|V_{\sigma^k}|}{|\sigma^k||\star\sigma^k|}\langle\alpha^k,\sigma^k\rangle\langle\ast\beta^k,\star\sigma^k\rangle\, .
\end{align*}
When we eliminate common factors from both sides, we obtain the expression,
\[\frac{1}{|\sigma^k|}\langle \beta^k,\sigma^k\rangle = \frac{1}{|\star\sigma^k|}\langle\ast\beta^k,\star\sigma^k\rangle\, ,\]
which is the expression we previously obtained in Definition~\ref{dec:defn:hodge_star} of \S\ref{dec:sec:hodge}.
\end{remark}

The $L^2$ norm for discrete differential forms is given below.

\begin{definition}\label{dec:defn:discrete_L2_norm}
Given a primal discrete $k$-form $\alpha^k\in\Omega_d^k(K)$, its \textbf{discrete $L^2$ norm}\index{norm} is given by
\begin{align*}
\| \alpha^k \|^2_d &= \langle \alpha^k,\alpha^k\rangle_d\\
&=\frac{1}{n} \sum_{\sigma^k\in K}\langle\alpha^k,\sigma^k\rangle\langle\ast\alpha^k,\star\sigma^k\rangle\\
& = \frac{1}{n}\sum_{\sigma^k\in K}\frac{|\star \sigma^k|}{|\sigma^k|}\langle\alpha^k,\sigma^k\rangle^2\, .
\end{align*}
\end{definition}

Given these definitions, we can now reproduce some computations that were originally shown in \cite{Ca2003}.

\paragraph{Harmonic Functions.}\index{Harmonic functions}
Harmonic functions $\phi:M\rightarrow\mathbb{R}$ can be characterized in a variational fashion as extremals of the following action functional,\index{harmonic functions!action functional}\index{variational problems!harmonic functions|see{harmonic functions}}
\[\mathcal{S}(\phi)=\frac{1}{2}\int_M \|\d\phi\|^2\mathbf{v},\]
where $\mathbf{v}$ is a Riemannian volume-form in $M$.
The corresponding Euler--Lagrange equation is given by\index{harmonic functions!Euler--Lagrange}
\[\ast\d\ast\d\phi=-\Delta \phi =0,\]
which is the familiar characterization of harmonic functions in terms of the Laplace--Beltrami operator.

\begin{figure}[htbp]
\begin{center}
\subfigure[Vertex $v$]{\includegraphics[scale=0.4,clip=true]{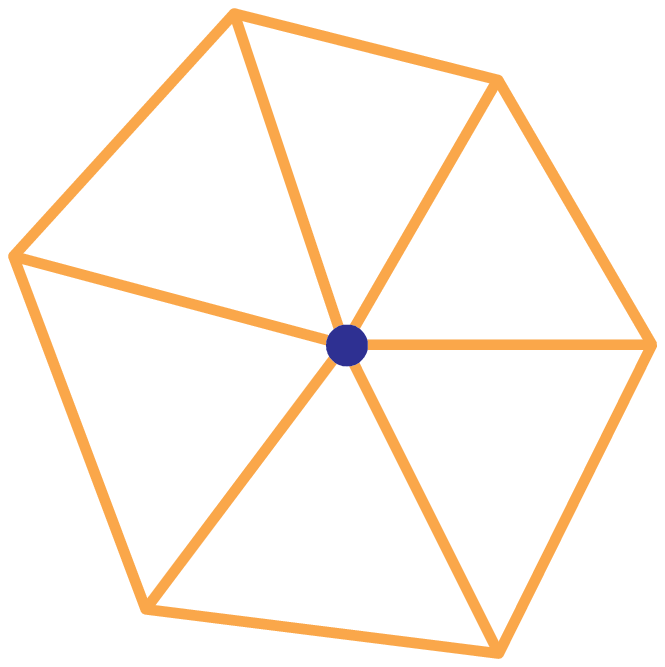}}\qquad
\subfigure[$\star v$]{\includegraphics[scale=0.4,clip=true]{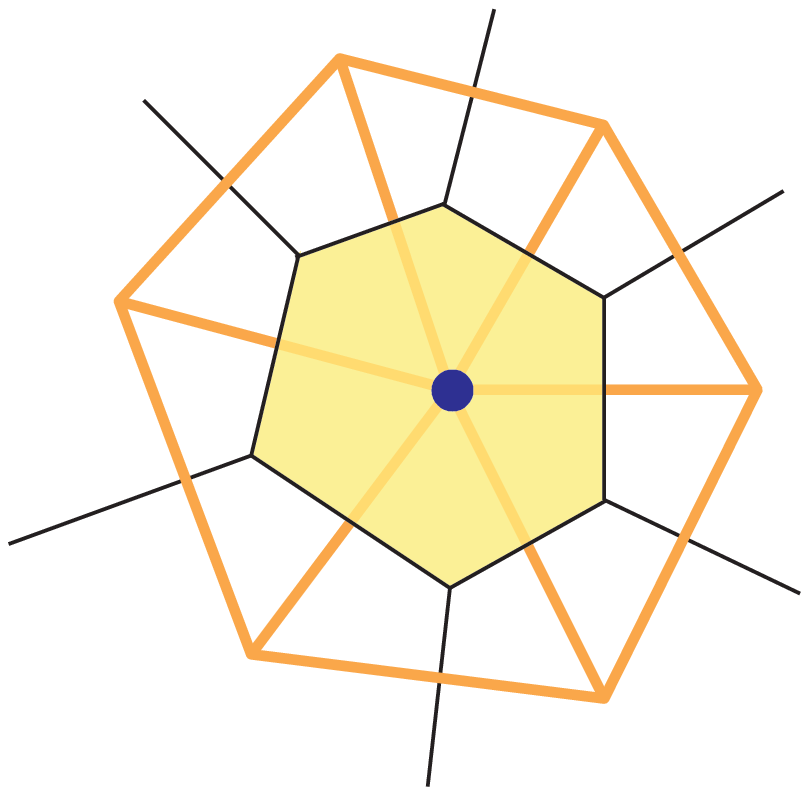}}\qquad
\subfigure[$\partial \star v$]{\includegraphics[scale=0.4,clip=true]{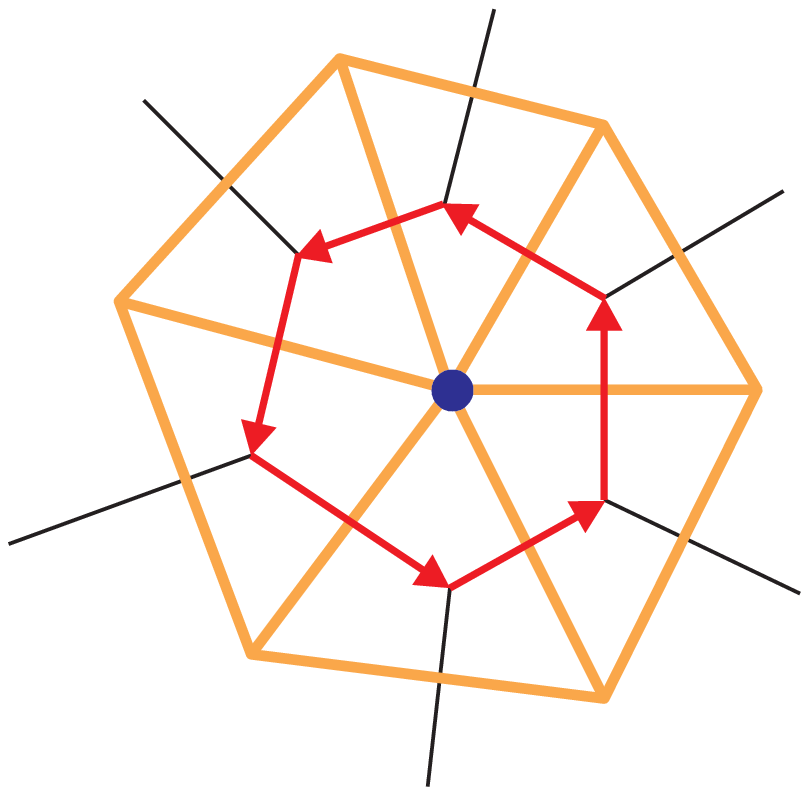}}\\
\subfigure[$\sigma^1\succ v$]{\includegraphics[scale=0.4,clip=true]{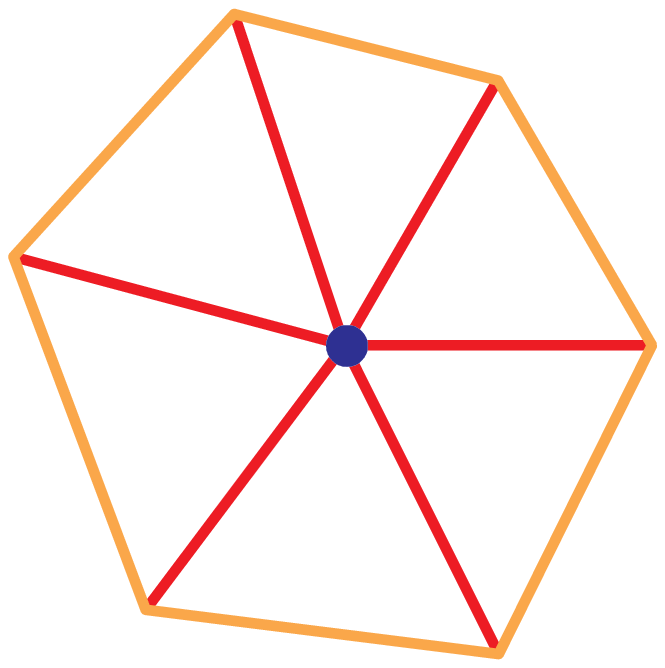}}\qquad
\subfigure[$\star\sigma^1$]{\includegraphics[scale=0.4,clip=true]{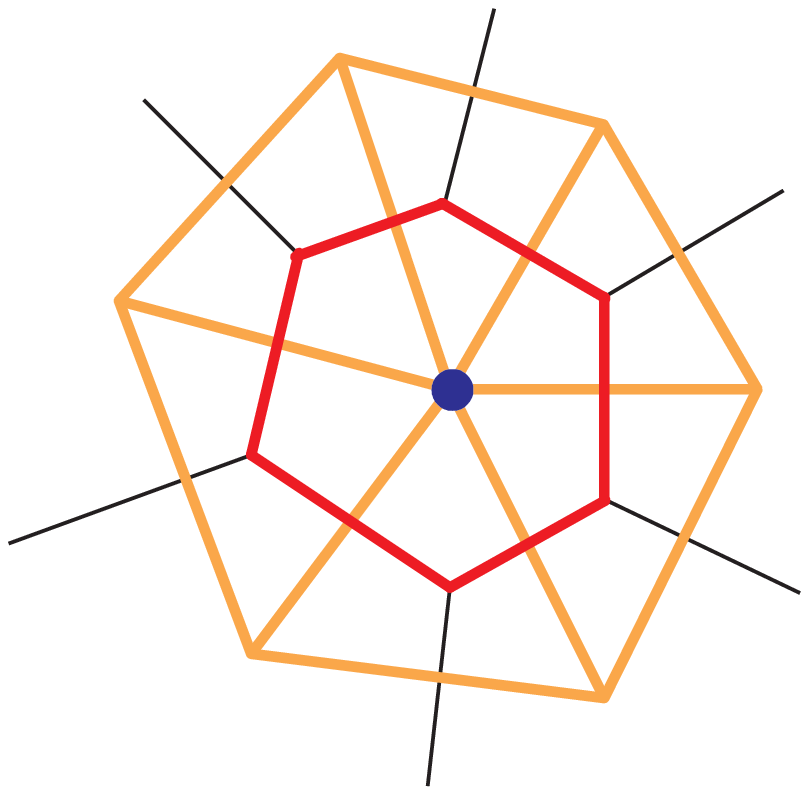}}
\end{center}
\caption{\label{dec:fig:harmonic}Boundary of a dual cell.}
\end{figure}

The discrete action functional can be expressed in terms of the $L^2$ norm we introduced above for discrete forms,
\begin{align*}
\mathcal{S}_{d}(\phi )&=\frac{1}{2}\|\d\phi\|_d^2\\
&=\frac{1}{2n}\sum_{\sigma ^{1}\in K}
\frac{\left| \star \sigma ^{1}\right| }{\left| \sigma ^{1}\right| } \langle
\d\phi,\sigma ^{1}\rangle ^{2}.
\end{align*}
The basic variations needed for the determination of the discrete
Euler--Lagrange operator are obtained from variations that vary the value of the function $\phi$ at a given vertex $v_0$, leaving the other values fixed. These variations have the form,
\[
\phi _{\varepsilon }=\phi +\varepsilon \tilde{\eta},
\]
where $\tilde{\eta}\in \Omega^0(M;\mathbb{R})$ is such that $\langle \tilde{\eta},v_{0}\rangle =1$, and $\langle \tilde{\eta},v\rangle =0$, for any $v\in K^{(0)}-\{v_{0}\}$. This family of variations
is enough to establish the variational principle. That is, we have
\begin{align}
0 &=\left. \frac{d}{d\varepsilon}\right|_{\varepsilon =0}\mathcal{S}_d(\phi _{\varepsilon})\notag\\
&=\frac{1}{n}\sum_{\sigma ^1\in K}\frac{\left| \star \sigma^1\right| }{\left| \sigma^1\right| }
\langle \d\phi,\sigma ^1\rangle \langle \d\tilde{\eta},\sigma^1\rangle\notag\\
&= \frac{1}{n}\sum _{v_0\prec \sigma^1}\frac{\left| \star \sigma^1\right| }{\left| \sigma^1\right|}\langle \d\phi,\sigma^1\rangle \operatorname{sgn}(\sigma^1;v_0),\label{dec:eq:harmonic_del}
\end{align}
where $\operatorname{sgn}(\sigma ^{1};v)$ stands for the sign of $\sigma^1$
with respect to $v$. Which is to say, $\operatorname{sgn}(\sigma^1;v)=1$ if $\sigma^1=[v^\prime ,v]$, and $\operatorname{sgn}(\sigma ^{1};v)=-1$ if $\sigma^1=[v,v^\prime]$. On the other hand,
\begin{align*}
\langle \ast \d\ast \d\phi ,v_0\rangle &=\frac{1}{|\star v_0|}
\langle \d\ast \d\phi ,\star v_0\rangle\\
&=\frac{1}{|\star v_0| }\langle \ast \d\phi ,\partial \star v_0\rangle\\
&=\frac{1}{|\star v_0| }\sum_{v_0\prec
\sigma^1}\langle \ast \d\phi ,\star \sigma^1\rangle \operatorname{sgn}(\sigma^1;v_0) \\
&=\frac{1}{|\star v_0|}\sum_{v_0\prec \sigma^1}\frac{|\star \sigma^1| }{|\sigma^1|}\langle \d\phi
,\sigma^1\rangle \operatorname{sgn}(\sigma^1;v_0),
\end{align*}
where in the second to last equality, one has to note that the
border of the dual cell of a vertex $v_{0}$ consists, up to orientation, in
the dual of all the $1$-simplices starting from $v_{0}$. This is illustrated in Figure~\ref{dec:fig:harmonic}, and follows from a general expression for the boundary of a dual cell that was given in Definition~\ref{dec:def:dual_boundary}.

The sign factor
comes from the relation between the orientation of the dual of the $1$-simplices and that of $\partial \ast v_{0}$. From this, we conclude that the variational discrete equation, given in Equation \ref{dec:eq:harmonic_del}, is equivalent to the vanishing of the discrete Laplace--Beltrami operator,
\begin{equation*}
\ast \d\ast \d\phi = -\Delta = 0.
\end{equation*}

\paragraph{Maxwell Equations.}\index{Maxwell equations}\index{variational problems!Maxwell equations|see{Maxwell equations}}
We can formulate the Maxwell equations of electromagnetism in a
covariant fashion by considering the $1$-form $A$ (the potential) as our fundamental variable in a Lorentzian manifold $X$. The action functional for a Lagrangian formulation of electromagnetism is given by,\index{Maxwell equations!action functional}
\[ \mathcal{S}(A)=\frac{1}{2} \int_X \| \d A\|^2 \mathbf{v}\, ,\]
where $\|\cdot\|$ is the norm on forms induced by the Lorentzian metric on $X$, and $\mathbf{v}$ is the pseudo-Riemannian volume-form. The $1$-form $A$ is related to the $4$-vector potential encountered in the relativistic formulation of electromagnetism (see, for example~\cite{Jackson1998}).

The Euler--Lagrange equation corresponding to this action functional is given by
\[ \ast\d\ast\d A=0\, .\]
In terms of the field strength, $F=\d A$, the last equation is usually rewritten as\index{Maxwell equations!Euler--Lagrange}
\[\d F=0\, ,\qquad \ast\d\ast F =0\, ,\]
which is the geometric formulation of the Maxwell equations.

For the purposes of simplicity of exposition, we consider the special case where the Lorentzian manifold decomposes into $X=M\times\mathbb{R}$, where $(M,g)$ is a compact Riemannian $3$-manifold. In formulating the discrete version of this variational problem, we need to generalize the notion of a discrete Hodge dual to take into account the pseudo-Riemannian metric structure. This can be subtle in practice, and to overcome this, we consider a special family of complexes instead.

\begin{figure}[htbp]
\WARMprocessMoEPS{prism_new}{eps}{bb}
\renewcommand{\xyWARMinclude}[1]{\scaledfig{.5}{#1}}
\begin{center}
\leavevmode
\begin{xy}
\xyMarkedImport{}
\xyMarkedTextPoints{1-2}
\end{xy}
\end{center}
\renewcommand{\xyWARMinclude}[1]{\includegraphics{#1}}
\caption{\label{dec:fig:prism}Prismal cell complex decomposition of space-time.}
\end{figure}
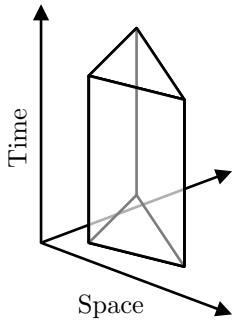

Let $K^{\prime }$ be a simplicial complex modelling $M$. For the sake of simplicity we consider $M=\mathbb{R}^3$ although this is not strictly necessary. We now consider a discretization $%
\{t_n\}_{n\in \mathbb{Z}}$ of $\mathbb{R}$. We define the complex $K$, modelling $X=\mathbb{R}^4$, the cells of which are the sets $\sigma =\sigma^\prime\times \{t_n\}\subset \mathbb{R}^3 \times \mathbb{R}$, and $\sigma =\sigma^\prime\times (t_n,t_{n+1})\subset \mathbb{R}^3 \times \mathbb{R}$ for any $\sigma^\prime\in K^\prime$ and $n\in \mathbb{Z}$. Of course,
this is not a simplicial complex but rather a ``prismal'' complex\index{complex!prismal}, as shown in Figure~\ref{dec:fig:prism}.

The advantage of these cell complexes is the existence of the
Voronoi dual. More precisely, given any prismal cells $\sigma^\prime\times\{t_n\}\in K$ and $\sigma^\prime\times (t_n,t_{n+1})\in K$, the Lorentz orthonormal to any of its edges coincide
with the Euclidean one in $\mathbb{R}^{4}$ and the existence of the
circumcenter is thus guaranteed. In other words, the Lorentz circumcentric dual $%
\star K$ to $K$ is the same as the Euclidean one in $\mathbb{R}^{4}$.

\begin{remark}
Much of the construction above can be carried out more generally by considering arbitrary cell complexes in $\mathbb{R}^{4}$ that are not necessarily prismal, as long as none of its $1$-cells are lightlike. This causality condition is necessary to ensure that the circumcentric dual complex is well-behaved. However, it is sufficient for computational purposes that the complex is well-centered, in the sense that the Lorentzian circumcenter of each cell is contained inside the cell. These issues will be addressed in future work.\end{remark}

Recall that the Hodge star $\ast$ is uniquely defined by satisfying the following expression,
\[ \alpha\wedge \ast \beta = \langle\!\langle \alpha, \beta \rangle\!\rangle \mathbf{v}\, ,\]
for all $\alpha,\beta\in\Omega^k(X)$. The upshot of this is that
the Hodge star operator depends on the metric, and since we have a
pseudo-Riemannian metric, there is a sign that is introduced in
our expression for the discrete Hodge star
(Definition~\ref{dec:defn:hodge_star}) that depends on whether the
cell it is applied to is either spacelike or timelike. The
discrete Hodge star for prismal complexes in Lorentzian space is
given below.

\begin{definition}
The \textbf{discrete Hodge star for prismal complexes in
Lorentzian space}\index{Hodge star!Lorentizan space} is a map
$\ast :\Omega _d^k (K)\rightarrow \Omega_d^k (\ast K)$ defined by
giving its action on cells in a prismal complex as follows,
\[
\frac{1}{|\star\sigma^k|}\langle \ast \alpha^k ,\star \sigma^k\rangle =
\kappa (\sigma^k)\frac{1}{|\sigma^k| }\langle \alpha^k
,\sigma^k\rangle ,
\end{equation*}
where $|\cdot| $ stands for the volume and the \textbf{causality sign}\index{causality sign}\index{$\kappa$|see{causality sign}} $%
\kappa (\sigma^k)$ is defined to be $+1$ if all the edges of $\sigma^k$
are spacelike, and $-1$ otherwise.
\end{definition}

The causality sign of $2$-cells in a $(2+1)$-space-time is summarized in Table~\ref{dec:table:causality}. We should note that the causality sign for a $0$-simplex, $\kappa(\sigma^0)$, is always 1. This is because a $0$-simplex has no edges, and as such the statement that \textit{all} of its edges are spacelike is trivially true.

\begin{table}[h!]
\caption{\label{dec:table:causality}\index{causality sign!example}Causality sign of $2$-cells in a $(2+1)$-space-time.}
\begin{center}
\begin{tabular}{|c|c|c|c|c|c|}
  \hline
  \raisebox{8ex}[0pt]{$\sigma^2$} &
  \includegraphics[scale=0.35]{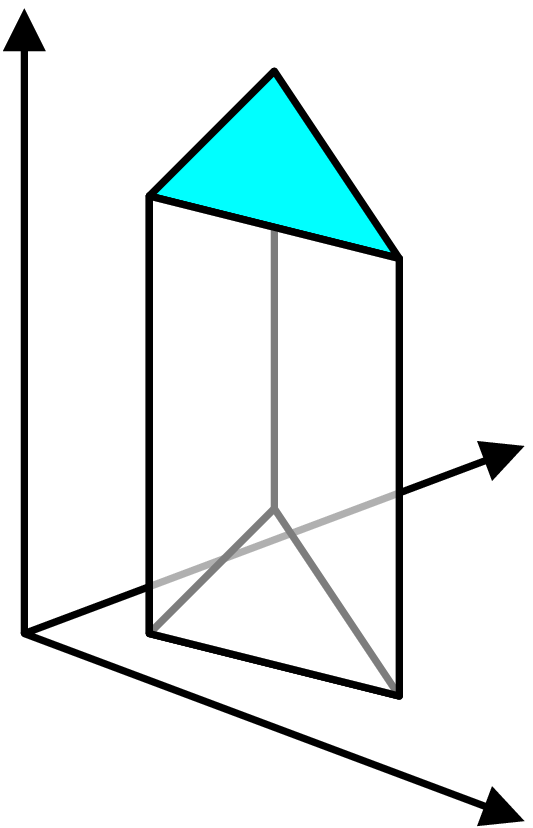} &
  \includegraphics[scale=0.35]{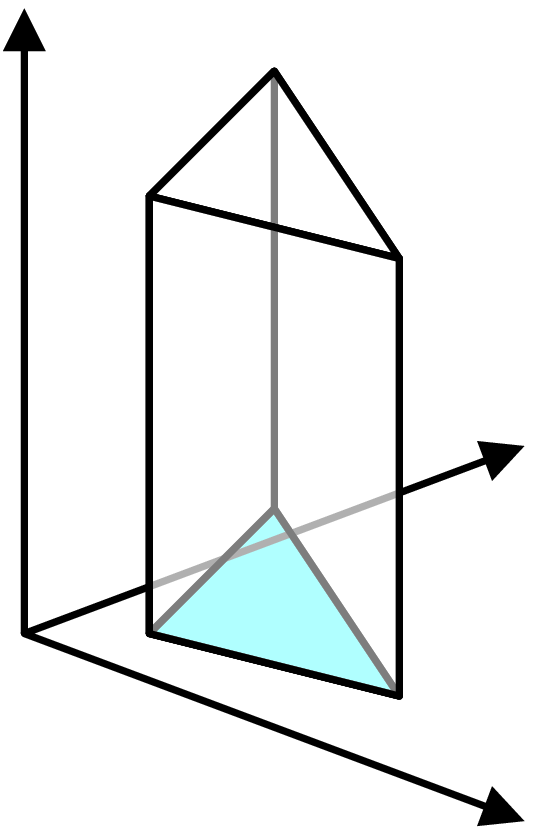} &
  \includegraphics[scale=0.35]{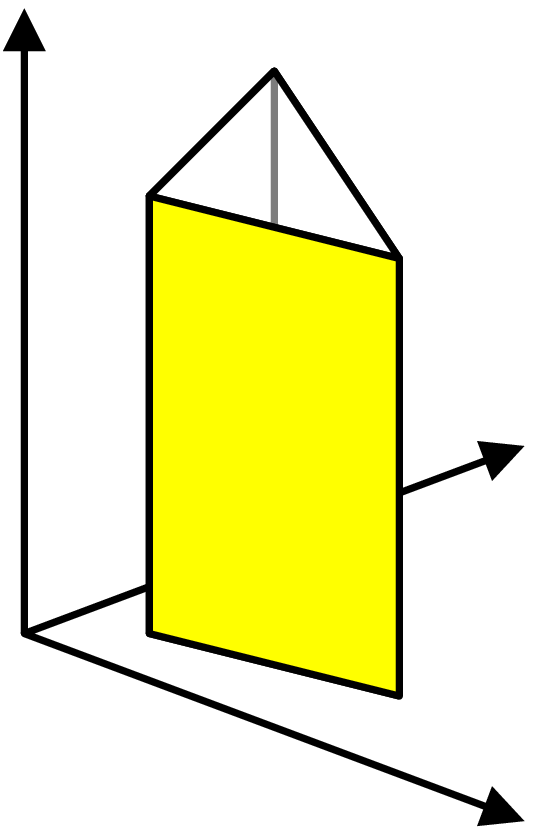} &
  \includegraphics[scale=0.35]{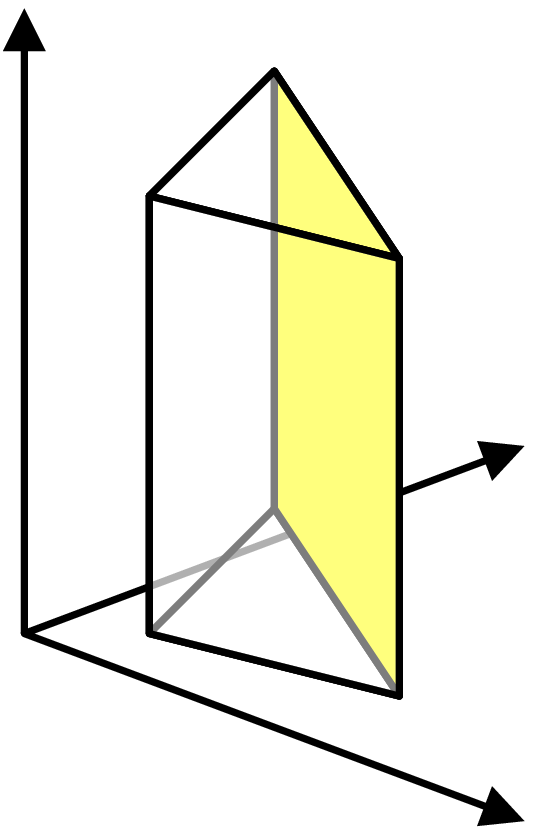} &
  \includegraphics[scale=0.35]{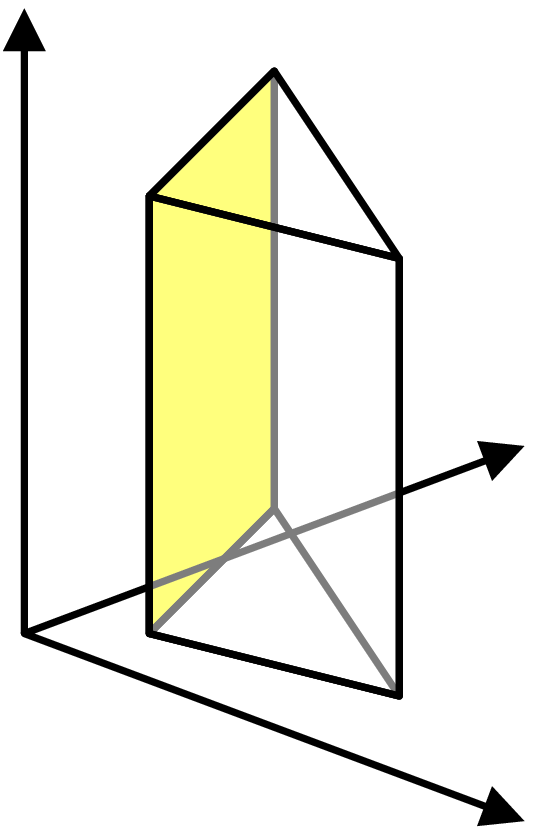}
  \\ \hline
  $\kappa(\sigma^2)$ & +1 & +1 & $-1$ & $-1$ & $-1$\\ \hline
\end{tabular}
\end{center}
\end{table}

This causality term in the discrete Hodge star has consequences for the expression for the discrete norm (Definition~\ref{dec:defn:discrete_L2_norm}), which is now given.

\begin{definition}
Given a primal discrete $k$-form $\alpha^k\in\Omega_d^k(K)$, its \textbf{discrete $L^2$ Lorentzian norm}\index{norm!Lorentzian} is given by,
\begin{align*}
\|\alpha^k\|^2_{\operatorname{Lor},d} &=\frac{1}{n} \sum_{\sigma^k\in K}\langle\alpha^k,\sigma^k\rangle\langle\ast\alpha^k,\star\sigma^k\rangle\\
& = \frac{1}{n}\sum_{\sigma^k\in K}\kappa(\sigma^k)\frac{|\star\sigma^k|}{|\sigma^k|}\langle\alpha^k,\sigma^k\rangle^2\, .
\end{align*}
\end{definition}

Having defined the discrete Lorentzian norm, we can express the discrete action as
\begin{align*}
\mathcal{S}_d(A) &= \frac{1}{2}\| \d A\|^2_{\operatorname{Lor},d}\\
&= \frac{1}{8}\sum_{\sigma^2\in K} \langle \d A, \sigma^2 \rangle \langle \ast\d A, \star \sigma^2 \rangle\\
&= \frac{1}{8}\sum_{\sigma^2\in K} \kappa(\sigma^2) \frac{|\star\sigma^2|}{|\sigma^2|}\langle \d A, \sigma^2 \rangle^2\, .
\end{align*}
The basic variations needed to determine the discrete Euler--Lagrange operator are obtained from variations that vary the value of the $1$-form $A$ at a given $1$-simplex $\sigma^1_0$, leaving the other values fixed. These variations have the form,
\[
A_\varepsilon=A_\varepsilon+\varepsilon \tilde{\eta},
\]
where $\tilde{\eta}\in \Omega^1_d (K)$ is given by $\langle \tilde{\eta}
,\sigma^1_0\rangle =1$ for a fixed interior $\sigma^1_0\in K$
and $\langle \tilde{\eta},\sigma^1\rangle =0$ for $\sigma^1\neq \sigma
_0^1$. The derivation of the variational principle gives
\begin{align*}
\left. \frac{d}{d\varepsilon}\right| _{\varepsilon=0}\mathcal{S}_d(A_{\varepsilon })
&=\frac{1}{4}\sum_{\sigma^2\in K}\frac{|\star\sigma^2|}
{|\sigma^2|}\kappa (\sigma^2)\langle \d A,\sigma^2\rangle \langle \d\tilde{\eta},\sigma^2\rangle\\
&=\frac{1}{4}\sum_{\sigma_0^1\prec \sigma^2}\frac{|\star\sigma^2|}{|\sigma^2|}\kappa(\sigma ^2)\langle \d A,\sigma^2\rangle \langle \d\tilde{\eta},\sigma^2\rangle \\
&=\frac{1}{4}\sum_{\sigma _0^1\prec \sigma^2}\frac{|\star\sigma^2|}{|\sigma^2|}\kappa(\sigma ^2)\langle \d A,\sigma^2\rangle \operatorname{sgn}(\sigma^2,\sigma_0^1),
\end{align*}
which vanishes for all the basic variations above.
On the other hand, we now expand the discrete $1$-form $\ast \d\ast \d A$. For
any $\sigma_0^1\in K$, we have that
\begin{align*}
\langle \ast \d\ast \d A,\sigma _0^1\rangle &=\frac{| \sigma_0^1| }{| \star \sigma_0^1| }\kappa (\sigma_0^1)\langle \d\ast \d A,\star \sigma_0^1\rangle\\
&=\frac{| \sigma_0^1| }{| \star \sigma _0^1| }\kappa (\sigma_0^1)\langle \ast \d A,\partial \star \sigma_0^1\rangle\\
&=\frac{| \sigma_0^1| }{|\star\sigma_0^1| }\kappa (\sigma_0^1)\sum_{\sigma_0^1\prec
\sigma^2}\operatorname{sgn}(\sigma^2,\sigma_0^1)\langle \ast \d A,\star\sigma^2\rangle\\
&=\frac{| \sigma_0^1| }{|\star \sigma_0^1| }\kappa (\sigma_0^1)\sum_{\sigma_0^1\prec \sigma ^2}\frac{|\star \sigma^2| }{| \sigma^2| }\kappa (\sigma^2)\langle \d A,\sigma^2\rangle \operatorname{sgn}(\sigma^2,\sigma_0^1),
\end{align*}
where the sign $\operatorname{sgn}(\sigma^2,\sigma^1)$ stands for the
relative orientation between $\sigma^2$ and $\sigma^1$. Which is to say, $\operatorname{sgn}(\sigma ^2,\sigma^1)=1$ if the orientation induced by $\sigma^2$ on $\sigma^1$ coincides with the orientation of $\sigma^1$, and $\operatorname{sgn}(\sigma^2,\sigma^1)=-1$ otherwise. For the second to last equality, one has to note that the border of the
dual cell of an edge $\sigma_0^1$ consists, conveniently oriented with
the $\operatorname{sgn}$ operator, of the union of the duals of all the $2$-simplices containing $\sigma _0^1$. This statement is the content of Definition~\ref{dec:def:dual_boundary}, which gives the expression for the boundary of a dual cell, and was illustrated in Figure~\ref{dec:fig:harmonic} for the case of $n$-dimensional dual cells.

By comparing the two computations, we find that for an arbitrary choice of $\sigma_0^1\in K$, $\langle \ast\d\ast \d A,\sigma^1_0\rangle$ is equal (up to a non-zero constant) to $\delta \mathcal{S}_d(A)$, which always vanishes. It follows that the variational discrete equations obtained above is equivalent to the discrete Maxwell equations,
\begin{equation*}
\ast \d\ast \d A=0.
\end{equation*}

\section{Extensions to Dynamic Problems}\index{dynamic problems}

It is desirable to leverage the exactness properties of the operators of discrete exterior calculus to construct numerical algorithms with discrete conservation properties. For these purposes, it is appropriate to extend the scope of DEC to incorporate dynamical behavior, by addressing the issue of discrete diffeomorphisms and flows.

As discussed in the previous section, DEC and discrete mechanics have interesting synergistic properties, and in this section we will explore a groupoid interpretation of discrete mechanics that is particularly appropriate to formulating the notion of pull-back and push-forward of discrete differential forms.

\subsection{Groupoid Interpretation of Discrete Variational
Mechanics}\index{discrete mechanics!groupoid}\index{groupoid!discrete mechanics|see{discrete mechanics, groupoid}}

The groupoid formulation of discrete mechanics is particularly fruitful and natural, and it serves as a unifying tool for understanding the variational formulation of discrete Lagrangian mechanics, and discrete Euler--Poincar\'e reduction, as discussed in the work of \cite{We1996} and \cite{MaPeSh1999, MaPeSh2000}.

The groupoid interpretation of discrete mechanics is most clearly illustrated if we consider the discretization of trajectories on $TQ$ in two stages. Given a curve
$\gamma:\mathbb{R}^{+}\rightarrow TQ$, we consider a discrete
sampling given by
\[
g_i=\gamma(ih)\in TQ.
\]
We then approximate $TQ$ by $Q\times Q$, and associate to $g_i$
two elements in $Q$. We denote this by
\[
g_i\mapsto (q_i^0,q_i^1).
\]
Or equivalently, in the language of groupoids, see \cite{CaWe1999,
We2001}, we have
\[\NoTips
\xymatrix{G \ar@/_/@<-1ex>[d]_\alpha \ar@/^/@<1ex>[d]^\beta\\
Q\\}
\]
where $\alpha$ is the {\bfi source}\index{groupoid!source} map,
and $\beta$ is the {\bfi target}\index{groupoid!target} map. Then,
\begin{align*}
g_i\mapsto (\alpha(g_i),\beta(g_i))=(q_i^0,q_i^1).
\end{align*}
This can be visualized as
\[\NoTips\begin{xy}
0*{q_i^0=\alpha(g_i)}, (30,0)*{q_i^1=\beta(g_i)},
(0,4)*{\bullet}="a", (30,4)*{\bullet}, \ar@{<-}@/_20pt/"a"_{g_i}
\end{xy}\]
A product $\cdot:G^{(2)}\rightarrow G$ is defined on the set of
composable pairs,
\begin{align*}
G^{(2)}:=\{(g,h)\in G\times G \mid \beta(g)=\alpha(h)\}.
\end{align*}
The {\bfi groupoid composition} $g\cdot
h$\index{groupoid!composition} is defined by
\begin{align*}
\alpha(g\cdot h)&=\alpha(g),\\
\beta(g\cdot h)&=\beta(h).
\end{align*}
This can be represented graphically as follows,
\[\NoTips\begin{xy} <2cm,0cm>:
(0,0)*+@{*}="a"*+!U{\alpha(g)=\alpha(g\cdot h)},
(2,0)*+@{*}="b"*+!U{\beta(g)=\alpha(h)},
(4,0)*+@{*}="c"*+!U{\beta(h)=\beta(g\cdot h)}, \ar@/^6ex/_{g}
"a";"b", \ar@/^6ex/_{h} "b";"c", \ar@<1ex>@/^12ex/^{g\cdot h}
"a";"c"
\end{xy}\]
The set of composable pairs is the discrete analogue of the set of second-order curves on $TQ$. A curve $\gamma:\mathbb{R}^{+}\rightarrow TQ$ is said to be second-order
if there exists a curve $q:\mathbb{R}^{+}\rightarrow Q$, such that,
\[
\gamma(t)=(q(t),\dot{q}(t)).
\]
The corresponding condition for discrete curves is that given a
sequence of points in $Q\times Q$,
$(q_1^0,q_1^1),\ldots,(q_p^0,q_p^1)$, we require that
\[
q_i^1=q_{i+1}^0.
\]
This implies that the discrete curve on $Q\times Q$ is derived
from a $(p+1)$-pointed curve $(q_0,\ldots,q_p)$ on $Q$, where
\[
q_i=
\begin{cases}
q_{i+1}^0, & \text{if $0\le i < p$;}\\
q_i^1,     & \text{if $i=p$.}
\end{cases}
\]
This condition has a direct equivalent in groupoids,
\[
\beta(g_i)=q_i^1=q_{i+1}^0=\alpha(g_{i+1}).
\]
Which is to say that the sequence of points in $Q\times Q$ are composable. In general, this hierarchy of sets is denoted by
\[
G^{\left(  p\right)  } :=\left\{  \left(
g_{1},\ldots,g_{p}\right)\in G^p  \mid \beta\left(  g_{i}\right)
=\alpha\left(  g_{i+1}\right)  \right\},
\]
where $G^{(0)}\simeq Q$.

In addition, the {\bfi groupoid inverse}\index{groupoid!inverse}
is defined by the following,
\begin{align*}
\alpha(g^{-1})&=\beta(g),\\
\beta(g^{-1})&=\alpha(g).
\end{align*}
This is represented as follows,
\[\NoTips\begin{xy}
0*{\beta(g^{-1})=\alpha(g)}, (50,0)*{\alpha(g^{-1})=\beta(g)},
(10,4)*+@{*}="a", (40,4)*+@{*}="b", \ar@/^20pt/^{g}"a";"b",
\ar@{-->}@/^20pt/^{g^{-1}}"b";"a"
\end{xy}\]

\begin{figure}[htbp]
\WARMprocessMoEPS{groupoid_new}{eps}{bb}
\begin{center}
\leavevmode
\begin{xy}
\xyMarkedImport{}
\xyMarkedTextPoints{1-8}
\end{xy}
\end{center}
\caption{\label{dec:fig:groupoids}Groupoid composition and inverses.}
\end{figure}
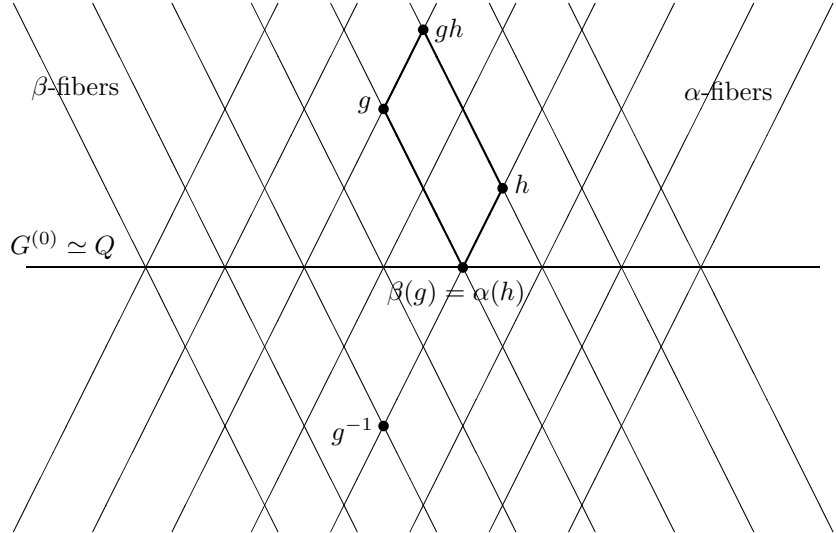

\paragraph{Visualizing Groupoids.}\index{groupoid!visualizing}
In summary, composition of groupoid elements, and the inverse of
groupoid elements can be illustrated by Figure~\ref{dec:fig:groupoids}.
As we will see in the next subsection, representing discrete diffeomorphisms as pair groupoids is the natural method of ensuring that the mesh remains nondegenerate.

\subsection{Discrete Diffeomorphisms and Discrete Flows}
We will adopt the point of view of representing a discrete diffeomorphism as a groupoid, which was first introduced in \cite{PeWe2003}, and appropriately modify it to reflect the simplicial nature of our mesh. In addition, we will address the induced action of a discrete diffeomorphism on the dual mesh.

\begin{definition}Given a complex $K$ embedded in $V$, and its
corresponding abstract simplicial complex $M$, a \textbf{discrete
diffeomorphism},\index{diffeomorphism} $\varphi\in\operatorname{Diff}_d(M)$, is a pair of simplicial complexes $K_1$, $K_2$, which are realizations of $M$ in
the ambient space $V$. This is denoted by $\varphi(M)=(K_1,K_2)$.
\end{definition}
\begin{definition}
A \textbf{one-parameter family of discrete diffeomorphisms}\index{diffeomorphism!one-parameter family} is a
map $\varphi:I\rightarrow \operatorname{Diff}_d(M)$, such that,
\[\pi_1(\varphi(t))=\pi_1(\varphi(s)),\qquad \forall s,t\in I.\]
\end{definition}
Since we are concerned with evolving equations represented by
these discrete diffeomorphisms, and mesh degeneracy causes the
numerics to fail, we introduce the notion of non-degenerate
discrete diffeomorphisms,
\begin{definition}
A \textbf{non-degenerate discrete diffeomorphism}\index{diffeomorphism!non-degenerate} $\varphi=(K_1,
K_2)$ is such that $K_1$ and $K_2$ are non-degenerate realizations
of the abstract simplicial complex $M$ in the ambient space $V$.
\end{definition}

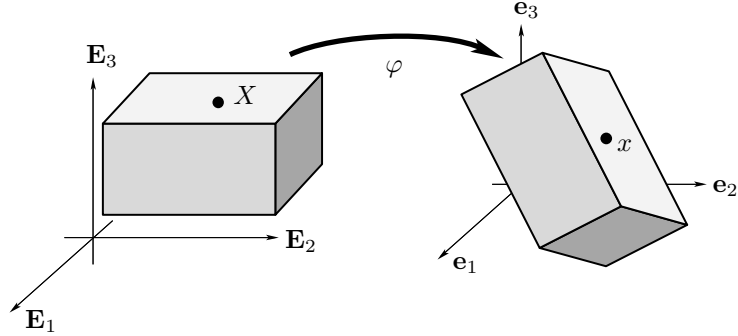
\begin{figure}[htbp]
\WARMprocessMoEPS{spatial_material}{eps}{bb}
\begin{center}
\leavevmode
\begin{xy}
\xyMarkedImport{}
\xyMarkedMathPoints{1-9}
\end{xy}
\end{center}
\caption{\label{dec:fig:spatial_material}Spatial and material representations.}
\end{figure}

Notice that it is sufficient to define the discrete diffeomorphism
on the vertices of the abstract complex $M$, since we can extend
it to the entire complex by the relation
\[
  \varphi([v_0,...,v_k]) = ([\pi_1\varphi(v_0),...,\pi_1\varphi(v_k)],[\pi_2\varphi(v_0),...,\pi_2\varphi(v_k)]).
\]
If $X\in K^{(0)}$ is a material vertex of the manifold,
corresponding to the abstract vertex $w$, that is to say,
$\pi_1\varphi_t(w)=X, \forall t\in I$, the corresponding
trajectory followed by $X$ in space is $x=\pi_2\varphi_t(w)$. Then,
the {\bfi material velocity}\index{velocity!material} $V(X,t)$ is
given by
\[
V(\pi_1(w),t)=\left.\frac{\partial\pi_2\varphi_s(w)}{\partial
s}\right|_{s=t},\] and the {\bfi spatial
velocity}\index{velocity!spatial} $v(x,t)$ is given by
\[v(\pi_2(w),t)=V(\pi_1(w),t)=\left.\frac{\partial\varphi_s(\varphi^{-1}_t(x))}{\partial
s}\right|_{s=t}.
\]
The distinction between the spatial and material representation is illustrated in Figure~\ref{dec:fig:spatial_material}.

The material velocity field can be thought of as a discrete vector
field with the vectors based at the vertices of $K$, which is to
say that $T\varphi_t\in\mathfrak{X}_d(K)$, is a discrete primal
vector field. Notice that $\varphi_t$ on $K$ induces a map
$\star\varphi_t$ on the vertices of the dual $\star K$, by the
following,
\[\star\varphi_t(c[v_0,\ldots,v_n])=(c[\pi_1\varphi_t(v_0),\ldots,\pi_1\varphi_t(v_n)],c[\pi_2\varphi_t(v_0),\ldots,\pi_2\varphi_t(v_n)]).\]
Similarly then, $T\star\varphi_t\in\mathfrak{X}_d(\star K)$ is a
discrete dual vector field.

\paragraph{Comparison with Interpolatory Methods.}\index{diffeomorphism!interpolatory methods}
At first glance, the groupoid formulation seems like a cumbersome way to define a
one-parameter family of discrete diffeomorphisms, and one may be
tempted to think of extending $\varphi_t$ to the ambient space. We
would then be thinking of $\varphi_t:V\rightarrow V$. This is
undesirable since given $\varphi_t$ and $\psi_s$ which are
non-degenerate flows, their composition $\varphi_t\circ\psi_s$,
which is defined, may result in a degenerate mesh when applied to
$K$. Thus, non-degenerate flows are not closed under this notion
of composition.

If we adopt groupoid composition instead at the level of vertices,
we can always be sure that if we compose two nondegenerate
discrete diffeomorphisms, they will remain a nondegenerate
discrete diffeomorphism.

\paragraph{Discrete Diffeomorphisms as Pair Groupoids.}\index{diffeomorphism!groupoid}\index{groupoid!diffeomorphism|see{groupoid, diffeomorphism}}
The space of discrete diffeomorphisms naturally has the structure
of a pair groupoid. The discrete analogue of
$T\operatorname{Diff}(M)$ from the point of view of temporal
discretization is the pair groupoid $\operatorname{Diff}(M)\times
\operatorname{Diff}(M)$. In addition, we discretize
$\operatorname{Diff}(M)$ using $\operatorname{Diff}_d(M)$, which
is in turn a pair groupoid involving realizations of an abstract
simplicial complex in an ambient space.

\subsection{Push-Forward and Pull-Back of Discrete Vector Fields and Discrete Forms}
For us to construct a discrete theory of exterior calculus that admits dynamic problems, it is critical that we introduce the notion of push-forward and pull-back of discrete vector fields and discrete forms under a discrete flow.

\paragraph{Push-Forward and Pull-Back of Discrete Vector Fields.}

The push-forward of a discrete vector field satisfies the following commutative diagram,
\[\xymatrix{
K \ar[r]^\star \ar[d]^f & \star K \ar[r]^X \ar[d]^{\star f}& \mathbb{R}^N \ar[d]^{Tf}\\
L \ar[r]^\star & \star L \ar[r]^{f_\ast X} & \mathbb{R}^N
}\]
and the pull-back satisfies the following commutative diagram,
\[\xymatrix{
K \ar[r]^\star \ar[d]^f & \star K \ar[r]^{f^\ast X} \ar[d]^{\star f}& \mathbb{R}^N \ar[d]^{Tf}\\
L \ar[r]^\star & \star L \ar[r]^{X} & \mathbb{R}^N
}\]
By appropriately following the diagram around its boundary, we obtain the following expressions for the push-forward and pull-back of a discrete vector field.

\begin{definition}
The \textbf{push-forward of a dual discrete vector field}\index{push-forward!vector field} $X\in\mathfrak{X}_d(\star K)$, under the map $f: K\rightarrow L$, is given by its evaluation on a dual vertex $\hat\sigma_0=\star\sigma^n\in (\star L)^{(0)}$,
\[
  f_\ast X(\star \sigma^n) = Tf\cdot X(\star(f^{-1}(\sigma^n))).
\]
\end{definition}

\begin{definition}
The \textbf{pull-back of a dual discrete vector field}\index{pull-back!vector field} $X\in\mathfrak{X}_d(\star L)$, under the map $f: K\rightarrow L$, is given by its evaluation on a dual vertex $\hat\sigma_0=\star\sigma^n\in (\star K)^{(0)}$,
\[
  f^\ast X(\star \sigma^n) = (f^{-1})_\ast X(\star \sigma^n)= T(f^{-1})\cdot X(\star(f(\sigma^n))).
\]
\end{definition}

\paragraph{Pull-Back and Push-Forward of Discrete Forms.}
A natural operation involving exterior calculus in the context of dynamic problems is the pull-back of a differential form by a flow. We define the pull-back of a discrete form as follows.

\begin{definition}
The \textbf{pull-back of a discrete form}\index{pull-back!form} $\alpha^k \in \Omega_d^k(L)$, under the map $f: K
\rightarrow L$, is defined so that the \textbf{change of variables formula}\index{change of variables} holds,
\[
  \langle f^* \alpha^k, \sigma^k\rangle = \langle \alpha^k, f (\sigma^k)
  \rangle,
\]
where $\sigma^k \in K$.
\end{definition}

We can define the push-forward of a discrete form as its pull-back under the inverse map as follows.

\begin{definition}
The \textbf{push-forward of a discrete form}\index{push-forward!form} $\alpha^k \in \Omega_d^k(K)$, under the map $f: K
\rightarrow L$ is defined by its action on $\sigma^k\in L$,
\[
  \langle f_* \alpha^k, \sigma^k\rangle = \langle (f^{-1})^\ast \alpha^k, \sigma^k\rangle
  = \langle \alpha^k, f^{-1} (\sigma^k) \rangle.
\]
\end{definition}

\paragraph{Naturality under Pull-Back of Wedge Product.}\index{wedge product!naturality} We find that the discrete wedge product we introduced in \S\ref{dec:sec:wedge} is not natural under pull-back, which is to say that the relation
\[f^\ast (\alpha\wedge\beta) = f^\ast \alpha \wedge f^\ast\beta\, ,\]
does not hold in general. However, a metric independent definition that is natural under pull-back was proposed in \cite{Ca2003}.

\begin{definition}[\cite{Ca2003}]\label{dec:def:wedge_alternative}
Given a primal discrete $k$-form $\alpha^k \in \Omega^k_d(K)$, and
a primal discrete $l$-form $\beta^l \in \Omega^l_d(K)$, the
\textbf{natural discrete primal-primal wedge product},\index{wedge product!natural} $\wedge :
\Omega^k_d(K) \times \Omega^l_d(K) \rightarrow \Omega^{k+l}_d(K)$, is
defined by its evaluation on a $(k+l)$-simplex $\sigma^{k+l} =
[v_0, \ldots, v_{k+l}]$,
\[
\langle \alpha^k \wedge \beta^l, \sigma^{k+l} \rangle =
\frac{1}{(k+l+1)!} \sum_{\tau \in S_{k+l+1}}
\operatorname{sign}(\tau)
    \alpha \smile \beta (\tau(\sigma^{k+l})) \, .
\]
\end{definition}

In contrasting this definition to that given by Definition~\ref{dec:def:wedge}, we see that the geometric factor
\[ \frac{|\sigma^{k+l}\cap\star v_{\tau(k)}|}{|\sigma^{k+l}|}\, , \]
has been replaced by
\[ \frac{1}{k+l+1}\]
in this alternative definition. By replacing the geometric factor which is metric dependent with a constant factor, Definition~\ref{dec:def:wedge_alternative} becomes natural under pull-back.

The proofs in \S\ref{dec:sec:wedge} that the discrete wedge product is anti-commutative, and satisfies a Leibniz rule, remain valid for this alternative discrete wedge product, with only trivial modifications. As for the proof of the associativity of the wedge product for closed forms, we note the following identity,
\[ \sum_{\tau\in S_{k+l+1}}\frac{|\sigma^{k+l}\cap\star v_{\tau(k)}|}{|\sigma^{k+l}|} = \sum_{\tau\in S_{k+l+1}}\frac{1}{k+l+1}=(k+l)!\, ,\]
which is a crucial observation for the original proof to apply to the alternative wedge product.

\section{Remeshing Cochains and Multigrid Extensions}\index{remeshing}\index{multigrid}
It is sometimes desirable, particularly in the context of
multigrid, multiscale, and multiresolution computations, to be able
to represent a discrete differential form which is given as a
cochain on a prescribed mesh, as one which is supported on a new
mesh. Given a differential form $\omega^k\in\Omega^k(K)$, and a
new mesh $M$ such that $|K|=|M|$, we can define it at the level of
cosimplices,
\[\forall\tau^k\in M^{(k)},\qquad \langle \omega^k, \tau^k \rangle =
\sum_{\sigma^k\in K^{(k)}} \operatorname{sgn}(\tau^k,\sigma^k)
\frac{|V_{\tau^k}\cap V_{\sigma^k}|}{|V_{\sigma^k}|}\langle
\omega^k, \sigma^k \rangle,\] and extend this by linearity to
cochains. Here, $ \operatorname{sgn}(\tau^k,\sigma^k)$ is $+1$ if
the orientation of $\tau^k$ and $\sigma^k$ are consistent, and
$-1$ otherwise. Since $k$-skeletons of meshes that are not related
by subdivision may not have nontrivial intersections,
intersections of support volumes are used in the remeshing
formula, as opposed to intersections of the $k$-simplices.

We denote this transformation at the level of cochains as,
$T_{K,M}:C^k(K)\rightarrow C^k(M)$. This has the natural property
that if we have a $k$-volume $U^k$ that can be represented as a
chain in either the complex $K$ or the complex $M$, that is to
say, $U^k=\sigma^k_1+\ldots+\sigma^k_l=\tau^k_1+\ldots+\tau^k_l$,
then we have
\begin{align*}
\langle \omega^k, \tau^k_1+\ldots+\tau^k_m \rangle &= \sum_{i=1}^m
\langle \omega, \tau^k_i \rangle = \sum_{i=1}^m \sum_{\sigma^k\in
K^{(k)}} \operatorname{sgn}(\tau^k_i,\sigma^k)
\frac{|V_{\tau^k_i}\cap V_{\sigma^k}|}{|V_{\sigma^k}|} \langle
\omega^k,
\sigma^k \rangle\\
&= \sum_{i=1}^m \sum_{j=1}^l
\operatorname{sgn}(\tau^k_i,\sigma^k_j) \frac{|V_{\tau^k_i}\cap
V_{\sigma^k_j}|}{|V_{\sigma^k_j}|} \langle \omega^k, \sigma^k_j
\rangle\\
&= \sum_{j=1}^l \sum_{i=1}^m
\operatorname{sgn}(\tau^k_i,\sigma^k_j) \frac{|V_{\tau^k_i}\cap
V_{\sigma^k_j}|}{|V_{\sigma^k_j}|} \langle \omega^k,
\sigma^k_j \rangle\\
&= \sum_{j=1}^l \langle \omega^k, \sigma^k_j \rangle = \langle
\omega^k, \sigma^k_1+\ldots+\sigma^k_l\rangle.
\end{align*}
Which is to say that the integral of the differential form over
$U^k$ is well-defined, and independent of the representation of
the differential form.

Note that, in particular, if we choose to coarsen the mesh, the
value the form takes on a cell in the coarser mesh is simply the
sum of the values the form takes on the old cells of the fine mesh
which make up the new cell in the coarser mesh.

\paragraph{Non-Flat Manifolds.}\index{manifolds!non-flat}
The case of non-flat manifolds presents a challenge in remeshing
akin to that encountered in the discretization of differential
forms. In particular, if the two meshes represent different
discretizations of a non-flat manifold, they will in general
correspond to different polyhedral regions in the embedding space,
and not have the same support region.

We assume that our discretization of the manifold is sufficiently
fine that for every simplex, all its vertices are contained in
some chart. Then, by using these local charts, we can identify
support volumes in the computational domain with $n$-volumes in
the manifold, and thereby make sense of the remeshing formula.

\section{Conclusions and Future Work}
We have presented a framework for discrete exterior calculus using the cochain representation of discrete differential forms, and introduced combinatorial representations of discrete analogues of differential operators on discrete forms and discrete vector fields. The role of primal and dual cell complexes in the theory are developed in detail. In addition, extensions to dynamic problems and multi-resolution computations are discussed.

In the next few paragraphs, we will describe some of the future directions that emanate from the current work on discrete exterior calculus.

\paragraph{Relation to Computational Algebraic Topology}
Since we have introduced a discrete Laplace-deRham operator, one can hope to develop a discrete Hodge-deRham theory, and relate the deRham cohomology of a simplicial complex to its simplicial cohomology.

\paragraph{Extensions to Non-Flat Manifolds.}
The intrinsic notion of what constitutes the discrete tangent space to a node on a non-flat mesh remains an open question. It is possible that this notion is related to a choice of discrete connection on the mesh, and it is an issue that deserves further exploration.

\paragraph{Generalization to Arbitrary Tensors.}
The discretization of differential forms as cochains is particularly natural, due to the pairing between forms and volumes by integration. When attempting to discretize an arbitrary tensor, the natural discrete analogue is unclear. In particular, while it is possible to expand an arbitrary tensor using the tensor product of covariant and contravariant one-tensors,  this would be cumbersome to represent on a mesh. In \cite{LeMaWe2003}, which is on discrete connections, we will see Lie group-valued discrete $1$-forms, and one possible method of discretizing a $(p,q)$-tensor that is alternating in the contravariant indices, is to consider it as a $(0,q)$-tensor-valued discrete $p$-form.

It would be particularly interesting to explore this in the context of the elasticity complex (see, for example,~\cite{Arnold2002}),
\[\xymatrix{\mathfrak{se}(3)\,\ar@{^{(}->}[r] & C ^\infty(\Omega,\mathbb{R}^3)\ar[r]^{\epsilon} & C^\infty(\Omega,\mathbb{S}) \ar[r]^{J} & C^\infty(\Omega,\mathbb{S})\ar[r]^{\operatorname{div}} & C^\infty(\Omega,\mathbb{R}^3)\ar[r] & 0\, ,
}\] where $\mathbb{S}$ is the space of $3\times 3$ symmetric
matrices. One approach to discretize this was suggested
in~\cite{Arnold2002}, which cites the use of the
Bernstein--Gelfand--Gelfand resolution in \cite{Eastwood2000} to
derive the elasticity complex from the deRham complex. Alternatively, it might be appropriate in the context of the elasticity complex to consider Lie algebra-valued discrete differential forms. 

\paragraph{Convergence and Higher-Order Theories.}
The natural question from the point of view of numerical analysis would be to carefully analyze the convergence properties of these discrete differential geometric operators. In addition, higher-order analogues of the discrete theory of exterior calculus are desirable from the point of view of computational efficiency, but the cochain representation is attractive due to its conceptual simplicity and the elegance of representing discrete operators as combinatorial operations on the mesh.

It would therefore be desirable to reconcile the two, by ensuring that high-order interpolation and combinatorial operations are consistent. As a low-order example, Whitney forms, which are used to interpolate differential forms on a simplicial mesh, have the nice property that taking the Whitney form associated with the coboundary of a simplicial cochain is equal to taking the exterior derivative of the Whitney form associated with the simplicial cochain. As such, the coboundary operation, which is a combinatorial operation akin to finite differences, is an exact discretization of the exterior derivative, when applied to the degrees of freedom associated to the finite-dimensional function space of Whitney forms.

It would be interesting to apply subdivision surface techniques to construct interpolatory spaces that are compatible with differential geometric operations that are combinatorial operations on the degrees of freedom. This will result in a massively simplified approach to higher-order theories of discrete exterior calculus, by avoiding the use of symbolic computation, which would otherwise be necessary to compute the action of continuous exterior differential operators on the polynomial expansions for
differential forms.

\bibliographystyle{plainnat}
\bibliography{umich_dec}

\end{document}